\documentclass[titlepage,final,11pt]{article}
\usepackage{theorem}
\usepackage{enumerate}
\usepackage{array}
\usepackage[reqno]{amsmath}
\usepackage{amssymb}
\usepackage{mathtools}
\usepackage{latexsym}
\usepackage{makeidx}
\usepackage{fancybox}
\input epsf.sty
\usepackage[dvips]{graphicx,color}
\usepackage{showlabels}
\usepackage{subcaption}
\usepackage{url}
\usepackage{xcolor}

\usepackage{epsfig}

\newcommand{\drawing}[4]{
	\begin{figure}[!hbt]
	\begin{center}
	\leavevmode
	\epsfxsize=#2
	\epsfbox{#1}
	\caption{\small #3}
	\label{#4}
	\end{center}
	\end{figure}}

\theoremstyle{change}
\newtheorem{proclaim}{PROCLAIM}[section]
\newtheorem{theorem}[proclaim]{Theorem}
\newtheorem{definition}[proclaim]{Definition}

\newtheorem{proposition}[proclaim]{Proposition}

\newtheorem{example}[proclaim]{Example}

\numberwithin{equation}{section}

\outer\def\proclaim #1. #2\par{\medbreak \noindent{\bf#1.\enspace}{\sl#2}\par
  \ifdim\lastskip<\medskipamount
  \removelastskip\penalty55\medskip\fi}
\def\state #1. { \noindent{\bf#1.\enspace}}
\def\algo #1. { \noindent{\bf#1.\enspace}}

\DeclareMathOperator{\con}{con}

\DeclareMathOperator{\dom}{dom}

\DeclareMathOperator{\nt}{int}

\DeclareMathOperator{\prj}{prj}

\DeclareMathOperator{\bprob}{b-prob}
\DeclareMathOperator{\srsk}{s-rsk}

\newcommand{\comp}{\,{\raise 1pt \hbox{$\scriptstyle\circ$}}\,}

\newcommand{\reals}{\mathbb{R}}

\newcommand{\natnums}{{{\rm l} \kern -.13em {\rm N} }}

\newcommand{\snats}{{I\kern -.29em N}}
\newcommand{\rats}{{Q\kern -.64em \raise 1pt \hbox{$\scriptstyle |$}\;\,}}
\newcommand{\srats}
	{{Q\kern -.56em \raise 1.2pt \hbox{$\scriptscriptstyle /$}\,}}
\newcommand{\ints}{Z\kern -.46em Z}

\newcommand{\pluss}{\hskip1pt \raise1pt\vbox{\hrule width6pt \vskip1pt \hrule
                    width6pt} \kern-4pt{\lower1pt\hbox{\vrule height6pt
		    \kern1pt\vrule height6pt}}\hskip5pt}
\newcommand{\eop}
	{\hfill{$\vcenter{\hrule height1pt \hbox{\vrule width1pt height5pt
   	 \kern5pt \vrule width1pt} \hrule height1pt}$} \medskip}
\newcommand{\half}
	{{\raisebox{1pt}{$\frac{1}{2}$}}}

\newcommand{\fourthird}
	{{\raisebox{1pt}{$\frac{4}{3}$}}}

\newcommand{\proxfacnu}
	{{\raisebox{1pt}{$\frac{1}{2\lambda^\nu}$}}}

\newcommand{\lplus}{{\scriptscriptstyle +}}

\newcommand{\setd}{{ d \kern -.15em l}}
\newcommand{\hatsetd}{ d \hat{\kern -.15em l }}

\renewcommand{\epsilon}{\varepsilon}
\renewcommand{\phi}{\varphi}

\newcommand{\tto}{\;{\lower 1pt \hbox{$\rightarrow$}}\kern -12pt
           \hbox{\raise 2.5pt \hbox{$\rightarrow$}}\;}
\newcommand{\overto}[1]{\,{\raise 0pt\hbox{$\rightarrow$}}\kern -9pt
     \hbox{\lower 3pt \hbox{$\scriptscriptstyle#1$}}\hskip6pt}
\newcommand{\underto}[1]{\,{\lower 1pt\hbox{$\rightarrow$}}\kern -9pt
     \hbox{\raise 4pt \hbox{$\,\scriptscriptstyle#1$}}\hskip7pt}
\newcommand{\bigoverto}[1]{{\raise 0pt\hbox{$\,\longrightarrow$}}\kern -16pt
     \hbox{\lower 3pt \hbox{$\scriptscriptstyle#1$}}\hskip4pt}
\newcommand{\bigunderto}[1]{\,{\lower 1pt\hbox{$\longrightarrow$}}\kern -16pt
     \hbox{\raise 4pt \hbox{$\,\scriptscriptstyle#1$}}\hskip6pt}
\newcommand{\bigbigto}[2]{\,{\raise 0pt\hbox{$\,\longrightarrow$}}\kern -16pt
     \hbox{\lower 3pt \hbox{$\scriptscriptstyle#2$}}\kern -10pt
     \hbox{\raise 4pt \hbox{$\,\scriptscriptstyle#1$}}\hskip7pt}
\newcommand{\downto}{{\raise 1pt \hbox{$\scriptscriptstyle \,\searrow\,$}}}
\newcommand{\upto}{{\raise 1pt \hbox{$\scriptscriptstyle \,\nearrow\,$}}}

\newcommand{\notimply}
	{\quad\hbox{$\Longrightarrow \kern -14pt {/}$}\hskip6pt\quad}

\newcommand{\lto}{\,{\lower 1pt\hbox{$\rightarrow$}}\kern -10pt
     \hbox{\raise 4pt \hbox{$\, \scriptstyle l$}}\hskip7pt}
\newcommand{\eto}{\,{\lower 1pt\hbox{$\rightarrow$}}\kern -11pt
     \hbox{\raise 4pt \hbox{$\, \scriptstyle e$}}\hskip7pt}
\newcommand{\hto}{\,{\lower 1pt\hbox{$\rightarrow$}}\kern -11pt
     \hbox{\raise 4pt \hbox{$\, \scriptstyle h$}}\hskip7pt}
\newcommand{\pto}{\,{\lower 1pt\hbox{$\rightarrow$}}\kern -11pt
     \hbox{\raise 4.5pt \hbox{$\, \scriptstyle p$}}\hskip7pt}
\newcommand{\cto}{\,{\lower 1pt\hbox{$\rightarrow$}}\kern -11pt
     \hbox{\raise 4pt \hbox{$\, \scriptstyle c$}}\hskip7pt}
\newcommand{\gto}{\,{\lower 1pt\hbox{$\rightarrow$}}\kern -11pt
     \hbox{\raise 4.5pt \hbox{$\, \scriptstyle g$}}\hskip7pt}
\newcommand{\sto}{\,{\lower 1pt\hbox{$\rightarrow$}}\kern -11pt
     \hbox{\raise 4pt \hbox{$\, \scriptstyle s$}}\hskip7pt}
\newcommand{\awto}{\,{\lower 1pt\hbox{$\rightarrow$}}\kern -15pt
     \hbox{\raise 4pt \hbox{$\, \scriptstyle aw$}}\hskip7pt}
\def\Nto{\,{\raise 1pt\hbox{$\rightarrow$}}\kern -13pt
     \hbox{\lower 3pt \hbox{$\, \scriptstyle N$}}\hskip7pt}
\def\Cto{\,{\raise 1pt\hbox{$\rightarrow$}}\kern -14pt
     \hbox{\lower 3pt \hbox{$\, \scriptstyle C$}}\hskip7pt}
\def\fto{\,{\raise 1pt\hbox{$\rightarrow$}}\kern -14pt
     \hbox{\lower 3pt \hbox{$\, \scriptstyle f$}}\hskip7pt}

\newcommand{\low}[1]{{\lower1pt \hbox{$\scriptstyle #1$}}}
\newcommand{\loww}[1]{{\lower2pt \hbox{$\scriptstyle #1$}}}
\newcommand{\high}[1]{{\raise1pt \hbox{$\scriptstyle #1$}}}

\newcommand{\cA}{{\cal A}}

\newcommand{\cC}{{\cal C}}
\newcommand{\cD}{{\cal D}}
\newcommand{\cE}{{\cal E}}
\newcommand{\cF}{{\cal F}}

\newcommand{\cL}{{\cal L}}

\newcommand{\cP}{{\cal P}}
\newcommand{\cQ}{{\cal Q}}
\newcommand{\cR}{{\cal R}}
\newcommand{\cS}{{\cal S}}

\newcommand{\cV}{{\cal V}}

\newcommand{\nsum}{\mathop{\sum}\nolimits}

\newcommand{\nsup}{\mathop{\rm sup}\nolimits}

\newcommand{\nnmin}{\mathop{\rm minimize}}

\newcommand{\nargmax}{\mathop{\rm argmax}\nolimits}
\newcommand{\nargmin}{\mathop{\rm argmin}\nolimits}

\newcommand{\var}{{\mathop{\rm var}\nolimits}}
\newcommand{\std}{\hspace{0.04cm}{\mathop{\rm std}\nolimits}\hspace{-0.04cm}}
\newcommand{\prob}{{\mathop{\rm prob}\nolimits}}

\newcommand{\bfxi}{\mbox{\boldmath $\xi$}}
\newcommand{\bfeta}{\mbox{\boldmath $\eta$}}
\newcommand{\bfzeta}{\mbox{\boldmath $\zeta$}}

\newcommand{\bfpi}{\mbox{\boldmath $\pi$}}

\newcommand{\bfz}{\mbox{\boldmath $z$}}

\newcommand{\bfnull}{\mbox{\boldmath $0$}}
\newcommand{\bfone}{\mbox{\boldmath $1$}}

\newcommand{\lwdy}[2]{\mathrel{\mathop
        {\raisebox{0.1ex}{\null$#1$}}{\hbox{\kern -1.0em
	{\raisebox{-0.8ex}{$\scriptstyle{\;\to #2}$}}}}}}
\newcommand{\lwwdy}[2]{\mathrel{\mathop
        {\raisebox{0.2ex}{\null$#1$}}{\hbox{\kern -1.0em
	{\raisebox{-1.1ex}{$\scriptstyle{\;\to #2}$}}}}}}
\newcommand{\slwdy}[2]{\scriptsize{{\mathrel{\mathop
        {\raisebox{0.1ex}{\null$#1$}}{\hbox{\kern -1.0em
	{\raisebox{-0.8ex}{$\scriptstyle{\;\to #2}$}}}}}}}}
\newcommand{\slwwdy}[2]{\scriptsize{{\mathrel{\mathop
        {\raisebox{0.2ex}{\null$#1$}}{\hbox{\kern -1.0em
	{\raisebox{-1.1ex}{$\scriptstyle{\;\to #2}$}}}}}}}}

\definecolor{lightgray}{gray}{0.75}
\definecolor{myred}{rgb}{0.55,0,0}
\definecolor{myblue}{rgb}{0,0,0.5} 
\definecolor{mygreen}{rgb}{0,0.5,0} 
\definecolor{purple}{rgb}{0.5,0,0.5} 
\definecolor{turq}{rgb}{0,0.805,0.816} 
\definecolor{maroon}{rgb}{0.51,0,0}
\definecolor{MAROON}{rgb}{0.51,0,0}
\definecolor{redor}{rgb}{0.78,0.078,0.078}
\definecolor{dgreen}{rgb}{0,0.3,0}

\newcommand{\Ex}{\mathbb{E}}

\newcommand{\bcdot}{\,{\raise .2ex \hbox{$\centerdot$}}\,}

\newcommand{\bbI}{\mathbb{I}}

\newcommand{\bbP}{\mathbb{P}}

\begin{document}


\begin{center}
\begin{large}
{\bf Risk-Adaptive Approaches to Stochastic Optimization: A Survey}
\smallskip
\end{large}
\vglue 0.2truecm
\begin{tabular}{cc}
  \begin{large} {\sl Johannes O. Royset} \end{large}\\
  Daniel J. Epstein Department of Industrial \& Systems Engineering\\
  University of Southern California\\
  royset@usc.edu
\end{tabular}

\vskip 0.2truecm

\end{center}

\noindent {\bf Abstract}. \quad Uncertainty is prevalent in engineering design, data-driven problems, and decision making broadly. Due to inherent risk-averseness and ambiguity about assumptions, it is common to address uncertainty by formulating and solving conservative optimization models expressed using measures of risk and related concepts. We survey the rapid development of risk measures over the last quarter century. From their beginning in financial engineering, we recount the spread to nearly all areas of engineering and applied mathematics. Solidly rooted in convex analysis, risk measures furnish a general framework for handling uncertainty with significant computational and theoretical advantages. We describe the key facts, list several concrete algorithms, and provide an extensive list of references for further reading. The survey recalls connections with utility theory and distributionally robust optimization, points to emerging applications areas such as fair machine learning, and defines measures of reliability.\\

\vskip 0.1truecm



\noindent {\bf Date}:\quad \today

\baselineskip=15pt

\tableofcontents

\section{Introduction}\label{sec:intro}

In many areas of applied mathematics, we aim to determine a best decision, design, estimate, or model in the presence of uncertainty about the values of key parameters and data in the problem. The uncertainty stems from our incomplete knowledge about the state of a system, the inaccuracy of approximations that might have been adopted, and the human inaptitude to predict the future. One might hope that historical observations about varying quantities can inform us and ideally eliminate the uncertainty. Despite living in the era of big data, this aspiration rarely comes to fruition. Real-world data sets tend to be noisy, biased, corrupted, or simply insufficiently large. It is therefore prudent to optimize decisions, designs, and models while accounting {\em conservatively} for the uncertainty. Even when our grasp of the uncertainty in an application is fairly good, human nature and our desire to avoid exceptionally ``bad'' outcomes motivate conservativeness in decision making and modeling broadly. Thus, financial decisions, operational plans, military strategies, system controls, engineering designs, and statistical models are often selected conservatively.

While the various fields have approached decision making under uncertainty somewhat differently, overarching themes now emerge: the need to capture risk-averseness and ambiguity about uncertainty, promote computations and analysis, and enable explanations of algorithmic outcomes. The concept of {\em risk measures} provides a mathematical framework for unifying and understanding the various threads and how they connect. The framework encapsulates many common problems in robust control and optimization, distributionally robust optimization, adversarial machine learning, statistical estimation, financial risk management, utility maximization, attacker-defender games, and reliability-based design optimization. Supported by convex analysis, risk measures furnish a rich area for nonlinear analysis as well as opportunities for efficient computations. The vast literature on risk measures developed over the last 25 years is a testimony to the potency of the framework, both theoretically and practically.

In this survey, we discuss risk measures and related concepts with a focus on advances over the last quarter century. The canonical problem involves a {\em quantity of interest} given by a function $f(\xi,x)$, which depends on a {\em parameter} $\xi$ and a {\em decision (control)} $x$. For example, $f(\xi,x)$ might quantify the performance of an engineering system designed according to our decision $x$, given an environmental condition represented by $\xi$. In supervised learning, $f(\xi,x)$ might specify the prediction error of a neural network designed according to $x$, given feature and label data $\xi$. Without fully knowing $\xi$, the problem is to determine a decision $x$ such that $f(\xi,x)$ is minimized or, alternatively, $f(\xi,x)$ does not exceed a given threshold. The problem is ill-posed because $\xi$ is unsettled. The uncertainty about $\xi$ could stem from our incomplete knowledge of a ``true'' value, such as the demand for a product tomorrow, and/or from inherent variability in the value as exemplified by the many feature-label pairs a neural network needs to handle accurately.

There are several ways to proceed. One can estimate the value of $\xi$ and adopt that value in the subsequent optimization and decision making. However, this fails to account for the uncertainty associated with $\xi$. The approach is especially problematic when $\xi$ varies inherently and our decision $x$ needs to perform satisfactorily under different values of $\xi$. Alternatively, if there is a set $\Xi$ of possible values of $\xi$, then one may consider the quantity of interest in the worst case across these values, i.e.,  $\sup_{\xi\in \Xi} f(\xi,x)$. The problem shifts to determining a decision $x$ that minimizes this conservative quantity or makes it sufficiently low, which we refer to as {\em robust optimization}; see, e.g., \cite{BentalElghaouiNemirovski.09,BertsimasBrownCaramanis.11}. Yet another possibility is to model the uncertainty associated with $\xi$ using a probability distribution. This brings in the vast and sophisticated tools of probability theory and statistics. In assessing a decision $x$, one can leverage the expected value of $f(\,\cdot\,,x)$ computed with respect to the adopted probability distribution. Thus, the problem becomes to find a decision that is satisfactory on average. This classical approach is the central tenet of {\em stochastic programming}; see, e.g., \cite{BirgeLouveaux.11,ShapiroDentchevaRuszczynski.21,primer}.

Risk measures capture all these possibilities and many more. Through choices of probability distributions, risk measures allow us to incorporate data and other information about the possible values of $\xi$ and their likelihoods. They reflect a wide variety of concerns and preferences a decision maker may have toward various outcomes.  When restricted suitably to the class of {\em regular measures of risk}, they exhibit theoretically and computationally advantageous properties. In particular, convexity, linearity, and smoothness of $f(\xi,x)$ in $x$ commonly carry over to the resulting optimization problem. Many risk measures also address the fact that any adopted probability distribution is an imperfect model of the uncertainty associated with $\xi$. Thus, they account for ambiguity about the model of uncertainty, which leads to {\em distributionally robust optimization problems}, as exemplified by \cite{WiesemannKuhnSim.14,RoysetWets.16b,MohajerinesfahaniKuhn.18,ShapiroDentchevaRuszczynski.21}, and connections with {\em stochastic dominance}; see, e.g., \cite{DentchevaRuszczynski.03}.

Risk measures originated in financial engineering as an approach to quantify the reserves banks, insurance companies, and other financial institutions need to cover potential future losses; see for example the influential paper \cite{ArtznerDelbaenEberHeath.99}, the textbooks \cite{McNeilFreyEmbrechts.15,FollmerSchied.16}, and the review article \cite{FollmerWeber.15}. In this survey, we discuss the spread of risk measures {\em beyond} financial engineering to operations management, reliability analysis, engineering design, defense planning, statistics, and machine learning.

Among the early reviews of risk measures, \cite{Rockafellar.07} stands out for a concise description of the important advances taking place around the turn of the century. The tutorials \cite{SarykalinSerrainoUryasev.08,KrokhmalZabarankinUryasev.11} explain key concepts and introduce connections with statistics; see also the monograph \cite{ZabarankinUryasev.14}. With a focus on expected utility theory and dual utility theory, \cite{Ruszczynski.13} reviews the deep mathematical relations between risk measures and these earlier theories as well as stochastic dominance; see also the recent paper \cite{FrohlichWilliamson.22}.
While surveying risk measures and related concepts, \cite{RockafellarUryasev.13} breaks new ground by connecting concepts of risk, regret, deviation, and error and thereby relating risk management to statistics.
The paper \cite{RockafellarRoyset.15b} emphasizes the applicability of risk measures in reliability-based engineering design. For PDE-constrained optimization and uncertainty quantification, \cite{KouriShapiro.18} discusses the needed technical assumptions in such infinite-dimensional settings as well as algorithms. The paper \cite{WangChapman.22} surveys applications for autonomous systems.
Reviews of superquantile risk measures appear in \cite{FilippiGuastarobaSperanza.20}, with a focus on applications in supply chain management, scheduling, networks, energy, and medicine, and in \cite{LaguelPillutlaMalickHarchaoui.21}, dealing with machine learning applications. Books covering risk measures broadly include \cite{PflugRomisch.07,FollmerSchied.16,ShapiroDentchevaRuszczynski.21,primer}.

This survey provides a succinct introduction to the area of risk measures and related concepts without devolving into the more technical aspects. We discuss the increasing interest in risk-averse approaches to statistical applications, with an updated review of the risk quadrangle proposed in \cite{RockafellarUryasev.13} and refined in \cite{RockafellarRoyset.15}. The survey recounts the historical development of superquantiles (a.k.a. conditional value-at-risk, average value-at-risk, tail value-at-risk, and expected shortfall) and the central role they now play in many areas of operations research, engineering, and statistics. This includes a behind-the-scene description of the derivations that led to an influential formula for superquantiles. With minimal overhead from probability theory, we describe a duality theory and connections with distributionally robust optimization. We introduce the terminology {\em measure of reliability} for failure probabilities, buffered failure probabilities, buffered probabilities of exceedance, and related concepts. Throughout, the survey focuses on computational aspects and implementable algorithms.

To avoid technical distractions, we mainly focus on finite-dimensional decision and uncertainty spaces, i.e. $x$ and $\xi$ are finite-dimensional vectors. We limit the scope to ``here-and-now'' decisions, but include a brief summary of multi-stage decision processes in Section \ref{sec:extensions}.

The survey begins in Section \ref{sec:decision} with examples of decision making under uncertainty and the introduction of risk measures. Section \ref{sec:super} discusses superquantiles and algorithms for their optimization and estimation. Section \ref{sec:regretrisk} reviews measures of risk broadly and their companions: measures of regret.  Section \ref{sec:errordev} shifts the focus to the statistical domain by discussing measures of error and deviation, and their application in generalized regression modeling. With their roots in convex analysis, many measures of risk are equivalently expressed by ``dual formulas'' that furnish additional computational possibilities as well as deep connections with distributionally robust optimization; see Section \ref{sec:duality}. Section \ref{sec:addexamples} provides additional examples and Section \ref{sec:comptool} reviews computational tools for risk-adaptive optimization. The survey ends with extensions and open questions in Section \ref{sec:extensions}.

\section{Decision Making under Uncertainty}\label{sec:decision}

This section starts by illustrating how problems of decision making under uncertainty arise in diverse areas. We define measures of risk, give concrete examples, and make connections with utility theory.

\subsection{Engineering Design and Control}

The future responses, levels of damage, rates of deterioration, life-cycle costs, and numerous other quantities of interest for an engineering system are invariably uncertain. This is caused by unknown environmental conditions, load patterns, as well as our imperfect modeling of underlying mechanisms. Thus, design and control of engineering systems give rise to a wealth of decision making problems under uncertainty. We discuss three concrete examples, starting with a system described by explicit formulas for concreteness and ending with a high-level description of a multi-disciplinary design problem.

\drawing{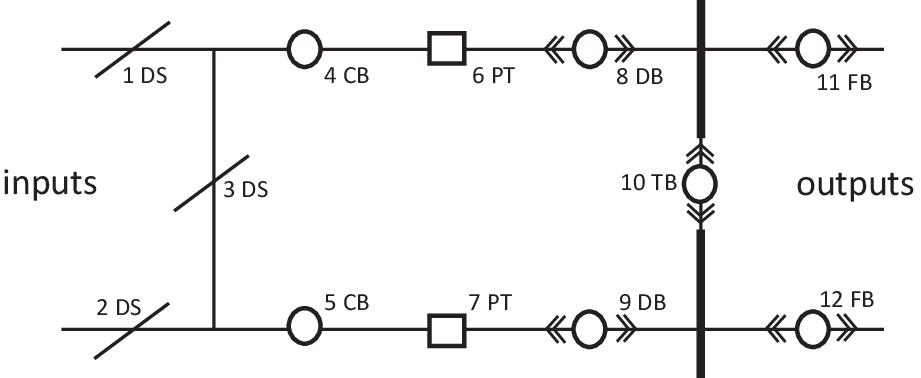}{3.4in} {Outline of electrical substation (reproduced from \cite{ByunDeoliveiraRoyset.22}).}{fig:substation}

\begin{example}{\rm (electrical substation design)}. Figure \ref{fig:substation} describes a planned two-transmission-line electrical substation  with 12 components of various kinds; disconnect switches (DS), circuit breakers (CB), power transformers (PT), drawout breakers (DB), tie breakers (TB), and feeder breakers (FB). In a reliability growth model, the failure rate of component $i$ is defined as $\lambda_i(x_{i}) = \alpha \beta e^{-\beta x_{i}}$, where $\alpha=9$, $\beta=2$, if it undergoes $x_i$ days of testing prior to installation; see \cite{ByunDeoliveiraRoyset.22}. Under decision $x_i$ and outcome $\xi_i$ of a uniformly distributed random variable on $[0,1]$, we assess the performance of component $i$ by the number of days it falls short of the desired lifetime of $\tau=365$, which becomes
\[
f_i(\xi_i, x_i) = \tau + \frac{\ln \xi_i}{\lambda_i (x_{i})}
\]
under an exponential failure distribution. We assess the overall system performance by identifying which component failures would trigger a disconnect between an input and an output in Figure \ref{fig:substation}. As identified in \cite{ByunDeoliveiraRoyset.22}, there are 25 such failure modes producing the quantity of interest
\[
f(\xi,x) = \max_{k=1, \dots, 25} \min_{i \in \bbI_k} f_i(\xi_i,x_i),
\]
where $\xi = (\xi_1, \dots, \xi_{12})$, $x = (x_1, \dots, x_{12})$, and the 25 sets $\bbI_k$ are  $\{1, 2\}$, $\{4, 5\}$, $\{4, 7\}$, $\{4, 9\}$, $\{5, 6\}$, $\{6, 7\}$, $\{6, 9\}$, $\{5, 8\}$, $\{7, 8\}$, $\{8, 9\}$, $\{11, 12\}$, $\{1, 3, 5\}$, $\{1, 3, 7\}$, $\{1, 3, 9\}$, $\{2, 3, 4\}$, $\{2, 3, 6\}$, $\{2, 3, 8\}$, $\{4, 10, 12\}$, $\{6, 10, 12\}$, $\{8, 10, 12\}$, $\{5, 10, 11\}$, $\{7, 10, 11\}$, $\{9, 10, 11\}$, $\{1, 3, 10, 12\}$, $\{2, 3, 10, 11\}$. Thus, $f(\xi,x)$ is the number of days the overall system falls short of the desired lifetime of $\tau=365$ under testing decision $x$ and outcome $\xi$ of the inherent randomness in the components' lifetimes. An engineer would like to make a constrained decision about $x$ such that $f(\xi,x)$ is sufficiently low, ideally below $0$, but needs to account for the uncertainty associated with $\xi$.
\end{example}

\begin{example}{\rm (additive manufacturing)}.
Advances in additive manufacturing enable optimization of novel structures tailored to specific loads and boundary conditions. For example, as described in \cite{Kouri.21}, we may seek to design a three-dimensional object represented by a density $z:D\to [0,1]$, where $D\subset\reals^3$ is the physical domain. Under the assumption of linear elastic material behavior, the displacement $u:D\to \reals^3$ of an object with density $z$ is determined by the weak form of the equations
\begin{align*}
  -\mbox{div}(E(z) : \epsilon) & = F  \mbox{ in } D\\
  \epsilon & = \tfrac{1}{2}(\nabla u + \nabla u^\top)  \mbox{ in } D\\
  \epsilon n & = t \mbox{ on } \Gamma_n; ~~~ u = g \mbox{ on } \Gamma_d,
\end{align*}
where $E(z)$ is the stiffness tensor (certainly dependent on $z$), $\epsilon$ is the strain tensor, $F$ is an external load, $t$ is a boundary condition on the Neumann boundary $\Gamma_n$, $n$ is the outward normal of $D$, and $g$ is a boundary condition on the Dirichlet boundary $\Gamma_d$. Typically, the material properties $E(z)$ is uncertain, even for a fixed density $z$. The external load $F$ and the boundary conditions specified by $t$ and $g$ are also uncertain. Thus, a solution $u$ of the equations actually depends on the values of these uncertain parameters as well as the design decision $z$. A quantity of interest to be minimized could be the compliance expressed by $\int_D F \cdot u + \int_{\Gamma_n} t \cdot u$, which depends directly on the uncertain $F$ and $t$ as well as indirectly on the choice of density $z$ and the uncertain parameters through the solution $u$.
\end{example}

\begin{example}{\rm (supercavitating hydrofoil)}.\label{eHydrofoil}
Stretching back to the pioneering work by Enrico Forlanini in 1906, naval architects have aspired to essentially lift vessels out of the water using hydrofoils and thus dramatically reduce drag and increase speed. A major challenge, however, is the inception of cavitation, i.e., water vaporization caused by the pressure on the suction side of a hydrofoil dropping below the surface tension. This causes flow unsteadiness and trigger erosion with resulting loss of material integrity. Recent efforts aim to overcome these challenges by optimizing the shape of the hydrofoil so that its hydrodynamical and structural performance is satisfactory even in the presence of uncertainty about material properties and operating conditions; see \cite{BonfiglioRoyset.19} and references therein. This involves simulating the performance through numerical solution of unsteady Reynolds-averaged Navier-Stokes equations and linear elasticity equations. For instance, \cite{BonfiglioRoyset.19} considers five uncertain parameters $\xi = (\xi_1, \dots, \xi_5)$: cavitation index $\xi_1$ (which is a surrogate for speed) and material properties represented by Young's modulus $\xi_2$, Poisson's ratio $\xi_3$, density $\xi_4$, and yield stress $\xi_5$. A vector $x\in \reals^{17}$ specifies the shape of the hydrofoil. The quantities of interest are lift force, drag-to-lift ratio, von Mises stress, and displacement at the hydrofoil tip, and thus give rise to four functions $f_1(\xi,x), \dots, f_4(\xi,x)$.
\end{example}

We refer to \cite{GeiheLenzRumpfSchultz.13,KolvenbachLassUlbrich.18} for further applications in the area of shape optimization, especially of electrical engines and elastic structures. In design of a heterogeneous lens under manufacturing uncertainty, \cite{AlghamdiChenKaramehmedovic.22} optimizes photonic nanojets. The paper \cite{AntilKouriPfefferer.21} studies optimal control of elliptic PDEs with uncertain fractional exponents. A discussion of two-dimensional advection–diffusion equations and one-dimensional Helmholtz equations under uncertainty appears in \cite{ZouKouriAquino.19}. Welded beam structures and truss bridge structures are the subjects of \cite{ByunRoyset.22}. 
The paper \cite{ChenHabermanGhattas.21} highlights the large size of many engineering problems. In the design of an acoustic cloak under uncertainty, the resulting optimization problem involves one million design variables and half a million uncertain parameters. Thus, it becomes paramount to avoid adding significant complexity while modeling uncertainty as the problem with (assumed) known parameters is already challenging.

\subsection{Stochastic Optimization in Statistics and Machine Learning}\label{subsec:stat}

Some problems in statistics and machine learning can also be viewed as attempting to make a decision under uncertainty. In supervised learning, we aim to find a statistical model that best predicts an unknown output based on a given input. Typically, a statistical model is specified by a vector $c$ of coefficients, with the resulting prediction being $g(\xi; c)$ for input $\xi$.
If the input $\xi = (\xi_1, \dots, \xi_n) \in\reals^n$ and the statistical model is affine specified by the coefficients $c= (c_0, c_1, \dots, c_n)\in\reals^{1+n}$, then $g(\xi;c) = c_0 + \nsum_{i=1}^n c_i \xi_i$. 
But, $g(\xi;c)$ may take many other forms, possibly utilizing neural networks. Given an input-output pair $(\xi,\eta)$, we would like the model's prediction $g(\xi;c)$ to be close in some sense to the actual output value $\eta$. In regression analysis, the output quantity is a scalar and ``close'' is commonly quantified using 
\[
f\big((\xi,\eta), c\big) = \big(\eta - g(\xi;c) \big)^2~~~  \mbox{ or } ~~~f\big((\xi,\eta), c\big) = \big|\eta - g(\xi;c) \big|.
\]
If the sign of the prediction error is important, then the quantity of interest might be $f((\xi,\eta), c) = \eta - g(\xi;c)$, which distinguishes between overestimating and underestimating \cite{RockafellarUryasevZabarankin.08,RockafellarUryasev.13}.

In a $k$-class classification problem, the output $\eta \in \{1, \dots, k\}$ specifies a class and the statistical model $g(\xi;c) = (g_1(\xi;c), \dots, g_k(\xi;c))$ produces a $k$-dimensional vector of probabilities representing the likelihood that input $\xi$ corresponds to the various classes. In predicting $\eta$ from $\xi$, the cross-entropy of $g(\xi;c)$ relative to a probability mass function concentrated at $\eta\in \{1, \dots, k\}$  becomes the quantity of interest:
\[
f\big((\xi,\eta),c\big) = -\ln g_\eta(\xi;c).
\]
Regardless of the specific details, when selecting a statistical model we are uncertain about the input-output $(\xi,\eta)$ for which it should be accurate. In fact, we probably would like the statistical model to make accurate predictions for many input-output pairs. Consequently, we are faced with the problem of selecting $c$, under uncertainty about $(\xi,\eta)$, such that a quantity of interest $f((\xi,\eta),c)$ is ``optimized.''

\begin{example}{\rm (support vector machine with fairness constraint)}.
In binary classification, we seek to predict an output $\eta\in \{-1,1\}$ from an input $\xi$ using a statistical model $g(\xi;c)$, where $c$ is a vector of coefficients. Support vector machines achieve this by considering the hinge-loss
\[
f_1\big((\xi,\eta),c\big) = \max\big\{0, 1 - \eta g(\xi;c)\big\}
\]
as a quantity of interest. A concern in this setting is fairness; see, e.g., \cite{HashimotoSrivastavaNamkoongLiang.18,WilliamsonMenon.19}. When making mistakes, does the statistical model $g(\,\cdot\,;c)$ exhibit a bias in certain settings such as when being applied with an input $\xi$ corresponding to an under-represented group? To reduce bias in the statistical model, we may consider a secondary quantity of interest
\[
f_2\big((\xi,\zeta),c\big) = (\zeta - \bar \zeta) g(\xi;c),
\]
where $\zeta$ is an additional input representing sensitive attributes and $\bar \zeta$ is the average across a population of attributes; see, e.g., \cite{zafar2017fairness}. The secondary quantity of interest can be used to quantify the covariance between attributes and predictions produced by the statistical model and thus serves as a metric of fairness. Consequently, we are faced with a problem of selecting $c$, under uncertainty about $(\xi,\eta,\zeta)$, such that two quantities of interest are satisfactory.
\end{example}

\begin{example}{\rm (surrogate building and digital twins)}.\label{eSurrogate}
A computer model of a physical system, called a digital twin \cite{KapteynKnezevicWillcox.20,Thelenetal.22a,Thelenetal.22b}, as well as a coarse computer model of a high-fidelity simulation  \cite{PerdikarisVenturiRoysetKarniadakis.15,PeherstorferWillcoxGunzburger.18} require calibration to achieve acceptable predictions. These models lead to surrogates that predict an output $\eta$ from an input $\xi$. Given an input-output pair $(\xi,\eta)$, suppose that $g(\xi;c)$ is the prediction by a surrogate with coefficients $c$. Thus, a quantity of interest is the prediction error $f((\xi,\eta),c) = \eta - g(\xi;c)$.
Figure \ref{fig:lift} shows the results of low- and high-fidelity estimates of lift force (blue asterisks)  produced by a hydrofoil as described in Example \ref{eHydrofoil}. We see that the low-fidelity estimates do not match the high-fidelity estimates, but this can be addressed through calibration. With $\xi$ being a low-fidelity estimate, we may adopt $g(\xi;c) = c_0 + c_1 \xi$ as the prediction of the corresponding high-fidelity estimate. The tuning of the coefficients $c = (c_0, c_1)$ aims to bring the quantity of interest $f((\xi,\eta),c) = \eta - (c_0 + c_1 \xi)$ near zero, while accounting for the uncertainty (as illustrated by the spread in Figure \ref{fig:lift}) associated with $(\xi,\eta)$. Figure \ref{fig:lift} shows five such surrogates fitted with varying emphasis on conservativeness; dotted lines overestimate more than the dashed lines. Surrogates are often more refined than these affine models (see, e.g., \cite{VianaGoguGoel.21}), with neural networks emerging as key forms  \cite{NewmanRuthottoHartVanbloemenwaanders.21,WillardJiaXuSteinbachKumar.22,KarniadakisKevrekidisLuPerdikarisWangYang.22,CuomoDicolaGiampaoloRozzaRaissiPiccialli.22}. Regardless of the sophistication, the fundamental problem remains how to select a surrogate in the presence of uncertainty.
\end{example}

\drawing{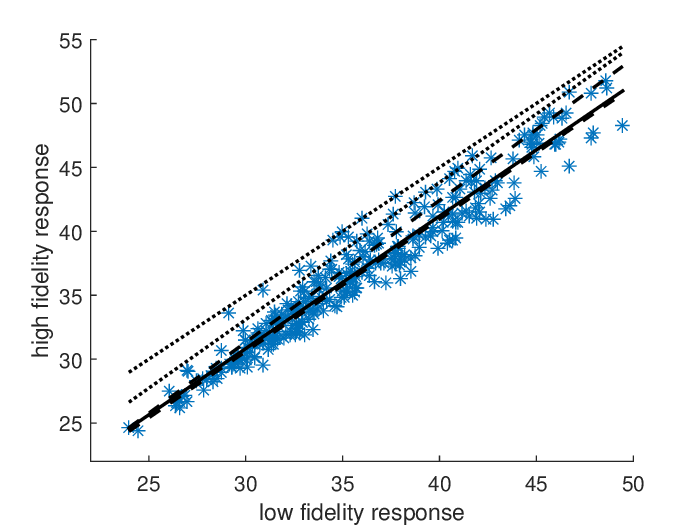}{3.7in} {Low- and high-fidelity estimates of lift force (blue asterisks) produced by a hydrofoil from Example \ref{eHydrofoil} \cite{BonfiglioRoyset.19} and five lines representing surrogates with varying degree of conservativeness.}{fig:lift}

\subsection{Other Application Areas}

Nearly all areas of human activity involve decision making under uncertainty. Financial engineering gives rise to many challenging problems in portfolio management and cashflow control; see, e.g., \cite{McNeilFreyEmbrechts.15,FollmerSchied.16}. The design of service networks takes place under uncertainty about future demand \cite{JiangBaiWallaceKendallLandasilva.21}. Demand is likewise uncertain in online and brick-and-mortar retail operations \cite{QiMakShen.20}. Transportation engineers account for uncertain travel times and travel patterns; see, e.g.,  \cite{CampbellGendreauThomas.11,ZhangWallaceGuoDongKaut.21,CominettiTorrico.22}. Planning of humanitarian relief has to be carried out without fully knowing where and how disasters may strike \cite{NoyanMerakliKucukyavuz.22}. The energy sector faces uncertainty in rainfall, which affects hydropower generation \cite{SeguinFletenCotePichlerAudet.17,AlaisCarpentierDelara.17}, but also uncertainty at the building level \cite{OttesenTomasgard.15} and in microgrids \cite{BagheriDagdouguiGendreau.22}. Decisions about when to charge and discharge an electrical storage unit face uncertainty about future supply and demand \cite{GangammanavarSen.16}.
Military jamming missions \cite{CommanderPardalosRyabchenkoUryasevZrazhevsky.07} and capital investments in defense technology \cite{TeterRoysetNewman.19} also take place under uncertainty. We refer to \cite{FilippiGuastarobaSperanza.20,LuShen.21} for a summary of applications in operations management broadly.
In this area, inventory control (see, e.g., \cite{AhmedCakmakShapiro.07,HanasusantoKuhnWallaceZymler.15}) gives rise to many challenging problems. A classical example illustrates the setting.

\begin{example}
\label{mNwsVnd}{\rm (newsvendor)}. In a contractual agreement, a newsvendor (a firm) places an order for the daily delivery of a fixed number of newspapers (perishable items) to meet a daily demand $\xi$. The newsvendor is charged $\gamma>0$ cents per paper ordered and sells each for $\delta>\gamma$ cents; unsold papers cannot be returned and are worthless at the end of the day. When ordering $x$ newspapers and $\xi$ is the demand, the loss (expense minus income) turns out to be
\[
f(\xi,x) =
\begin{cases}
   \gamma x -\delta x  &  \mbox{ if } \xi \geq x\\
   \gamma x - \delta \xi & \mbox{ otherwise,}
\end{cases}
\]
with negative values implying a profit. The loss becomes the quantity of interest that we seek to minimize while accounting for the uncertainty in $\xi$.
\end{example}

\subsection{Measures of Risk}\label{subsec:measureofrisk}

As exemplified in the previous subsections, we are often faced with a quantity of interest $f(\xi,x)$ and the need to manipulate it through the choice of $x$. However, $\xi$ is unsettled and we proceed by modeling its possible values and the associated likelihoods using a probability distribution. Just as $f(\xi,x)$ typically would be an imperfect representation of reality, the probability distribution is a {\em model} of the uncertainty associated with $\xi$. In a particular case, we construct the probability distribution using available observations of past values of $\xi$, expert opinions, and engineering judgement; see, e.g., \cite[Section 3.D]{primer}. The probability distribution might assign each outcome in a set the same likelihood (i.e., a uniform distribution) and even concentrate at a single point, which implies certainty about the value of $\xi$. Regardless of the situation, we indicate the shift from a parameter $\xi$ to a random quantity\footnote{As usual, a random quantity $\bfxi$ is a measurable mapping from some underlying probability space $(\Omega, \cA, \bbP)$ to a space of outcomes such as $\reals$, which makes $\bfxi$ a random variable, or to $\reals^m$, which makes $\bfxi$ a random vector.}  with a probability distribution by using boldface $\bfxi$. Outcomes of this random quantity are denoted by $\xi$.

For example, suppose that a quantity of interest is given by $f(\xi,x) = (\xi - 2/3)x - 1/3$, where $x \in \{-1, 1\}$, i.e., there are only two possible decisions. We model the uncertainty associated with $\xi$ using the triangular distribution on $[0,2]$ with mode at $0$. Thus, we are faced with the choice between the two random variables $f(\bfxi,1)$ and $f(\bfxi,-1)$; the first one quantifies the randomness under the decision $x = 1$ and the second one quantifies similarly  the randomness under $x = -1$. Since $\bfxi$ has the assumed triangular distribution, we can compute the probability density functions of $f(\bfxi,1)$ and $f(\bfxi,-1)$, denoted by $p_1$ and $p_2$, respectively; see Figure \ref{fig:tri2}.

Throughout the article, we assume that a quantity of interest represents cost, loss, damage, displacement in excess of a threshold, or another performance metric that we wish to be {\em low}. In this setting, what decision is better, $x=1$ or $x = -1$? One possible way of answering this question is to consider the expectations $\Ex[f(\bfxi,1)]$ and $\Ex[f(\bfxi,-1)]$, but they both come out as $-1/3$. Figure \ref{fig:tri2} highlights a significant difference between the decisions, however: $f(\bfxi,1)$ has values as high as $1$, while $f(\bfxi,-1)$ does not exceed $1/3$. Thus, from a worst-case point of view, $x = -1$ appears better.

\drawing{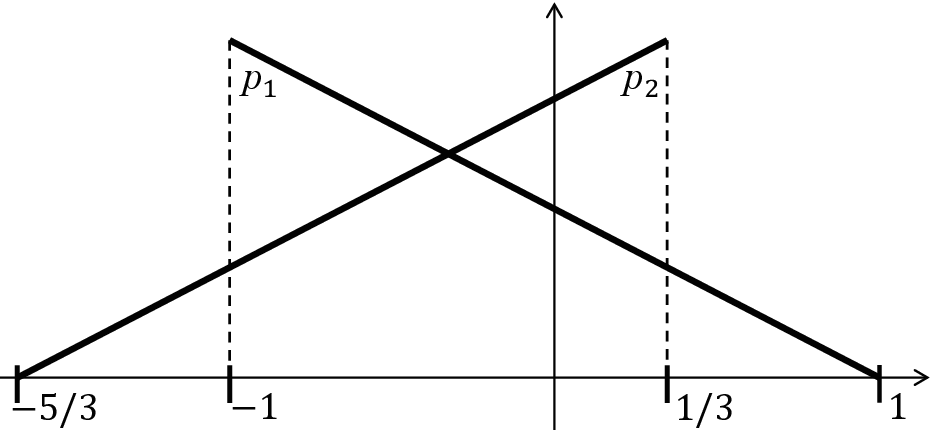}{3.4in} {Probability density functions $p_1$ and $p_2$ for random variables $f(\bfxi,1)$ and $f(\bfxi,-1)$, respectively.}{fig:tri2}


Expectations as well as worst-case outcomes can be viewed as defining measures of risk. A risk measure converts a random variable, say $f(\bfxi,1)$, into a scalar, which furnishes a basis for comparison with the scalar obtained from, say $f(\bfxi,-1)$. Thus, risk measures facilitate decision making when the quantity of interest is uncertain.

A risk measure assesses a random variable, which might have a complicated definition in terms of other random variables and decision variables. This is certainly the case when the random variable  represents a real-world quantity of interest as we have seen above. Even in the simple case of Figure \ref{fig:tri2}, risk measures assessed the random variables $f(\bfxi,1)$ and $f(\bfxi,-1)$, which in turn are defined by the triangularly distributed random variable $\bfxi$. Regardless of the circumstances, it is sometimes beneficial to ``hide'' this complexity and we often write $\bfxi$ for a generic random variable that we seek to assess using a risk measure.

\begin{definition}\label{dRiskmeasure}{\rm (risk measure)}.
A {\it measure of risk} (also called risk measure) $\cR$ assigns to a random variable $\bfxi$ a number $\cR(\bfxi) \in [-\infty,\infty]$ as a quantification of its risk, with this number being $\xi$ if $\prob\{\bfxi = \xi\} = 1$.
\end{definition}

Since we consistently prefer lower values to higher values, the choice between two random variables  $\bfxi_1$ and $\bfxi_2$ now reduces to comparing the two numbers $\cR(\bfxi_1)$ and $\cR(\bfxi_2)$ and selecting the random variable with the lower number.

The requirement $\cR(\bfxi) = \xi$ for a random variable $\bfxi$ concentrated at $\xi$ stems from the translation equivariance condition of \cite{ArtznerDelbaenEberHeath.99} and it is brought out explicitly in \cite{Rockafellar.07}. The requirement is certainly reasonable: a quantity of interest that has a specific value with certainty should indeed be assessed as that value. Subsection \ref{subsec:desirable} associates additional properties with risk measures. However, we prefer Definition \ref{dRiskmeasure} as a means to introduce the concept without imposing potentially limiting requirements.

Since a risk measure simply converts a random variable into a scalar, with a minor requirement for constant random variables, a large number of possibilities exist.\\

\noindent{\bf Expectation.} The choice $\cR(\bfxi) = \Ex[\bfxi]$, the expected value, defines a risk measure, but one that is insensitive to the possibility of high values as long as they are offset by low values. Thus, it is referred to as {\it risk-neutral}. This risk measure is meaningful in settings where a decision is implemented repeatedly and the focus is on good average performance, with little concern about fluctuations in performance. Moreover, one should be quite confident that the assumed probability distributions are accurate, at least in the sense that they produce the correct expected values. Stochastic programming originated with a focus on expectation minimization; see, e.g., \cite{Wets.74} for early developments. In statistics, expectation minimization arises in classification, regression, and other contexts.\\

Many real-life situations involve a decision that will only be executed once or at most a few times. Then, the average performance across many (hypothetical) instances is less important than guarantees about reasonable performance in the actual instances that will take place. This is often the case for major events such as an earthquake, a mars landing, or certain military missions. Even if a decision will be put to test often, we may not accept poor performance even though it is ``countered'' by excellent performance other times. This is common for engineering systems such as driverless cars where the performance needs to be sufficiently good nearly all the time.

We may also suspect that our models, producing the quantity of interest and probability distributions representing the uncertainty, are incomplete. The performance that we predict, even under assumed known conditions, could be erroneous. Thus, we may consider ``worst cases'' as a means to compensate for modeling deficiencies. These situations trigger the need for {\em risk-averse} measures of risk, which is the main focus of this article. They allow us to address decision makers that are willing to forego good average performance to achieve reduced possibility and severity of poor outcomes.

Generally, for a random variable $\bfxi$ with support\footnote{The support of $\bfxi$ is the smallest closed subset of $\reals$ containing outcomes of $\bfxi$ occurring with probability one.} $\Xi$ let $\sup \bfxi = \sup \Xi$, which can be interpreted as the largest possible value of $\bfxi$.\\

\noindent{\bf Worst-case risk.} The choice $\cR(\bfxi) = \sup \bfxi$ represents a conservative outlook, possibly overly so because $\cR(\bfxi) = \infty$ when $\bfxi$ has, for example, a normal or exponential distribution. This risk measure ignores the probabilities of the various outcomes of $\bfxi$ and utilizes only the support $\Xi$. Thus, it is especially suitable when there is little information for building a probability distribution. A worst-case risk measure underpins robust optimization \cite{BentalElghaouiNemirovski.09,BertsimasBrownCaramanis.11}, semiinfinite programming \cite[Chapter 3]{Polak.97}, robust control \cite[Chapter 4]{Polak.97}, adversarial training \cite{MadryMakelovSchmidtTsiprasVladu.18}, and diametrical risk minimization \cite{NortonRoyset.22}.\\

\noindent{\bf Mean-plus-standard-deviation risk.} A natural choice motivated by statistical confidence intervals is to adopt a mean-plus-standard-deviation risk measure defined by $\cR(\bfxi) = \Ex[\bfxi] + \lambda \std(\bfxi)$, where $\std(\bfxi)$ is the standard deviation of $\bfxi$ and $\lambda$ is a positive constant. While mean-plus-standard-deviation risk measures can be traced back to early efforts in portfolio management \cite{Markowitz.52}, where in fact the variance was used instead of the standard deviation, their close ties to normal distributions are problematic. In particular, the two random variables from Figure \ref{fig:tri2} have the same mean and standard deviation and thus the same mean-plus-standard-deviation risk. This is caused by an equal treatment of outcomes above the mean and below the mean in the calculation of standard deviations, which is counter to many decision makers' greater concern about variability above the mean than below the mean.\\

Beyond mean and standard deviations, quantiles and superquantiles are central quantities for summarizing a random variable. Most importantly, they treat high values differently than low values. For $\alpha\in (0,1)$, the {\em $\alpha$-quantile} of a random variable $\bfxi$ with cumulative distribution function\footnote{For a random variable $\bfxi$ defined on $(\Omega,\cA,\bbP)$, its cumulative distribution function is $P(\xi) = \bbP(\{\omega\in \Omega~|~\bfxi(\omega) \leq \xi \})$.} $P$ is given by
\begin{equation}\label{eqn:quantile}
Q(\alpha) = \min \big\{\xi\in\reals ~\big|~ P(\xi) \geq \alpha\big\}
\end{equation}
and the {\em $\alpha$-superquantile} of an integrable random variable $\bfxi$ is given by
\begin{equation}\label{eqn:superquantile}
\bar Q(\alpha) = Q(\alpha) + \frac{1}{1-\alpha}\Ex\big[\max\{0,\bfxi-Q(\alpha)\}\big].
\end{equation}
We also define $\bar Q(0) = \Ex[\bfxi]$ and $\bar Q(1) = \nsup \bfxi$, which indeed are the limits of $\bar Q(\alpha)$ as $\alpha\downto 0$ and $\alpha\upto 1$, respectively; see, e.g., \cite[Theorem 2]{RockafellarRoyset.14}.\\

\noindent{\bf Quantile risk.} The choice of risk measure $\cR(\bfxi) = Q(\alpha)$, the $\alpha$-quantile of $\bfxi$, is widely used in financial engineering (see, e.g., \cite{ArtznerDelbaenEberHeath.99}) under the name ``value-at-risk'' with typical $\alpha$ near 1. In Figure \ref{fig:tri2}, $Q(0.8) = 0.106$ for density function $p_1$ and $Q(0.8) = 0.122$ for $p_2$. Thus, $p_1$ (corresponding to decision $x = 1$) is better according to this measure of risk. The conclusion is problematic because $p_1$ extends further to the right in Figure \ref{fig:tri2} than $p_2$, which means that worst outcomes are possible under decision $x = 1$ than under $x = -1$. In fact, this measure of risk ignores the {\em magnitude} of outcomes above $Q(\alpha)$ and thus fails to reflect the severity of poor outcomes.\\

\noindent{\bf Superquantile risk.} The choice of risk measure $\cR(\bfxi) = \bar Q(\alpha)$, the $\alpha$-superquantile of $\bfxi$, covers the whole range of possibilities from the risk-neutral $\cR(\bfxi) = \Ex[\bfxi]$ to the worst-case risk measure $\cR(\bfxi) = \sup \bfxi$ by adjusting $\alpha$ from 0 to 1. In Figure \ref{fig:tri2}, $p_1$ and $p_2$ produce superquantiles $\bar Q(0.8) = 0.404$ and $\bar Q(0.8) = 0.230$, respectively. An assessment of the two decisions underpinning Figure \ref{fig:tri2} based on $0.8$-superquantiles thus gives a decisive advantage to $p_2$ and the decision $x = -1$. This reversal compared to the conclusion reached in the previous paragraph is caused by the fact that superquantiles account for the magnitude of poor outcomes and therefore flag decision $x=1$ due to the possibility of exceptionally poor outcomes under $p_1$.\\

These risk measures are just examples, with many others given below. We already see that the choice of risk measure profoundly affects a decision and thus should, ideally, reflect a decision maker's preferences and any model ambiguity that might be present. A risk measure should also be computationally tractable because any complication comes on top of the challenges already associated with a quantity of interest. Section \ref{sec:regretrisk} reviews desirable properties for risk measures that both promote computations and reflect typical preferences among decision makers.

In reliability analysis, one often considers the probability of a random variable $\bfxi$ exceeding a
threshold $\tau$, i.e., $\prob\{\bfxi > \tau\}$. Such probabilities do {\em not} define measures of risk because they fail to produce $\xi$ when the random variable takes the value $\xi$ with probability one. For this reason, we prefer to treat such expressions separately in Subsection \ref{subsec:reliability} under the name {\em measures of reliability}.

\subsection{Connections with Utility Theory}\label{subsec:utility}

How to make a choice between two random variables has a long history, especially in economics. The classical prescription is due to von Neumann and Morgenstern \cite{vonNeumannMorgenstern.44} and their {\em expected utility theory}: Faced with a choice between random variables $\bfxi_1$ and $\bfxi_2$, adopt a utility function $u:\reals\to \reals$, and make the choice using $\Ex[u(\bfxi_1)]$ and $\Ex[u(\bfxi_2)]$, with higher numbers being preferred. For brief summaries of expected utility theory, we refer to \cite{RockafellarRoyset.15b,AbbasCadenbach.18,FrohlichWilliamson.22}.
With our focus on achieving lower values (of cost, damage, shortfall relative to a target,  etc.), there is a mismatch with the orientation of expected utility theory. This is easily rectified by reversing the sign of utility functions and thus producing {\em disutility functions}; see, e.g., \cite{Rockafellar.20}. For random variables $\bfxi_1$ and $\bfxi_2$ representing quantities of interest oriented toward lower values, we can define a disutility function $v:\reals\to \reals$ from a utility function $u:\reals\to\reals$ by setting $v(\xi) = -u(-\xi)$. Then, we may deem $\bfxi_1$ preferable over $\bfxi_2$ if $\Ex[v(\bfxi_1)] \leq \Ex[v(\bfxi_2)]$.

In stochastic programming models (see, e.g., \cite{BirgeLouveaux.11,ShapiroDentchevaRuszczynski.21,primer}), one often seeks to match a quantity of interest $f(\xi,x)$ with some target level, say 0. This produces the constraint $f(\xi,x) = 0$, which rarely can be satisfied due to the uncertainty associated with $\xi$. Simple recourse models \cite{Wets.74} circumvent this difficulty by  removing the constraint,  assuming that the values of the uncertain parameters are governed by the random vector $\bfxi$, and adding the term $\Ex[v(f(\bfxi,x))]$ to the objective function. Here, $v:\reals\to\reals$ can be viewed as a disutility function, often of the form $v(\eta) = \max\{\delta \eta, \gamma \eta\}$ with $\delta \leq 0 < \gamma$ being fixed parameters; see, e.g., \cite[Section 3.H]{primer}.

While expected utility theory and related approaches such as prospect theory \cite{KahnemanTversky.79} are widely adopted tools for risk-averse decision making, they do not automatically lead to measures of risk because $\Ex[v(\bfxi)] \neq \xi$ for a random variable $\bfxi$ concentrated at $\xi$; see \cite{DenuitDhaeneGoovaertsKaasLaeven.06,Rockafellar.20} for efforts to bridge this gap. These theories have also come under criticism \cite{FriedmanEtAl.14,FrohlichWilliamson.22}. The main concern is that the ``right'' utility function (or disutility function) is challenging to determine in practice, with more intricate functions resulting in difficulties with explaining how one reached a particular decision. There may also be no utility function that fully captures a decision maker's preference; she may be thinking about possible recourse actions that cannot be captured by utility functions alone. These concerns worsen when there are multiple stakeholders, which is common in many engineering applications. Moreover, there might be ambiguity about the underlying models and probability distributions making expected disutility less meaningful; see \cite{FrohlichWilliamson.22} for a detailed discussion.

\begin{example}\label{eHedging}{\rm (disutility functions and beyond).} Suppose that a firm faces a choice between two uncertain costs in the future described by the random variables $\bfxi_1$ and
$\bfxi_2$ distributed uniformly between $[-3/2, 1]$ and $[-8,2]$, respectively. (A negative cost means that the firm receives money.) The comparison between $\Ex[\bfxi_1] = -1/4$  and $\Ex[\bfxi_2]= -3$ leads to the conclusion that the better choice is to select $\bfxi_2$; on average it results in a lower cost. However, this may not be the right decision if the firm is concern about the possibility of a high cost, i.e., it is risk-averse.
\end{example}
\state Detail. Following expected (dis)utility theory, we may adopt a disutility function $v(\xi) = \max\{0, \xi\}/(1-\alpha)$, where $\alpha \in (0,1)$ is a  fixed parameter. The disutility function reflects a preference for outcomes $\xi \in (-\infty,0]$, with increasing displeasure associated with positive outcomes. This risk-averse perspective leads to a comparison between $\Ex[\max\{0, \bfxi_1\}/(1-\alpha)]$  and $\Ex[\max\{0, \bfxi_2\}/(1-\alpha)]$. These quantities can be thought of as the expected levels of displeasure felt by the firm when facing the random costs $\bfxi_1$ and $\bfxi_2$, respectively. Regardless of $\alpha$, the comparison produces a tie between the two choices because $\Ex[\max\{0, \bfxi_1\}]=\Ex[\max\{0, \bfxi_2\}]=1/5$.

Neither of these two approaches for assessing the merit of adopting one cost over the other considers the possibility of mitigating actions by the firm. It turns out that such deeper considerations may change the decision.
Let us suppose that the costs $\bfxi_1$ and $\bfxi_2$ are given in present money and that displeasure is quantified as above. If the firm is more active and invests $\gamma$ amount of money in a
risk-free asset (bonds or bank deposit) now, then the future displeasure, as perceived now, is reduced from
$\Ex[\max\{0, \bfxi_i\}/(1-\alpha)]$ to $\Ex[\max\{0, \bfxi_i-\gamma\}/(1-\alpha)]$ as $\gamma$ will be available at the future point in time to offset the cost $\bfxi_i$. The upfront expense $\gamma$ also needs to be considered and the goal becomes to select the investment $\gamma$ such that
\begin{equation}\label{eqn:regretminexample}
\gamma + \frac{1}{1-\alpha}\Ex\big[\max\{0, \bfxi_i-\gamma\}\big] ~\mbox{ is minimized.}
\end{equation}
As we see in Theorem \ref{tsuperquantile} below, the resulting minimum value is the $\alpha$-superquantile of $\bfxi_i$. A comparison between the $0.8$-superquantile of $\bfxi_1$, which is $3/4$, and the $0.8$-superquantile of $\bfxi_2$, which is $1$, reveals that $\bfxi_1$ is preferred. When accounting for the possibility of mitigating future displeasure by investing in a risk-free asset, the advantage tilts decisively toward the first cost. We see that a superquantile measure of risk inherently incorporates in its assessment of a random variable the possibility of such mitigation.\eop

The example illustrates the difference between making decisions based on expected (dis)utility theory and based on superquantile risk. The latter turns out to be deeply rooted in {\em dual utility theory} \cite{Yaari.87}, which relies on axioms parallel to those of expected utility theory; see  \cite{DentchevaRuszczynski.13,Rockafellar.20} and our discussion in Subsection \ref{subsec:mixed}. Moreover, superquantile risk has the following axiomatic justification \cite{WangZitikis.20}: For a real-valued risk measure $\cR$, defined on the integrable random variables, with $\cR(\bfxi) = 1$ for constant random variables $\bfxi$ concentrated at 1, we have that
\begin{align*}
&\cR \mbox{ satisfies the monotonicity, law invariance, prudence, and no-reward-for-concentration axioms }\\
&\Longleftrightarrow~~~~ \cR \mbox{ is a superquantile risk measure given by } \bar Q(\alpha) \mbox{ for some } \alpha \in (0,1).
\end{align*}
Thus, superquantile risk measures are the {\em only} risk measures with the listed axiomatic properties. (We define monotonicity and law invariance in Subsection \ref{subsec:desirable}, and prudence relates to lower semicontinuity as also defined in that section. For details about prudence and no-reward-for-concentration, we refer to \cite[Sections 2.1-2.2]{WangZitikis.20}.) This theoretical foundation and their practical usage support the claim that superquantile risk measures are ``currently the most important risk measure in banking practice'' \cite{WangZitikis.20}.

\section{Superquantiles}\label{sec:super}

In the vast landscape of risk measures, superquantiles emerge as central. They capture risk-averseness and one-sided concerns about high values. They span the range of preferences from the risk-neutral perspective (by setting $\alpha = 0$) to a focus on worst-case outcomes (by setting $\alpha = 1$). With the dependence on a single parameter ($\alpha$), superquantiles are easy to explain to decision makers. In the following, we identify three other features: (i) superquantiles have mathematical properties that facilitate computational optimization, (ii) they capture ambiguity about the adopted probability distribution and thus connect with distributionally robust optimization, and (iii) they furnish the fundamental building blocks for many ``reasonable'' measures of risk. Thus, we proceed with a comprehensive review of superquantiles; Sections \ref{sec:regretrisk} and \ref{sec:addexamples} cover more general risk measures.

\subsection{Equivalent Formulas}\label{subsec:equivalence}

Superquantiles can at least be traced back to the concept of {\em optimized certainty equivalents} as developed in  \cite{BentalTeboulle.86} for random variables oriented toward high value as preferable. Converting the concept into the present orientation, this pioneering paper effectively defines the quantity
\[
\min_{\gamma \in \reals} \,\gamma + \Ex\big[v(\bfxi - \gamma)\big],
\]
where $v:\reals\to\reals$ is an increasing, strictly convex, and twice continuously differentiable disutility function, normalized with $v(0) = 0$ and derivative $v'(0) = 1$. While the increasing and smoothness properties as well as strict convexity are violated by $v(\xi) = \max\{0,\xi\}/(1-\alpha)$, this particular choice brings us to the minimization formula for an $\alpha$-superquantile used in conjunction with \eqref{eqn:regretminexample} and formally expressed as \eqref{eqn:superquantRU} below; see \cite{BentalTeboulle.07} for details about these connections.

Averages beyond a quantile, which are closely related to \eqref{eqn:superquantile}, are mentioned in \cite{ArtznerDelbaenEberHeath.97, ArtznerDelbaenEberHeath.99} under the name tail value-at-risk; see also \cite{Embrechts.99}. Using the term conditional value-at-risk (CVaR), \cite{RockafellarUryasev.00} also starts from an average beyond a quantile and derives the equivalence below between \eqref{eqn:superquant} and \eqref{eqn:superquantRU} under the assumption that the random variable $\bfxi$ is continuously distributed. The equivalence between \eqref{eqn:superquant} and \eqref{eqn:superquantRU} for general, integrable random variables is confirmed in \cite{Pflug.00}; see also \cite{RockafellarUryasev.02}, which further establishes a connection with the modified tail expectation \eqref{eqn:superquantTail}.
This modified tail expectation accounts for random variables with a positive probability of taking a value exactly at a quantile and is also the starting point for \cite{AcerbiTasche.02}. That paper proceeds by establishing equivalence with the formula \eqref{eqn:superquantAcerbi} from \cite{Acerbi.02} involving an integral of quantiles across different probability levels as well as with the minimization formula \eqref{eqn:superquantRU} from \cite{RockafellarUryasev.00}. Integrals of quantiles have a long history in statistics and underpin Lorenz curves \cite{Lorenz.05}; see, for example, the discussion in Section 3 of \cite{RockafellarRoyset.14} for additional connections. Further insight stems from \cite{Delbaen.00}, which expresses a superquantile of a continuous random variable as the worst-case expectation over a family of probability distributions. The situation for general distributions is hinted to in \cite{Delbaen.00}, but brought out more clearly in \cite{RockafellarUryasevZabarankin.02}. In summary, the flurry of activity around the turn of the century produced the following equivalent formulas for an $\alpha$-superquantile.

\begin{theorem}\label{tsuperquantile} {\rm (equivalent formulas for superquantiles).}
For $\alpha\in (0,1)$ and an integrable random variable $\bfxi$ with cumulative distribution function $P$ and quantile function $Q$, the following hold:
\begin{align}
&Q(\alpha) + \frac{1}{1-\alpha}\Ex\big[\max\{0,\bfxi-Q(\alpha)\}\big]\label{eqn:superquant}\\
& = \min_{\gamma\in \reals} \,\gamma + \frac{1}{1-\alpha}\Ex\big[\max\{0,\bfxi-\gamma\}\big]\label{eqn:superquantRU}\\
& = \frac{1}{1-\alpha}\int_\alpha^1 Q(\beta)d\beta\label{eqn:superquantAcerbi}\\
& = \mbox{expectation of the $\alpha$-tail distribution of $\bfxi$}\label{eqn:superquantTail},
\end{align}
where the $\alpha$-tail distribution is defined as having $P^{[\alpha]}(\xi) = \max\{ 0, P(\xi)-\alpha\}/(1-\alpha)$ as its cumulative distribution function. Thus, any of these formulas can be taken as the definition of the $\alpha$-superquantile of $\bfxi$, which we denote by $\bar Q(\alpha)$ or $\srsk_\alpha(\bfxi)$ to highlight the dependence on $\bfxi$.

If $\bfxi$ is square integrable, then 
\begin{equation}\label{eqn:argminsuper}
\bar Q(\alpha) = \nargmin_{\gamma\in \reals}  \bigg\{\gamma + \frac{1}{1-\alpha}\int_0^1 \max\big\{0, \bar Q(\beta) - \gamma\big\} d\beta \bigg\}.
\end{equation}

If $\bfxi$ has finite support $\Xi$ of cardinality $s$ and corresponding probabilities $p_\xi, \xi\in \Xi$, then
\begin{equation}\label{eqn:maxsuper}
\bar Q(\alpha) = \max_{\bar p \in \Delta_\alpha} \nsum_{\xi\in \Xi} \xi\, \bar p_\xi, ~\mbox{ where } ~
\Delta_\alpha = \Big\{ \bar p\in \reals^s~\Big|~ 0 \leq \bar p_\xi \leq \frac{p_\xi}{1-\alpha}, \xi\in \Xi, ~\nsum_{\xi\in \Xi} \bar p_\xi = 1 \Big\}.
\end{equation}
\end{theorem}
\state Proof. While originally developed by several researchers as discussed before the theorem, the four equivalences are concisely summarized in \cite{RockafellarRoyset.14}; see the discussion around equation (3.4) in that reference for the equivalence between \eqref{eqn:superquantAcerbi} and \eqref{eqn:superquantTail} above. Section 4 and Theorem 2 of \cite{RockafellarRoyset.14} confirm the other equivalences.

Theorem 7 in \cite{RockafellarRoyset.14} establishes the argmin formula \eqref{eqn:argminsuper}. A proof of the last expression \eqref{eqn:maxsuper} appears in \cite{RockafellarUryasevZabarankin.02}; see also our discussion in Section \ref{sec:duality}.\eop

In general, the $\alpha$-superquantile $\bar Q(\alpha)$ of a random variable $\bfxi$ is equal to neither $\Ex[\bfxi~|~\bfxi \geq Q(\alpha) ]$ nor $\Ex[\bfxi~|~\bfxi > Q(\alpha) ]$ and this sometimes causes confusion. The discrepancy is reflected in \eqref{eqn:superquantTail}, where any probability atom  at the $\alpha$-quantile $Q(\alpha)$ of $\bfxi$ is carefully ``split.''  Nevertheless, in the absence of such an atom, $\bar Q(\alpha) = \Ex[\bfxi~|~\bfxi \geq Q(\alpha) ] = \Ex[\bfxi~|~\bfxi > Q(\alpha) ]$. In particular, if $\bfxi$ has a density function $p$, then
\begin{align}
\bar Q(\alpha) &= Q(\alpha) + \frac{1}{1-\alpha}\int_{-\infty}^\infty \max\big\{0, \xi - Q(\alpha)\big\} p(\xi) \,d\xi\nonumber\\
&= Q(\alpha) + \frac{1}{1-\alpha}\int_{Q(\alpha)}^\infty \xi \, p(\xi) \,d\xi - \frac{Q(\alpha)}{1-\alpha}\int_{Q(\alpha)}^\infty  p(\xi) \,d\xi\nonumber\\
& = \frac{1}{1-\alpha}\int_{Q(\alpha)}^\infty \xi \, p(\xi) \,d\xi,\label{eqn:condsuper}
\end{align}
which coincides with $\Ex[\bfxi~|~\bfxi \geq Q(\alpha) ]$ and $\Ex[\bfxi~|~\bfxi > Q(\alpha) ]$.

These finer points may sometimes be glossed over. Based on our experience with practitioners, we recommend adopting \eqref{eqn:superquant} as the definition of a superquantile, with the supporting remark:
\[
\mbox{$\alpha$-superquantile of $\bfxi ~=\, $ average of the worst } (1-\alpha)100\% \mbox{ outcomes of } \bfxi.
\]
The word ``worst'' is intuitively understood by the practitioner and is sufficiently ambiguous to provide cover for the mathematician.

The importance of \eqref{eqn:superquantRU} in computations emerges in Subsection \ref{subsec:superinopt}. It is apparent from \eqref{eqn:superquant} and \eqref{eqn:superquantRU} that the $\alpha$-quantile $Q(\alpha)$ is a minimizer of the optimization problem over $\gamma\in \reals$ in the latter formula. However, it may not be the only minimizer. There are multiple minimizers if the cumulative distribution function $P$ for $\bfxi$ has a ``flat stretch'' to the right of $Q(\alpha)$ and $P(Q(\alpha)) = \alpha$. For any $\alpha\in (0,1)$, we obtain that
\begin{equation}\label{eqn:argminquant}
\big[Q(\alpha), Q^\lplus(\alpha)\big] = \nargmin_{\gamma\in \reals} \bigg\{\gamma + \frac{1}{1-\alpha}\Ex\big[\max\{0,\bfxi-\gamma\}\big]\bigg\},
\end{equation}
where $Q^\lplus(\alpha) =  \sup\{\gamma \in\reals ~| ~P(\gamma) \leq \alpha\}$; see, e.g., \cite[Equation (4.1)]{RockafellarRoyset.14}.

The insight leading to \eqref{eqn:superquantRU} in \cite{RockafellarUryasev.00}, as recounted in \cite{RockafellarRoyset.14}, was that an integrable random variable $\bfxi$ can be associated with the convex function given by $e(\gamma) = \Ex[\max\{\gamma,\bfxi\}]$ whose subgradients recover the cumulative distribution function of the random variable. The conjugate function $e^*$ of the convex function $e$, as given by $e^*(\alpha) = \sup_{\gamma\in\reals} \, \alpha \gamma - e(\gamma)$, 
turns out to furnish \eqref{eqn:superquantRU} after a scaling with $\alpha - 1$. Specifically, for $\alpha \in (0,1)$, one has
\[
\frac{e^*(\alpha)}{\alpha-1} = \min_{\gamma\in\reals} \, \gamma + \frac{1}{1-\alpha} \Ex\big[\max\{0, \bfxi - \gamma\}\big]
\]
as seen in \cite[Theorem 2]{RockafellarRoyset.14}. Similar relations between conjugate pairs and their connections with cumulative distribution functions and quantile functions were examined independently in \cite{OgryczakRuszczynski.02}; see also \cite{DentchevaMartinez.12}.

The argmin-formula \eqref{eqn:argminsuper} in Theorem \ref{tsuperquantile} may at first appear less useful as it requires the knowledge of all superquantiles to compute one of them. However, it is the linchpin for superquantile regression (see \cite{RockafellarRoysetMiranda.14} and our discussion in Subsection \ref{subsec:superregression}) and allows us to estimate one superquantile from a finite number of other ones via numerical integration. Additional argmin-formulas for superquantiles appear in \cite{Kouri.19b}.

As we see in Section \ref{sec:duality}, the formula \eqref{eqn:maxsuper}, referred to as the {\em dual formula}, holds much beyond finite distributions. Still, recording this special case is useful as it avoids all technical overhead while addressing important applications in data-driven optimization and learning. The key insight from \eqref{eqn:maxsuper} is that two seemingly different decision makers will make the same assessment of a random variable. Specifically, for $\alpha \in (0,1)$ and random variables with finite support, consider the two individuals:
\begin{quote}
Mr.~Averse has full confidence in the assumed probability distributions of random variables, but is inherently risk-averse and makes decision by comparing $\alpha$-superquantiles of the random variables.
\end{quote}
\begin{quote}
Ms.~Ambiguous is risk-neutral and makes decisions based on expectations, but is suspicious about the assumed probability distributions. She computes an expectation using the worst-case probability distribution obtained by scaling the assumed probabilities with factors between $0$ and $1/(1-\alpha)$.
\end{quote}
In effect, Mr.~Averse computes the left-hand side of \eqref{eqn:maxsuper} and Ms.~Ambiguous computes the right-hand side. Thus, they reach the same assessment. We conclude that superquantiles can be used to address a decision maker's inherent risk-averseness as well as ambiguity about probability distributions, for example due to lack of data, fear of contamination, or adversarial interference.

It follows almost immediately from the definitions that the superquantiles $\bar Q(\alpha)$ of an integrable random variable $\bfxi$ tend to the worst-case $\sup \bfxi = \bar Q(1)$ as $\alpha\upto 1$. A useful quantification of the difference $\bar Q(1)-\bar Q(\alpha)$ appears in \cite{AndersonXuZhang.20}. For further properties of superquantiles, we refer to \cite{RockafellarRoyset.14}.

With its origin in financial engineering, superquantiles are also known as {\em tail value-at-risk}, {\em expected shortfall}, {\em average value-at-risk}, and {\em conditional value-at-risk} (with CVaR as an abbreviation). Despite slight difference in definitions originally, we now accept these as synonyms for the quantity here called superquantiles. The name superquantile stems from \cite{RockafellarRoyset.10} and was motivated by a need for making this fundamental concept free from dependence on financial terminology. 

We give explicit formulas for superquantiles in three cases next; see \cite{NortonKhokhlovUryasev.21} for many other cases involving common probability distributions.

\begin{example}\label{eSuperTriang}{\rm (triangular distribution).} Given $\alpha\in (0,1)$, the
$\alpha$-quantiles and $\alpha$-superquantiles of a random variable $\bfxi$ with the triangular density function $p_1$ in Figure \ref{fig:tri2} are
\[
Q(\alpha) = 1-2\sqrt{1-\alpha} ~~~\mbox{ and }~~~ \bar Q(\alpha) = 1-\fourthird\sqrt{1-\alpha}.
\]
Moreover, $\bar Q(0) = -1/3$ and $\bar Q(1) = 1$.
\end{example}
\state Detail. Let $P$ be the cumulative distribution function of $\bfxi$. Since $p_1(\xi) = -\xi/2 + 1/2$ for $\xi\in [-1,1]$, the solution of the equation $P(\xi) = \int_{-1}^\xi p_1(\eta) \,d\eta = \alpha$ is $Q(\alpha) = 1-2\sqrt{1-\alpha}$. This formula and \eqref{eqn:condsuper} give $\bar Q(\alpha)$.\eop

\begin{example}\label{eSuperNormal}{\rm (normal distribution).} Given $\alpha \in (0,1)$, the
$\alpha$-quantiles and $\alpha$-superquantiles of a random variable that is normally
distributed with mean $\mu$ and variance $\sigma^2$ are
\[
Q(\alpha) = \mu + \sigma \Phi^{-1}(\alpha)~~~ \mbox{ and }~~~ \bar Q(\alpha) = \mu + \frac{\sigma\phi\big(\Phi^{-1}(\alpha)\big)}{1-\alpha},
\]
where $\phi$ is the standard normal density function and $\Phi^{-1}(\alpha)$ is the corresponding $\alpha$-quantile given by the standard normal cumulative distribution function $\Phi$. Moreover, $\bar Q(0) = \mu$ and $\bar Q(1) = \infty$.
\end{example}

\begin{example}\label{eSuperFinite}{\rm (finite distribution).} For $\alpha \in [0,1]$, the $\alpha$-superquantile of a finitely distributed random variable with values $\xi_1<\xi_2< \dots< \xi_r$, which occur with probabilities $p_1, p_2,$ $\dots$, $p_r$, respectively, is given by
\begin{equation*}
\bar Q(\alpha) = \begin{cases} \nsum_{j=1}^r p_j \xi_j & \mbox{if } \alpha = 0\\
\frac{1}{1-\alpha}\left(\Big(\big(\nsum_{j=1}^i p_j\big) - \alpha\Big)\xi_i + \nsum_{j=i+1}^n p_j \xi_j\right) & \mbox{if } \nsum_{j=1}^{i-1} p_j <\alpha\leq \nsum_{j=1}^i p_j<1\\
\xi_r & \mbox{if } \alpha > 1-p_r.
\end{cases}
\end{equation*}
\end{example}

\subsection{Superquantiles in Optimization Models}\label{subsec:superinopt}

Some of the formulas in Theorem \ref{tsuperquantile} might leave the impression that superquantiles are more complicated to compute and optimize than quantiles. However, this is not the case. The formula \eqref{eqn:superquantRU} decouples superquantiles from quantiles and this turns out to be especially important in optimization models.

Consider a quantity of interest $f:\reals^m\times \reals^n\to \reals$ and an $m$-dimensional random vector $\bfxi$. Suppose that for each $x\in\reals^n$, $f(\bfxi,x)$ is an integrable random variable. With
\[
\srsk_\alpha\big(f(\bfxi,x)\big) = \alpha\mbox{-superquantile of } f(\bfxi,x),
\]
we may seek a decision $x$ that minimizes $\srsk_\alpha(f(\bfxi,x))$. If the minimum value turns out to be $\tau$, then the obtained decision $x^\star$ has the guarantee that $f(\bfxi,x^\star) \leq \tau$ on average  across the worst $(1-\alpha)100\%$ outcomes. Moreover, the alternative formula \eqref{eqn:maxsuper} leads to the insight that $x^\star$ is a decision that minimizes the worst-case expected value of the quantity of interest across a set of probability distributions ``near'' the nominal one. The resulting decision would typically be rather different than those obtained by minimizing $\Ex[f(\bfxi,x)]$ under the nominal probability distribution, which pay no particular attention to the possibility of high values of $f(\bfxi,x)$ or deviations from the nominal distribution.

Superquantiles preserve the important convexity property as recognized in \cite{RockafellarUryasev.00}; here we recall \cite[Proposition 3.10]{primer}.

\begin{proposition}\label{pSuperquantiles2}{\rm (convexity of superquantile functions).}
For a random vector $\bfxi$ with support $\Xi\subset\reals^m$ and
$f:\reals^m\times \reals^n\to\reals$, suppose that
\begin{enumerate}[{\rm (a)}]

\item $f(\bfxi,x)$ is integrable for all $x\in\reals^n$

\item $f(\xi,\cdot\,)$ is convex for all $\xi\in \Xi$.

\end{enumerate}
Then, for any $\alpha\in [0,1)$, the function $x\mapsto \srsk_\alpha (f(\bfxi,x))$ is convex and real-valued.
\end{proposition}

If we ignore the uncertainty associated with $\bfxi$ and simply consider minimizing $f(\hat \xi, x)$ for some nominal value $\hat \xi$, then we would face a convex problem as long as $f(\hat \xi,\cdot\,)$ is convex. The proposition asserts that $x\mapsto \srsk_\alpha (f(\bfxi,x))$, which treats uncertainty much more comprehensively, is convex too when $f(\xi,\cdot\,)$ is convex for all $\xi\in \Xi$. Thus, optimization under uncertainty carried out in this manner does not muddle up convexity that might be present in $f$.

We next consider computational approaches. Given a feasible set $X\subset\reals^n$ and $\alpha \in (0,1)$, suppose that we seek to solve the optimization problem
\begin{equation}\label{eqn:sriskmin}
\nnmin_{x\in X} ~\srsk_\alpha\big(f(\bfxi,x)\big).
\end{equation}
The formula \eqref{eqn:superquantRU} enables us to reformulated this problem equivalently as
\begin{equation}\label{eqn:sriskmingamma}
\nnmin_{x\in X,\gamma\in\reals} ~\Ex\Big[\gamma + \frac{1}{1-\alpha}\max\big\{0,f(\bfxi,x)-\gamma\big\}\Big],
\end{equation}
which brings us back to minimizing an expectation function for which there are many algorithms including (stochastic) subgradient type methods (cf. Subsection \ref{eqn:algosuper}). Moreover, if $\bfxi$ is finitely distributed with support $\Xi$ and corresponding probabilities $\{p_{\xi}>0, \xi \in \Xi\}$, then the problem simplifies to
\begin{equation}\label{eqn:superquantExpanded0}
\nnmin_{x\in X,\gamma\in\reals} ~\gamma + \frac{1}{1-\alpha}\nsum_{\xi\in \Xi} p_{\xi}\max\big\{0,f(\xi,x)-\gamma\big\},
\end{equation}
which in turn can be reformulated as
\begin{equation}\label{eqn:superquantExpanded}
\nnmin_{x\in X,\gamma\in\reals,z\in\reals^{s}} ~\gamma + \frac{1}{1-\alpha}\nsum_{\xi\in \Xi} p_{\xi} z_\xi 
~\mbox{ subject to } ~f(\xi,x)-\gamma \leq z_\xi, ~~~ 0 \leq z_\xi ~~~\forall \xi\in \Xi,
\end{equation}
where $z = (z_\xi, \xi\in \Xi)\in \reals^s$ are additional variables and $s$ is the cardinality of $\Xi$. Both \eqref{eqn:superquantExpanded0} and \eqref{eqn:superquantExpanded} are convex as long as $X$ is convex and $f(\xi, \cdot\,)$ is convex for all $\xi\in \Xi$. In fact, if $f(\xi,\cdot\,)$ is affine for all $\xi\in \Xi$ and $X$ is polyhedral, then  \eqref{eqn:superquantExpanded} is a linear optimization problem.

Smoothness is also preserved. If $f(\xi,\cdot\,)$ is continuously differentiable for all $\xi\in \Xi$, then \eqref{eqn:superquantExpanded} involves inequality constraints of the kind commonly addressed by nonlinear programming algorithms.

In light of \eqref{eqn:argminquant}, we see that the $\gamma$-portion of any minimizer $(x^\star,\gamma^\star,z^\star)$ from \eqref{eqn:superquantExpanded} is the $\alpha$-quantile of $f(\bfxi,x^\star)$ or possibly a (slightly) larger quantity as specified by \eqref{eqn:argminquant}. The minimum value from \eqref{eqn:superquantExpanded} is of course the $\alpha$-superquantile of $f(\bfxi,x^\star)$. The advantages of \eqref{eqn:superquantExpanded} are therefore clear, but the reformulation only applies to finite distributions.

A superquantile appearing as a constraint is treated analogously.  For $f_0:\reals^n\to \reals$ and $\tau\in \reals$, the problem
\begin{equation}\label{eqn:sriskmincon}
\nnmin_{x\in X} ~f_0(x) ~~\mbox{ subject to } ~~\srsk_\alpha\big(f(\bfxi,x)\big) \leq \tau
\end{equation}
is equivalently stated as
\begin{align}\label{eqn:superquantExpandedcon}
\nnmin_{x\in X,\gamma\in\reals,z\in\reals^{s}} ~f_0(x) ~~ \mbox{ subject to } ~~\gamma + \frac{1}{1-\alpha}\nsum_{\xi\in \Xi} p_{\xi} z_\xi &\leq \tau\\
 f(\xi,x)-\gamma \leq z_\xi, ~~~ 0 & \leq z_\xi ~~~\forall \xi\in \Xi,\nonumber
\end{align}
again provided that $\bfxi$ has a finite distribution with $s$ outcomes. Following a similar pattern, we  achieve formulations for multiple quantities of interest as well.

The convexity property associated with superquantiles stands in sharp contrast to the situation for quantiles. The following illustration is taken from \cite[Example 3.12]{primer}.

\drawing{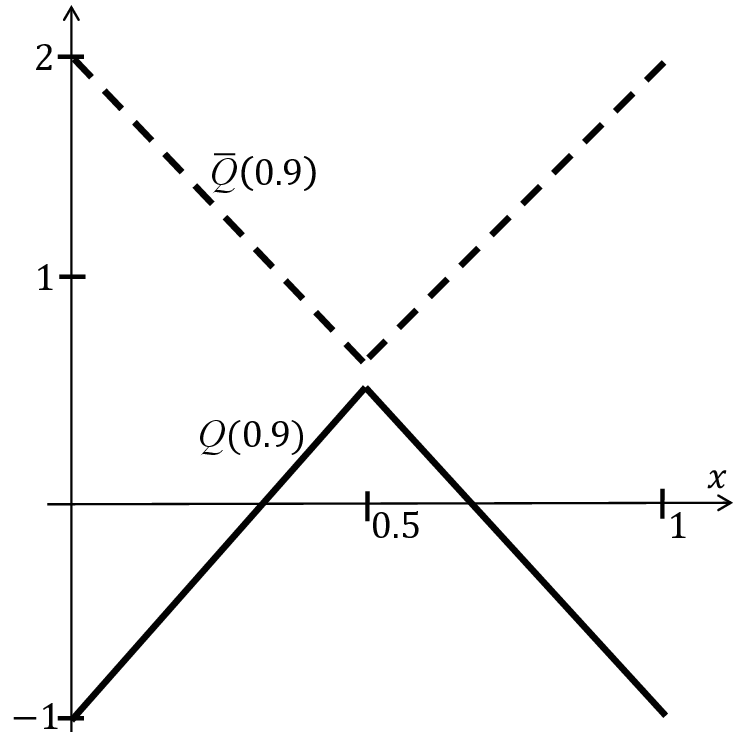}{2.8in} {Quantile $Q(0.9)$ (solid line) and superquantile
$\bar Q(0.9)$ (dashed) of $f(\bfxi,x)$ in Example \ref{eQuantvsSuperquant} as functions of $x$.}{fig:quantexamp}

\begin{example}{\rm (quantile vs superquantile minimization).}\label{eQuantvsSuperquant}
Suppose that we need to allocate funds between two financial instruments, both associated with uncertainty. The first instrument requires, with probability 0.1, a payment of 2 dollars per dollar committed and, with probability 0.9, yields one dollar per dollar committed. Let the random
variable $\bfxi_1$ model these losses; it takes the value 2 with probability 0.1 and the value $-1$ with probability 0.9. The second instrument is modeled by the random variable $\bfxi_2$, which has the same probability distribution as that of $\bfxi_1$ but the two random variables are statistically
independent. Suppose that we need to allocate 1 million dollars between these two instruments.  Since the instruments appear equally unappealing, we might be led to believe that any allocation is fine. This would be a big mistake, but one that remains hidden for an analyst examining quantiles. Superquantiles reveal that the best strategy would be to allocate half a million to each instrument.
\end{example}
\state Detail. Let $x \in [0,1]$ be the fraction of our million dollars allocated to the first instrument,
the remainder is allocated to the second instrument. Then, the quantity of interest describing our loss is  the random variable $f(\bfxi,x) = x\bfxi_1 + (1-x)\bfxi_2$. Since there are only four possible outcomes of $\bfxi = (\bfxi_1,\bfxi_2)$, we find that
\[
f(\bfxi,x) = \begin{cases}
  -1 & \mbox{ with probability } 0.81\\
  2-3x & \mbox{ with probability } 0.09\\
  3x-1 & \mbox{ with probability } 0.09\\
  2 & \mbox{ with probability } 0.01.
\end{cases}
\]
With $\alpha = 0.9$, the $\alpha$-quantile of $f(\bfxi,x)$ becomes
\[
Q(0.9) = \begin{cases}
3x -1 & \mbox{ if } x\in [0,1/2]\\
2-3x & \mbox{ otherwise},
\end{cases}
\]
which is depicted with a solid line in Figure \ref{fig:quantexamp}. By \eqref{eqn:superquant}, the $0.9$-superquantile of $f(\bfxi,x)$ becomes
\[
\bar Q(0.9) = \begin{cases}
-2.7x + 2 & \mbox{ if } x\in [0,1/2]\\
2.7x-0.7 & \mbox{ otherwise},
\end{cases}
\]
which is also shown in Figure \ref{fig:quantexamp} with a dashed line. As a function of $x$, the superquantiles define a convex function but the quantiles do not. The minimization of superquantiles has $x^\star = 0.5$ as minimizer, but the minimization of quantiles results in $x = 0$ and $x = 1$. The decision $x^\star$ involves {\em hedging}, a well-known strategy in finance to reduce risk. The decision lowers the probability of a loss of 2 million dollars from 0.1 to 0.01 compared to the choice $x = 0$ or $x = 1$. This comes at the
expense of reducing the probability of a loss of $-1$ million dollars from 0.9 to 0.81. The decision $x^\star$ also results in the possibility of a loss of 0.5 million dollars (with probability 0.18), but this might be much more palatable than a 2-million-dollar loss. The minimum value of the 0.9-superquantiles is 0.65 million dollars, which is a more reasonable and conservative assessment of the uncertain future loss than the wildly optimistic $-1$ million provided by the minimum value of the 0.9-quantiles.\eop

The advantages of superquantiles over quantiles are not limited to convexity properties and hedging strategies as indicated in the previous example. Using \cite[Example 3.13]{primer}, we next illustrate that superquantiles reveal the magnitude of ``typical'' poor outcomes and not only their likelihood.

\drawing{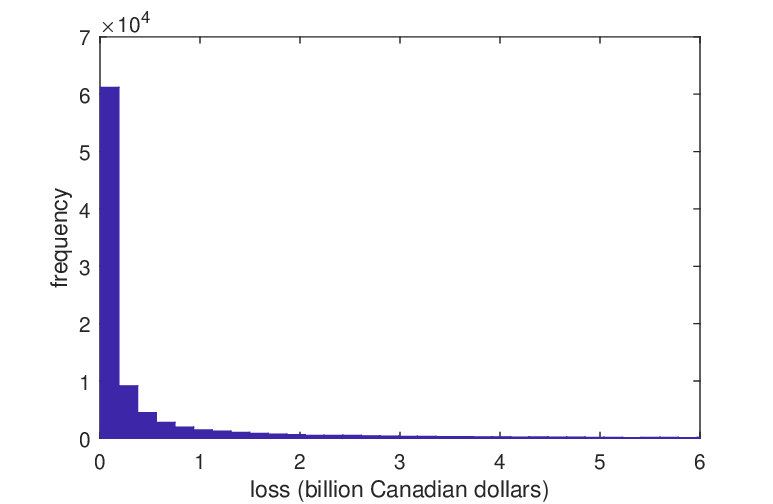}{4.0in}{Histogram of 100000 outcomes of the cumulative loss due to earthquake damage in Vancouver.}{fig:eqloss}

\begin{example}{\rm (revealing vulnerabilities).}\label{eResilience}
Superquantiles have the ability to reveal vulnerabilities that might remain hidden if the focus is on quantiles. From a study of the cumulative loss due to earthquake damage during the next 50 years in the greater Vancouver region in Canada \cite{Mahsuli.12}, we obtain 100000 loss values that one can view as outcomes for a random variable modeling the cumulative loss; see the  histogram of Figure \ref{fig:eqloss}. The cumulative loss is most likely less than one billion Canadian dollars. In fact, the 0.9-quantile of the cumulative loss is 4.39 billions. Although sizable, it seems that the region is quite resilient to earthquakes. It turns out that this conclusion is flawed and a study of superquantiles paints a gloomier picture.
\end{example}
\state Detail. The largest outcome is actually 373 billions; Figure \ref{fig:eqloss} should have been extended much to the right. Superquantiles quantify this long right-tail. Specifically, the 0.9-superquantile of the cumulative loss is 28.92 billions. This means that when the cumulative loss exceeds 4.39 billons, it does so substantially. The region might have a significant vulnerability after all.\eop

The advantageous properties of a superquantile, both computationally and from a modeling perspective, have brought this measure of risk to the forefront. From its start in financial applications \cite{RockafellarUryasev.00,KrokhmalPalmquistUryasev.02} (see \cite{AlotaibiDallavalleCraven.22,NasiniLabbeBrotcorne.22} for examples of recent work in that area), superquantiles have been used in retail operations \cite{LuShanthikumarShen.15}, in planning military jamming missions \cite{CommanderPardalosRyabchenkoUryasevZrazhevsky.07}, in supporting capital investment decisions \cite{TeterRoysetNewman.19}, and in inventory control \cite{AhmedCakmakShapiro.07,HanasusantoKuhnWallaceZymler.15}; see \cite{FilippiGuastarobaSperanza.20,LuShen.21} for a summary from the area of operations management broadly.
Applications of superquantiles in the energy sector are exemplified by \cite{GhasemiJamshidimonfaredLoniMarzband.21,XuanShenGuoSun.21}. Superquantiles appear also in the control of PDEs \cite{KouriSurowiec.16,GarreisSurowiecUlbrich.21} and design of physical systems \cite{RoysetBonfiglioVernengoBrizzolara.17,BonfiglioRoyset.19,ChaudhuriNortonKramer.20,ChaudhuriPeherstorferWillcox.20}.

Superquantiles are increasingly being used in machine learning. For computer vision problems,
\cite{LevyCarmonDuchiSidford.20,Yang.21} report promising empirical results. The paper \cite{PillutlaLaguelMalickHarchaoui.21} shows the potential in federated learning with simulations in the areas of character recognition and sentiment analysis; see \cite{LaguelPillutlaMalickHarchaoui.21b} for distributed learning on mobile devices involving heterogeneous distributions. Modeling using superquantiles also emerges as a tool to address fairness in machine learning \cite{WilliamsonMenon.19,FrohlichWilliamson.22} and, after slight adjustments, also outliers \cite{LiuPang.22}. Superquantiles in reinforcement learning are discussed in \cite{MorimuraSugiyamaKashimaHachiyaTanaka.10,SinghZhangChen.20,UrpiCuriKrause.21}. Additional recent efforts include \cite{SomaYoshida.20,CuriLevyJegelkaKrause.20}; see also \cite{FanLyuYingHu.17,KawaguchiLu.19} for closely related top $k$-approaches, with promising numerical results, and \cite{LiBeiramiSanjabiSmith.21} for a related tilted loss approach. We refer to the recent review article \cite{LaguelPillutlaMalickHarchaoui.21} for a summary of superquantile applications in machine learning. A general statistical estimation point of view, especially in the broad context of distributional shifts, appears in \cite{DuchiNamkoong.21}.

Connections between superquantiles and distributionally robust optimization with a Wasserstein ambiguity set emerge in \cite{HonguyenWright.22}; see also our discussion in Subsection \ref{subsec:distributionallyrobust}. Superquantiles relate to the classical newsvendor problem from Example \ref{mNwsVnd} \cite{GotohTakano.07} as well as support vector machines for binary classification \cite{TakedaSugiyama.08,GotohTakedaYamamoto.14,GotohUryasev.17}.

For applications in the area of two-stage stochastic mixed-integer programming, we refer to
\cite{SchultzTiedemann.06,vanBeestenRomeijnders.20} and the paper \cite{FernandezHinojosaPuertoSaldanhadagama.19}, which focuses on fixed-charge transportation problems; see also \cite{HemmatiSabooriSaboori.16} for applications to energy storage and transportation, \cite{ToumazisKwon.16} for routing hazardous material, and \cite{Noyan.12} for modeling of disaster relief.
Superquantiles even define norms as discussed in \cite{PavlikovUryasev.14,MafusalovUryasev.16,
GotohUryasev.16} and quantify the distance between cumulative distribution functions \cite{PavlikovUryasev.18}.

\subsection{Algorithms for Superquantile Minimization}\label{eqn:algosuper}

The basic approach for minimizing a superquantile risk over a feasible set $X$ takes us from \eqref{eqn:sriskmin} through a reformulation to \eqref{eqn:superquantExpanded}, with a similar treatment when a superquantile appears as a constraint; see \eqref{eqn:sriskmincon} and \eqref{eqn:superquantExpandedcon}.
The approach is appealing because it can leverage any state-of-the-art algorithm for solving the resulting problem. However, it breaks down if the random vector $\bfxi$ does not have a finite distribution or if its finite distribution involves a massively large number of outcomes. In this subsection, we survey more advanced algorithms for superquantile minimization.\\

\noindent {\bf Stochastic Subgradient Methods.} The superquantile minimization problem \eqref{eqn:sriskmin}  can always be written as \eqref{eqn:sriskmingamma} using the formula \eqref{eqn:superquantRU}. Thus, we are back in the familiar domain of expectation minimization. Let
\[
\psi\big(\xi,(x,\gamma)\big) = \gamma + \frac{1}{1-\alpha} \max\big\{0, f(\xi,x) - \gamma\big\},
\]
which thus defines the integrand in \eqref{eqn:sriskmingamma}, and $\phi(z) = \Ex[\psi(\bfxi,z)]$. We can then apply a standard stochastic subgradient method to $\phi$, among which the SGD\footnote{SGD stands for ``stochastic gradient descent'' but this is doubly misleading as it is neither a descent method nor involves gradients only; $\psi(\xi, \cdot\,)$ is nonsmooth and thus we need to consider subgradients.} method is well known.\\

\state SGD Method.

\begin{description}
  \item[Data.] ~~\,$z^0\in X\times\reals$ and step sizes $\lambda^\nu \in (0,\infty), ~\nu=0, 1, 2, \dots$.

  \item[Step 0.]  Set iteration counter $\nu = 0$.

  \item[Step 1.]  Generate an observation $\xi^\nu$ according to the distribution of $\bfxi$.

  \item[Step 2.]  Compute a subgradient $v$ of $\psi(\xi^\nu, \cdot\,)$ at the point $z^\nu$ and set
  \[
  z^{\nu+1} \in \prj_{X\times \reals}(z^\nu - \lambda^\nu v).
  \]

  \item[Step 3.]  Replace $\nu$ by $\nu +1$ and go to  Step 1.

\end{description}

\smallskip

Step 2 requires the projection of $z^\nu - \lambda^\nu v$ onto the set $X\times \reals$, which thus needs to be relatively easy to accomplish. There are several possible convergence results for this algorithm. Regardless, the algorithm is inherently random and we need to view the observations $\{\xi^0, \xi^1, \dots\}$ obtained in Step 1 as outcomes of the independent random vectors $\{\bfxi^0, \bfxi^1, \dots\}$ distributed as $\bfxi$. Then, the iterates $\{z^\nu\}_{\nu=1}^\infty$ produced by the algorithm and their averages become random vectors as indicated by switching to bold face: $\{\bfz^\nu\}_{\nu=1}^\infty$. In the convex case, the average of the first $\bar \nu$ iterates is a near-minimizer with a specific tolerance that is guaranteed in expectation.

\begin{theorem}\label{thm:stochsubgradMethod}{\rm (SGD method).} For a random vector $\bfxi$ with support $\Xi\subset \reals^m$, a function $f:\reals^m\times\reals^n\to \reals$, and a
nonempty, closed, and convex set $X\subset\reals^n$, suppose that
\begin{enumerate}[{\rm (a)}]

\item $f(\bfxi,x)$ is integrable for all $x\in X$

\item $f(\xi,\cdot\,)$ is convex for all $\xi\in \Xi$

\item for some $\beta\in\reals$,
\[
\nsup \Big\{\Ex\big[\|v(\bfxi)\|_2^2\big]~\Big|~ v(\bfxi) \mbox{ integrable}, \,v(\xi)\in \partial_z \psi(\xi,z)~\forall \xi\in \Xi, ~z\in X \times \reals\Big\}\leq \beta^2.
\]
\end{enumerate}
If the SGD method has generated $\{\bfz^\nu, \nu = 1, 2, \dots, \bar \nu\}$
using step size $\lambda^\nu = \lambda = \|z^{0} - z^\star\|_2/(\beta \sqrt{\bar \nu +1})$ for all $\nu$, then
\[
\Ex\big[ \phi(\bar\bfz)\big] - \inf_{z\in X\times\reals} \phi(z) \leq  \frac{\beta\|z^{0} - z^\star\|_2}{\sqrt{\bar\nu + 1}}
\]
where $z^\star \in \nargmin_{z\in X\times\reals} \phi(z)$ and $\bar \bfz = \frac{1}{\bar \nu + 1}(\nsum_{\nu = 1}^{\bar \nu} \bfz^\nu + z^0 )$.
\end{theorem}

A slight adjustment of the arguments in Sections 2.1 and 2.2 of \cite{NemirovskiJuditskyLanShapiro.09} leads to the theorem; see \cite[Section 3.G]{primer} for details. These arguments stem originally from \cite{NemirovskiYudin.78,NemirovskiYudin.83}. With the importance of expectation minimization in machine learning, there is a rapid development of closely related algorithms under milder assumptions, which we do not attempt to review systematically; see \cite{DavisDrusvyatskiy.18,DuchiRuan.18} for the setting of weakly convex functions, \cite{Ruszczynski.20} for functions that are differentiable in the generalized sense of V. Norkin, and \cite{ZhangLinJegelkaSraJadbabaie.20,DavisDrusvyatskiy.21} for modifications of Goldstein's subgradient method. There is also an extensive literature that only relies on function evaluations; see, e.g., \cite{LinZhengJordan.22}. The monograph \cite{Lan.20} provides an in-depth treatment, while the overview article \cite{GarrigosGower.23} offers accessible proofs. For estimates of the actual minimum value in the context of stochastic subgradient type algorithms, we refer to \cite{LanNemirovskiShapiro.12}. Theoretical justification for adapting the SGD method to \eqref{eqn:sriskmingamma}, with some  empirical evidence, appears in \cite{SomaYoshida.20}.  

While \eqref{eqn:sriskmingamma} indeed is just an expectation minimization problem over $z = (x,\gamma)$, it has a rather specific structure which may cause difficulties \cite{CuriLevyJegelkaKrause.20}. For example, if $\alpha$ is near one, then the optimal $\gamma$ for any fixed $x$ tends to be high; in fact it can be taken as the $\alpha$-quantile of $f(\bfxi,x)$; see Theorem \ref{tsuperquantile}. Thus, there will be ``few'' outcomes of $f(\bfxi,x)$ that exceed $\gamma$. A (small) minibatch $\{\xi^1, \dots, \xi^\nu\}$ may produce $\max\{0, f(\xi^i,x) -\gamma\} =0$ for all $i = 1, \dots, \nu$ and then a zero (sub)gradient with respect to $x$. This can be addressed, in part, by keeping $\gamma$ relatively low and/or increasing the minibatch size \cite{Yang.21}. Regardless, the main advantage of minimizing superquantiles via \eqref{eqn:sriskmingamma} is the possibility of leveraging the vast computational infrastructure for expectation minimization.

Another possibility, viable when the support of $\bfxi$ has low cardinality (say, less than 5,000), is to minimize a superquantile directly without passing through the formula \eqref{eqn:superquantRU}. This would require a subgradient of $x\mapsto \srsk_\alpha(f(\bfxi,x))$. Proposition 2 of \cite{LaguelPillutlaMalickHarchaoui.21} furnishes the following convenient expression.

\begin{proposition}{\rm (subgradients of superquantile functions).}\label{pSubgradsuper}
For $\alpha \in (0,1)$, $f:\reals^m\times \reals^n\to\reals$, and a random vector $\bfxi$ with finite support $\Xi\subset\reals^m$ of cardinality $s$ and associated probabilities $p_\xi = 1/s$ for all $\xi\in \Xi$, suppose that $f(\xi,\cdot\,)$ is continuously differentiable for all $\xi\in \Xi$. Consider the function given by
\[
  \phi(x) = \srsk_\alpha\big(f(\bfxi,x)\big).
\]
Then, then the set of subgradient\footnote{Subgradients are defined in a general (Mordukhovich) sense; see, e.g., \cite[Section 4.I]{primer}.} of $\phi$ at $x\in \reals^n$ is
\[
  \partial\phi(x) = \frac{1}{s(1-\alpha)}\sum_{\xi\in \Xi_+(x)} \nabla_x f(\xi,x) + \frac{\sigma - \alpha}{1-\alpha}\con\big\{\nabla_x f(\xi,x), ~\xi \in \Xi_0(x)\big\},
\]
where $\sigma = 1 - |\Xi_+(x)|/s$, $\con A$ is the convex hull of the set $A$, and
\[
\Xi_+(x) = \big\{\xi\in \Xi~\big|~f(\xi,x)>Q(\alpha)\big\} ~~\mbox{ and } ~~\Xi_0(x) = \big\{\xi\in \Xi~\big|~f(\xi,x)=Q(\alpha)\big\},
\]
with $Q(\alpha)$ being the $\alpha$-quantile of $f(\bfxi,x)$ and $|\Xi_+(x)|$ being the cardinality of $\Xi_+(x)$.
\end{proposition}

At least in the convex case, the proposition furnishes the necessary ingredient for a subgradient method (see, e.g., \cite[Subsection 2.I]{primer}) to minimize a superquantile. If $f(\xi,\cdot\,)$ is nonconvex but smooth with Lipschitz continuous gradients $\nabla_x f(\xi,\cdot\,)$ for each $\xi$ in the finite set $\Xi$, then $\phi$ is weakly convex and one can leverage developments in \cite{DavisDrusvyatskiyMacPheePaquette.18}; see also references therein. The calculation of a subgradient requires a quantile of $f(\bfxi,x)$, which essentially amounts to sorting of the $s$ values $f(\xi,x), \xi\in \Xi$. Thus, the worst-case complexity of computing a subgradient is of order $O(s \ln s)$. This means that a subgradient method for minimizing a superquantile has a relatively high, per-iteration computational cost if $s$ is large. For a more detailed complexity analysis and comparisons, we refer to \cite{LevyCarmonDuchiSidford.20} and references therein.\\

\noindent {\bf Primal Smoothing Methods.} Even if the quantity of interest $f(\xi,x)$ is continuously differentiable in $x$ for every $\xi \in \Xi$, the function $x\mapsto \srsk_\alpha(f(\bfxi,x))$ is only exceptionally smooth. (The cases $\alpha = 0$ and $\Xi$ being a singleton are examples of such exceptions.) The situation is similar as we pass to \eqref{eqn:sriskmingamma} after invoking the formula \eqref{eqn:superquantRU}; the max-term is typically nonsmooth. The simple form of this term, however, makes it easy to approximate using a continuously differentiable function; see \cite{AlexanderColemanLi.06,SomaYoshida.20} for a quadratic approximation and
\cite{BasovaRockafellarRoyset.11,KouriSurowiec.16,Yang.21} for exponential smoothing. The latter type of smoothing has the advantage of preserving any order of smoothness in $f(\xi,\cdot\,)$. Concretely, exponential smoothing (see, e.g., \cite[Example 4.16]{primer}) amounts to replacing $\gamma\mapsto \max\{0,\gamma\}$ by the function $h_\theta:\reals\to \reals$ defined as
\[
h_\theta(\gamma) = \frac{1}{\theta} \ln \big(1 + \exp(\theta \gamma)\big),
\]
where $\theta>0$ is a parameter. The error caused by smoothing is controlled by
\[
0 \leq h_\theta(\gamma) - \max\{0,\gamma\} \leq \frac{\ln 2}{\theta}~~~~\forall \gamma\in \reals.
\]
The function $h_\theta$ is differentiable any number of times, with easily accessible derivative formulas; see \cite[Example 4.16]{primer}. Thus, there is strong theoretical backing for approximating \eqref{eqn:sriskmingamma} by
\begin{equation}\label{eqn:smoothsrsk}
\nnmin_{x\in X,\gamma\in\reals} ~\Ex\Big[\gamma + \frac{1}{1-\alpha}h_\theta\big(f(\bfxi,x)-\gamma\big)\Big],
\end{equation}
which involves a continuously differentiable objective function provided that $f(\xi,\cdot\,)$ is continuously differentiable and mild additional assumptions hold; see, e.g., \cite[Subsection 9.2.5]{ShapiroDentchevaRuszczynski.21}. We refer to \cite{Royset.22} for general tools to justify such approximations and to \cite{Chen.12,KouriSurowiec.16} for many related smoothing schemes. Regardless, we are faced with the expectation minimization problem \eqref{eqn:smoothsrsk} and can leverage a vast array of existing algorithm, which may perform better practically and theoretically when applied to a smooth problem than to the (potentially) nonsmooth problem \eqref{eqn:sriskmingamma}.\\

\noindent {\bf Proximal Composite Method.} Minimizing a superquantile of $f(\bfxi,x)$ over $x \in X$ when $\bfxi$ is finitely distributed with support $\Xi$ is equivalent to \eqref{eqn:superquantExpanded0} as achieved by the formula \eqref{eqn:superquantRU}. This reformulation is well structured because it can be written in terms of the real-valued convex function $h:\reals^{1+s}\to \reals$ given by
\[
h(\gamma, u) = \gamma + \frac{1}{1-\alpha} \nsum_{\xi \in \Xi} \max\{0, u_\xi - \gamma\}, ~~~u = (u_\xi, \xi\in \Xi),
\]
and the mapping $F:\reals^{1+n}\to \reals^{1+s}$, where $F(\gamma,x) = (\gamma, (f(\xi,x), \xi\in \Xi))$; $s$ is the cardinality of $\Xi$. Thus, \eqref{eqn:superquantExpanded0} is equivalently expressed as minimizing $h(F(\gamma,x))$ over $\gamma\in \reals$ and $x\in X$. The main advantage of this perspective is that the potentially difficult functions $f(\xi, \cdot\,)$, $\xi\in \Xi$, are separated from the remaining parts of the problem represented by $h$ and $X$. If $f(\xi, \cdot\,)$ is smooth, then we may linearize it and potentially avoid evaluating it a large number of times. The following algorithm and its convergence result is a modification of \cite{LewisWright.16}; see \cite[Section 6.F]{primer} for details.\\

\state Proximal Composite Method.

\begin{description}

  \item[Data.] ~~$z^0\in \reals\times X$, $\tau\in (1,\infty), \sigma\in (0,1)$, $\bar \lambda \in (0,\infty)$, $\lambda^0 \in (0, \bar \lambda]$.

  \item[Step 0.]  Set $\nu = 0$.

  \item[Step 1.]  Compute
\[
  \bar z^\nu \in \nargmin_{z\in \reals\times X} \Big\{h\big(F(z^\nu) + \nabla F(z^\nu)(z-z^\nu)\big) + \proxfacnu\|z-z^\nu\|_2^2\Big\}.
\]
~~~~~If $\bar z^{\nu}=z^\nu$, then Stop.

  \item[Step 2.]  If
  \[
  h\big(F(z^\nu)  \big) - h\big(F(\bar z^\nu)  \big) \geq \sigma \Big( h\big(F(z^\nu)\big) - h\big(F(z^\nu) + \nabla F(z^\nu)(\bar z^\nu - z^\nu) \big) \Big),
  \]
~~~~~then set $\lambda^{\nu+1} = \min\{\tau \lambda^\nu, \bar\lambda\}$ and go to Step 3.

~~~~~Else, replace $\lambda^\nu$ by $\lambda^\nu/\tau$ and go to Step 1.

  \item[Step 3.] Set $z^{\nu+1} = \bar z^\nu$, replace $\nu$ by $\nu +1$, and go to  Step 1.
\end{description}

\begin{theorem}{\rm (proximal composite method).}\label{cProxComposite} For closed convex $X\subset\reals^n$ and twice continuously differentiable $f(\xi,\cdot\,)$, $\xi\in \Xi$, suppose that the proximal composite method has generated $\{z^\nu\}_{\nu=1}^\infty$ with a cluster point $(\gamma^\star,x^\star)$. Then, $(\gamma^\star,x^\star)$ satisfies a necessary optimality condition\footnote{The normal cone to a set $X$ at a point $x$ is denoted by $N_X(x)$; see, e.g., \cite[Section 4.G]{primer}.} for \eqref{eqn:superquantExpanded0}:
\[
\exists y = (y_\xi, \xi\in \Xi)\in \reals^s  ~\mbox{ such that }~ (0,y) \in \partial h\big(F(\gamma^\star,x^\star)\big) \,\mbox{ and } \,-\nsum_{\xi\in \Xi} y_\xi \nabla_x f(\xi,x^\star) \in N_X(x^\star).
\]
\end{theorem}

The advantage of this approach is that any convex optimization algorithm can be brought in to solve the subproblem in Step 1. If $X$ is polyhedral, then one may even solve the subproblem as a convex quadratic problem after introducing auxiliary variables. Step 2 tests whether the present linear approximation is sufficiently accurate and decreases $\lambda^\nu$ if it is not.\\

\noindent {\bf Dual Algorithms.} The dual formula \eqref{eqn:maxsuper} gives rise to the following approach. As recognized in  \cite{LaguelMalickHarchaoui.20,LaguelPillutlaMalickHarchaoui.21,PillutlaLaguelMalickHarchaoui.22}, we find that  
\[
\srsk_\alpha\big(f(\bfxi,x)\big) = \max_{\bar p\in \Delta_\alpha} \nsum_{\xi\in \Xi} \bar p_\xi f(\xi,x) \approx \max_{\bar p\in \Delta_\alpha} \nsum_{\xi\in \Xi} \bar p_\xi f(\xi,x) - \epsilon \phi(\bar p),
\]
where $\epsilon > 0$ and $\phi$ is a nonnegative strongly convex real-valued function. Proposition 3 of \cite{LaguelPillutlaMalickHarchaoui.21} asserts that the approximation error vanishes as $\epsilon \to 0$ and the approximation is continuously differentiable if $f(\xi,\cdot\,)$ is continuously differentiable for all $\xi\in \Xi$. Thus, one can minimize a superquantile approximately by minimizing a smooth approximation. A challenge is that each function and gradient computation of the approximating objective function requires the full set $\Xi$, which might be of high cardinality, especially in learning applications. We refer to \cite{PillutlaLaguelMalickHarchaoui.22} for further discussion of such issues. This dual smoothing approach can also be linked to primal smoothing methods; see \cite[Corollary 6]{LaguelPillutlaMalickHarchaoui.21}. A related possibility is mentioned in \cite[Example 8.32]{primer}.\\

\noindent {\bf Other Algorithmic Approaches.} By specializing the classical L-shaped method (see, e.g., \cite[Section 5.H]{primer}), the paper \cite{KunzibayMayer.06} achieves a decomposition algorithm for minimizing superquantiles involving quantities of interest that are affine in $x$; see \cite{Fabian.08,FabianWolfKobersteinSuhl.15} for further improvements. In the limited setting of simplex constraints and affine quantities of interest, \cite{LimSheraliUryasev.10} develops a three-phase algorithm that starts with a gradient descent heuristic applied to \eqref{eqn:superquantExpanded0}, continues with steps akin to the subgradient method, and ends with solving \eqref{eqn:superquantExpanded} from (hopefully) a good starting point using the simplex method. Simple active-set strategies to reduce the number of constraints in \eqref{eqn:superquantExpanded} generally appear highly beneficial \cite{BasovaRockafellarRoyset.11,ByunRoyset.22}.

Essentially all the algorithms discussed in this subsection require that we can evaluate the quantity of interest $f(\xi,x)$ as well as its gradients or subgradients with respect to $x$. If only function values are available, one can resort to black-box optimization \cite{AudetHare.17} and other zeroth-order methods, e.g., \cite{LinZhengJordan.22}. If function values are computationally costly to compute (possibly even with noise), then one might resort to surrogates; cf. Example \ref{eSurrogate}.

\subsection{Estimating Superquantiles}

Except when a random variable $\bfxi$ is finitely distributed or follows some other standard distribution (cf. Examples \ref{eSuperTriang}, \ref{eSuperNormal}, and \ref{eSuperFinite} as well as \cite{NortonKhokhlovUryasev.21}), there is no explicit formula for its $\alpha$-superquantile $\bar Q(\alpha)$. If $\bfxi$ is square integrable, then we have the bound from \cite[Proposition 1]{RockafellarRoysetMiranda.14} for any $\alpha \in [0,1)$:
\begin{equation}\label{eqn:barGbound}
  \Ex[\bfxi] \leq \bar Q(\alpha)  \leq \min\left\{\Ex[\bfxi] +
  \frac{\std(\bfxi)}{\sqrt{1-\alpha}}, ~\sup \bfxi\right\}.
\end{equation}

Much more accurate estimates of superquantiles are available via sampling. Suppose that $\bfxi_1, \bfxi_2, \dots$ are independent random variables with the same distribution as $\bfxi$, which is assumed to be integrable. Then, for $\alpha \in (0,1)$, the estimator
\begin{equation}\label{eqn:sqestimator}
{\bar{\bf Q}}^\nu(\alpha) = \min_{\gamma\in \reals} \, \gamma + \frac{1}{\nu(1-\alpha)} \sum_{i = 1}^\nu \max\{0, \bfxi_i - \gamma\}
\end{equation}
is strongly consistent, i.e., ${\bar{\bf Q}}^\nu(\alpha) \to \bar Q(\alpha)$ as $\nu\to \infty$ almost surely (\cite[Section 6.5.1]{ShapiroDentchevaRuszczynski.09} and \cite{WozabalWozabal.09}); see also \cite[Example 8.57]{primer}. As that example shows, the consistency carries over in the sense of epi-convergence to the case when $\bfxi$ is replaced by a quantity of interest $f(\xi,x)$ that is continuous in $x$. Thus, passing from \eqref{eqn:sriskmingamma} to \eqref{eqn:superquantExpanded0} by generating a finite sample independently is fundamentally sound, with cluster points of minimizers of the sample average approximation \eqref{eqn:superquantExpanded0} being minimizers of the actual problem \eqref{eqn:sriskmingamma} almost surely. For superquantiles and many other risk measures, \cite{Shapiro.13} furnishes a comprehensive treatment of strong consistency and \cite{GuiguesKratschmerShapiro.18} develops asymptotics for minimum values and minimizers obtained through solving sample average approximations as well as associated hypothesis testing methodology. Asymptotics for ${\bar{\bf Q}}^\nu(\alpha)$ appeared already in \cite[Section 6.5.1]{ShapiroDentchevaRuszczynski.09}.

There is a wealth of concentration inequalities for bounding the probability that a sample average deviates from the mean with some $\epsilon$ for a fixed sample size. Hoeffding's inequality addresses bounded random variables, but extensions to subgaussian random variables and beyond are available; see, for example, \cite{BoucheronLugosiMassart.13}. These concentration inequalities can be made to hold uniformly in some sense across values of the auxiliary variable $\gamma$ in the optimization problems giving ${\bar{\bf Q}}^\nu(\alpha)$ and $\bar Q(\alpha)$; see \cite[Section 9.2.11]{ShapiroDentchevaRuszczynski.21} for a general discussion. An early effort to capitalize on these possibilities is \cite{Brown.07}, which assumes that the support $\Xi$ of $\bfxi$ is a bounded interval $[0,\beta]$. Then, one has for every $\epsilon \in (0,\infty)$ and $\alpha \in (0,1)$ that
\begin{equation}\label{eqn:Brown1}
\prob\Big\{ {\bar{\bf Q}}^\nu(\alpha) \geq \bar Q(\alpha) + \epsilon  \Big\} \leq \exp\big(-2 (\alpha\epsilon/\beta)^2 \nu \big)
\end{equation}
and, provided that $\bfxi$ is continuously distributed, one also has
\[
\prob\Big\{ {\bar{\bf Q}}^\nu(\alpha)\leq \bar Q(\alpha) - \epsilon  \Big\} \leq 3\exp\big(-\alpha (\epsilon/\beta)^2 \nu /5 \big).
\]
The square dependence of $\alpha$ in \eqref{eqn:Brown1} is improved to linear dependence in \cite{WangGao.10}; see also \cite{GaoWang.11} for a discussion of Berry–Essen bounds, the law of iterated logarithm, and large deviation results. Recent efforts for bounded random variables also include \cite{Akellaetal.22}.
For unbounded subgaussian continuously distributed random variables, \cite{Kumaretal.19} achieves concentration bounds as long as the cumulative distribution functions are strictly increasing; see also \cite{BathPrashanth.19} for a slight tightening and \cite{Kumaretal.20} for a treatment of random variables with only finite $p$th moment for $p \in (1,2]$. Deviating from \eqref{eqn:sqestimator}, \cite{ThomasLearnedmiller.19} considers a more complicated estimator and achieves a strict improvement compared to \cite{Brown.07}, with strong empirical performance, but still in the setting of bounded random variables. The high-water mark for \eqref{eqn:sqestimator} appears to be the recent PAC-Bayesian bound in \cite{MhammediGuedjWilliamson.20}. The paper \cite{HongHuLiu.14} reviews developments up to 2014. Uniform bounds holding across a set of decisions appear in \cite{LeeParkShin.20,KhimLeqiPrasadRavikumar.20}. The paper \cite{GlynnFanFuHuPeng.20} examines central limit theorems for superquantiles and related conditional expectations.

Going beyond independent samples, \cite{LuoOu.17} establishes strong laws of large numbers and associated convergence rates for $\alpha$-mixing\footnote{The $\alpha$ here is of course unrelated to the $\alpha$ specifying which superquantile is being considered.} sequences and \cite{Luo.20} considers a kernel estimator also for $\alpha$-mixing sequences. In the context of computationally expensive models of physical systems, \cite{HeinkenschlossKramerTakhtaganovWillcox.18,HeinkenschlossKramerTakhtaganov.20} discuss reduced-order models to guide a choice of importance sampling and reduce the sample size needed in \eqref{eqn:sqestimator}. Similarly, \cite{ChaudhuriPeherstorferWillcox.20} explores variance reduction via the cross-entropy method, \cite{GilesHajiali.19} examines efficient multi-level nested Monte Carlo simulations, \cite{LeeKramer.23a,LeeKramer.23b} leverage multi-fidelity simulations and polynomial chaos expansion, and \cite{ChenKim.16} uses metamodels based on stochastic kriging with further extensions in \cite{DylanFredericYuan.19,KhayyerVinelKennedy.24} that also bring in extreme value theory. The paper \cite{FairbrotherTurnerWallace.22} discusses how to construct a finite distribution that results in an accurate approximation of an actual problem involving superquantile minimization; see also \cite{ArponHomemdemelloPagnoncelli.18} for other scenario reduction techniques.

\section{Risk and Regret}\label{sec:regretrisk}

As seen in Subsection \ref{subsec:measureofrisk}, there are many possible measures of risk with superquantiles furnishing main examples. In this section, we discuss what constitutes a ``good'' measure of risk, both from modeling and computing points of view. We introduce measures of regret as key building blocks for constructing measures of risk and discuss how superquantiles generate many, if not most, meaningful measures of risk. The section ends with a review of the role risk measures play in achieving {\em fairness}.

\subsection{Desirable Properties}\label{subsec:desirable}

The purpose of a risk measure is to model risk-averseness in a meaningful way and especially avoid embarrassing paradoxes that might discredit the modeling effort in the eyes of a decision maker. We would like to achieve this without introducing excessive computational complexity. So what properties would we like a measure of risk to possess? To formally answer this question, we view risk measures as functionals defined on a space of random variables.

For a given probability space $(\Omega, \cA, \bbP)$, we consider all square-integrable random variables:
\[
\cL^2 = \big\{\bfxi:\Omega\to\reals~\big|~\Ex\big[\bfxi^2\big]<\infty\big\}.
\]
As seen in \cite{RuszczynskiShapiro.06,CheriditoLi.08,FollmerSchied.16,DelbaenOwari.19}, one may consider alternative spaces, with the choice possibly depending on the risk measures of interest. For a  detailed discussion of suitable spaces on which to define risk measures, see the papers \cite{Pichler.13,FrohlichWilliamson.22}.  In this survey, we follow \cite{RockafellarUryasev.13} and limit the scope to $\cL^2$, which simplifies the exposition significantly.

We equip $\cL^2$ with the standard {\em norm} $\|\bfxi\|_{\cL^2} = (\Ex[\bfxi^2])^{1/2}$.
For $\{\bfxi,\bfxi^\nu\in \cL^2, \nu=1, 2, \dots\}$, the random variables $\bfxi^\nu$ {\em converge} to the random variable $\bfxi$, written $\bfxi^\nu\to \bfxi$, when $\|\bfxi^\nu - \bfxi\|_{\cL^2}\to 0$ or, equivalently, when $\Ex[(\bfxi^\nu - \bfxi)^2] \to 0$.  A subset $\cC\subset\cL^2$ is {\em closed} if $\{\bfxi^\nu\in \cC, \nu=1, 2, \dots\}$ and $\bfxi^\nu\to \bfxi$ imply $\bfxi\in \cC$. The set $\cC$ is {\em convex} if $(1-\lambda)\bfxi_0 + \lambda \bfxi_1 \in \cC$ for all $\lambda \in [0,1]$ and $\bfxi_0,\bfxi_1\in \cC$. We recall that a random variable $\bfxi\in \cL^2$ is {\em constant} if $\prob\{\bfxi = \alpha\} = \bbP(\{\omega\in \Omega~|~\bfxi(\omega) = \alpha\}) = 1$ for some $\alpha \in \reals$. The constant random variable with value $0$ is denoted by $\bfnull$. The underlying probability space $(\Omega, \cA, \bbP)$ is {\em finite} if the cardinality of $\Omega$ is finite, which means that any random variable in $\cL^2$ has a finite number of possible values, i.e., is finitely distributed. Throughout, we adopt the usual extended real-valued arithmetic rules such as $0 \cdot \infty = 0$, $0 \cdot (-\infty) = 0$, $-\infty + \infty = \infty$, and $\infty - \infty = \infty$; see, e.g., \cite[Section 1.E]{VaAn} or \cite[Section 1.D]{primer}. 

A functional $\cF:\cL^2\to [-\infty,\infty]$ may satisfy any of the following properties:
\begin{align*}
&\mbox{{\em  Constancy:} } &&\cF(\bfxi) = \alpha \mbox{ whenever $\bfxi$ is constant with value $\alpha\in\reals$}.\\
& \mbox{{\em Averseness:} } &&\cF(\bfxi) > \Ex[\bfxi] \mbox{ for all nonconstant $\bfxi$}.\\
& \mbox{{\em Convexity:} } &&\cF\big((1-\lambda)\bfxi_0+\lambda \bfxi_1\big) \leq (1-\lambda)\cF(\bfxi_0)+\lambda\cF(\bfxi_1) ~~\forall\bfxi_0,\bfxi_1, ~\lambda\in [0,1].\\
&\mbox{{\em  Lower semicontinuity:}} &&\{ \cF \leq \alpha \} = \big\{\bfxi\in \cL^2~\big|~\cF(\bfxi) \leq \alpha \big\} \mbox{ is closed}~~ \forall \alpha\in\reals.\\
& \mbox{{\em Positive homogeneity:} } &&\cF(\lambda\bfxi) = \lambda \cF(\bfxi) ~~~\forall \lambda \in [0,\infty).\\
&\mbox{{\em Monotonicity:} } &&\cF(\bfxi_0) \leq \cF(\bfxi_1) ~\mbox{ when }~ \bfxi_0(\omega) \leq \bfxi_1(\omega) ~\mbox{ for } \bbP\mbox{-almost every } \omega\in \Omega.\\
&\mbox{{\em Law invariance:} } &&\cF(\bfxi_0) = \cF(\bfxi_1) ~\mbox{ when }~ \bfxi_0, \bfxi_1 \mbox{ have the same distribution}.
\end{align*}

The constancy requirement is certainly reasonable and in fact a prerequisite for a measure of risk; see Definition \ref{dRiskmeasure}. The averseness property excludes the possibility $\cF(\bfxi) = \Ex[\bfxi]$, which is better treated separately. The purpose of a risk measure is after all to model risk-averseness.

The convexity property is key to make sure that a risk measure is computationally attractive. As we see from the comparison between quantiles and superquantiles in Example \ref{eQuantvsSuperquant}, it is also important in promoting hedging. Generally, the measure of risk $\cR(\bfxi) = Q(\alpha)$, where $Q(\alpha)$ is the $\alpha$-quantile of $\bfxi$, fails the convexity requirement. For example, consider two statistically independent random variables $\bfxi_0,\bfxi_1$ with common density function value $0.9$ on $[-1,0]$, value $0.05$ on $(0,2]$, and value zero otherwise. Then, the $0.9$-quantile is $0$ for both $\bfxi_0$ and $\bfxi_1$, but the $0.9$-quantile for $\bfxi_0/2 + \bfxi_1/2$ is approximately $0.25$.

Lower semicontinuity (lsc) is a technical condition beneficial in convex analysis (as seen in Subsection \ref{sec:duality}). We know even from finite-dimensional optimization that the minimum of a function over a compact set may not be attained if the function is not lsc. The lsc property holds for example when $\cF$ is real-valued and convex and either $\cF$ is monotone or the probability space is finite. In fact, then $\cF$ is continuous; see \cite[Equation (3.5)]{RockafellarUryasev.13} for a summary of this claim which in turn relies on the fundamental Proposition 3.1 about continuity of convex and monotone functionals in \cite{RuszczynskiShapiro.06}.

Positive homogeneity implies a certain invariance to scaling. If we convert the quantity of interest from dollar to euro, then the associated risk should not fundamentally change. However, this property fails for the variance $\cF(\bfxi) = \var(\bfxi) = (\std(\bfxi))^2$.

Monotonicity is a natural requirement because risk would typically be deemed less for $\bfxi_0$ than for $\bfxi_1$ when, for every pair of outcomes $\{\bfxi_0(\omega),\bfxi_1(\omega)\}$, the random variable $\bfxi_0$ never comes out worse than $\bfxi_1$. The mean-plus-standard-deviation risk measure $\cR(\bfxi) = \Ex[\bfxi] + \lambda \std(\bfxi)$ fails this requirement because the standard deviation does not distinguish between variability above and below the mean. Thus, a random variable with a high variability below the mean is deemed ``high risk,'' even though low values represent no ``real'' risk under our orientation. For example, consider the random variables $\bfxi_0$ and $\bfxi_1$ and their joint distribution that assigns probability $1/2$ to the outcome $(0,0)$ and probability $1/2$ to the outcome $(-2, -1)$. Thus, for every outcome the random variables have  either the same value or $\bfxi_0$ has a value below that of $\bfxi_1$. Still, with $\lambda = 2$, we obtain $\Ex[\bfxi_0] + \lambda \std(\bfxi_0) = 1$ and $\Ex[\bfxi_1] + \lambda \std(\bfxi_1) = 1/2$. 

In elementary probability theory, random variables are thought of as being ``fully'' described by their distributions. However, random variables are (measurable) functions from a sample space (here denoted by $\Omega$) to the reals, with many random variables potentially having the same distribution. Thus, the law invariance property is not automatically satisfied but certainly appears reasonably in many situations.

For $\alpha \in (0,1]$, the superquantile risk measure $\srsk_\alpha$ satisfies all the seven properties above \cite{RockafellarUryasev.13}; $\srsk_0 = \Ex[\cdot]$ misses the averseness requirement. It is also real-valued when $\alpha \in [0,1)$ by \eqref{eqn:barGbound}.

Constancy, convexity, and lsc properties together imply that 
\begin{equation*}\label{eqn:constancy}
\cF(\bfxi + \alpha) = \cF(\bfxi) + \alpha ~~~~\forall \bfxi\in \cL^2, ~\alpha\in\reals.
\end{equation*}
Thus, $\cF$ quantifies, for example, risk in a manner that is translation invariant; adding a constant amount to an uncertain future cost changes the risk by that amount.

Pioneering works in finance \cite{KijimaOhnishi.93,ArtznerDelbaenEberHeath.97,ArtznerDelbaenEberHeath.99,Delbaen.02} identify many of these properties as desirable. Although originally expressed slightly differently, \cite{ArtznerDelbaenEberHeath.99} defines a real-valued risk measure to be {\em coherent} if it satisfies the constancy, convexity, positive homogeneity, and monotonicity properties. The initial development in  \cite{ArtznerDelbaenEberHeath.99} was for finite probability spaces, with an extension to spaces of bounded random variables in \cite{Delbaen.02}.
The coherency of the quantity in \eqref{eqn:superquantRU} (which we now of course call a superquantile) was first established in \cite{Pflug.00} and, as the complete picture about the equivalences in Theorem \ref{tsuperquantile} emerged, also in \cite{AcerbiTasche.02,RockafellarUryasev.02}.

Further developments under the names {\em convex measures of risk} and {\em convex risk functions} relaxed the positive homogeneity condition \cite{FollmerSchied.02b,RuszczynskiShapiro.06}.
We adopt the concept of  {\em regularity} from \cite{RockafellarUryasev.13} (with refinements in \cite{RockafellarRoyset.15}), which is also used in \cite{KouriSurowiec.20}.

\begin{definition}{\rm (regular measure of risk).}\label{dRegMeasRegretRisk} A {\em regular measure of risk} $\cR$ is a functional from $\cL^2$ to $(-\infty,\infty]$ that is lsc, convex and also satisfies the constancy and averseness properties.
\end{definition}

As we see from the following development, these requirements align well with convex analysis while permitting a wide set of risk measures including those that may assign $\infty$ to some random variables and that may lack positive homogeneity and monotonicity. 

While measures of risk take a central role, there are other classes of functionals as well. In Subsection \ref{subsec:utility}, we discuss utility and disutility functions. It turns out that particular disutility functions, which can be traced back to studies of optimized certainty equivalents \cite{BentalTeboulle.86} (see also \cite{BentalTeboulle.07}), emerge as important building blocks for regular measures of risk. Here, we adopt the definition in \cite{RockafellarRoyset.15}.

\begin{definition}{\rm (regular measure of regret).}\label{dRegMeasRegret} A {\em regular measure of regret} $\cV$ is a functional from $\cL^2$ to $(-\infty,\infty]$ that is lsc and convex, with
\[
\cV(\bfnull) = 0, ~\mbox{ but }~\cV(\bfxi) > \Ex[\bfxi] ~~\forall \bfxi\neq \bf0.
\]
The quantity $\cV(\bfxi)$ is the {\em regret} of $\bfxi$.
\end{definition}

We note that the definition of a regular measure of regret in \cite{RockafellarUryasev.13} includes a limiting condition that is shown to be superfluous in \cite{RockafellarRoyset.15}.

\begin{example}{\rm (measures of regret and risk).}\label{eMeasRegRisk}  We have the following examples of regular measures of regret $\cV$ and regular measures of risk $\cR$; see, e.g., \cite[Examples 8.8 and 8.10]{primer}:
\begin{enumerate}[{\rm (a)}]
\item {\em Penalty regret} and {\em superquantile risk} with $\alpha \in (0,1)$:
\begin{flalign*}
& \cV(\bfxi) = \frac{1}{1-\alpha}\Ex\big[\max\{0, \bfxi\}\big] &  \cR(\bfxi) = \srsk_\alpha(\bfxi).
\end{flalign*}

\item {\em Worst-case regret} and {\em risk}:
\begin{flalign*}
&\cV(\bfxi) = \begin{cases}
0 & \mbox{ if } \nsup \bfxi \leq 0\\
\infty & \mbox{ otherwise}
\end{cases}
& \cR(\bfxi) = \sup \bfxi.
\end{flalign*}

\item {\em Moment-based regret} and {\em mean-plus-standard-deviation risk} with $\lambda>0$:
\begin{flalign*}
& \cV(\bfxi) = \Ex[\bfxi] + \lambda \sqrt{\Ex[\bfxi^2]} & \cR(\bfxi) = \Ex[\bfxi] +\lambda \std(\bfxi).
\end{flalign*}

\item {\em Mixed regret} and {\em risk} with $\lambda_1, \dots, \lambda_q>0$ and $\nsum_{i=1}^q \lambda_i = 1$:
\begin{flalign*}
& \cV(\bfxi) = \inf_{\gamma_1, \dots, \gamma_q} \Big\{\nsum_{i=1}^q \lambda_i \cV_i(\bfxi-\gamma_i)~\Big|\nsum_{i=1}^q \lambda_i \gamma_i = 0 \Big\} & \cR(\bfxi) = \nsum_{i=1}^q \lambda_i \cR_i(\bfxi)\\
& \mbox{when $\cV_1, \dots, \cV_q$ ($\cR_1, \dots, \cR_q$) are regular measures of regret (risk)}.
\end{flalign*}
\end{enumerate}
\end{example}

Several additional examples appear in Section \ref{sec:addexamples}. The distinction between regular measures of regret and risk emerges from the assessment of constant random variables. While $\cV(\bfnull) = \cR(\bfnull) = 0$, a constant random variable $\bfxi$ with value $\alpha\neq 0$ is treated differently:
\[
\cV(\bfxi) > \alpha, ~\mbox{ but } ~\cR(\bfxi) = \alpha.
\]
This seemingly minor discrepancy results in profoundly different roles, with regular measures of regret providing a means to construct regular measures of risk and thereby extending the useful
relationship between $\frac{1}{1-\alpha}\Ex[\max\{0,\bfxi\}]$ and
$\srsk_\alpha(\bfxi)$ in the formula \eqref{eqn:superquantRU} for superquantiles.

\begin{theorem}\label{tRiskregret} {\rm (regret-risk).} For a regular measure of regret $\cV$,
a regular measure of risk $\cR$ emerges from 
\[
  \cR(\bfxi) = \min_{\gamma\in\reals} \,\gamma + \cV(\bfxi - \gamma).
\]
For every regular measure of risk $\cR$ there exists a regular measure of regret $\cV$, not
necessarily unique, that constructs $\cR$ through this minimization formula.
\end{theorem}
\state Proof. Under slightly more restrictive assumptions, \cite{RockafellarUryasev.13} establishes the first conclusion. The present form is taken from \cite[Theorem 2.2]{RockafellarRoyset.15}; see \cite[Theorem 8.9]{primer} for the existence claim.\eop

The theorem generalizes the situation for superquantiles as established by the formula \eqref{eqn:superquantRU}: a regular measure of risk finds a ``best'' way of reducing displeasure with a mix of uncertain outcomes by optimizing a scalar $\gamma$. Thus, a regular measure of risk provides a deeper assessment of a random variable compared to a regular measure of regret. We note that by writing ``min'' instead of ``inf'' we assert that the minimum value in the theorem is attained for some $\gamma$.

As seen from \cite[Example 8.10]{primer}, which in part relies on \cite{RockafellarUryasev.13}, the regular measures of regret and risk in Example \ref{eMeasRegRisk}(a,b,c) pair in the sense of Theorem \ref{tRiskregret}. Moreover, if $\cV_i$ and $\cR_i$ pair by satisfying the formula in Theorem \ref{tRiskregret}, then $\cV$ and $\cR$ in Example \ref{eMeasRegRisk}(d) also satisfy that formula.

Generally, Theorem \ref{tRiskregret} provides a path to constructing new regular measures of risk from regular measures of regret, which in turn might be motivated by disutility functions.

\subsection{Risk Minimization}\label{subsec:riskmin}

A regular measure of risk offers key computational advantages when implemented in minimization problems. Suppose that $f:\reals^m\times \reals^n \to \reals$ represents a quantity of interest, $\bfxi$ is an $m$-dimensional random vector, $X\subset\reals^n$, $\cR$ is a regular measure of risk, and $f(\bfxi,x) \in \cL^2$ for all $x\in \reals^n$. Then, the problem
\begin{equation}\label{eqn:riskmin}
\nnmin_{x\in X}~\phi(x) = \cR\big( f(\bfxi, x) \big)
\end{equation}
is convex provided that $X$ is convex and $f(\xi, \cdot\,)$ is affine for all $\xi\in \Xi$, the support of $\bfxi$; see, e.g., \cite[Example 8.11]{primer} for a proof. Thus, a regular measure of risk does not disrupt the convexity of the affine functions $f(\xi,\cdot\,)$.

For any regular measure of regret $\cV$ that pairs with $\cR$ via Theorem \ref{tRiskregret}, the problem \eqref{eqn:riskmin} is equivalent to
\begin{equation}\label{eqn:regmin}
\nnmin_{x\in X,\gamma\in \reals}~\gamma + \cV\big( f(\bfxi, x) - \gamma \big),
\end{equation}
which might be computationally more attractive; cf. the situation for superquantiles in \eqref{eqn:sriskmin} and \eqref{eqn:sriskmingamma}.

If the regular measure of risk in \eqref{eqn:riskmin} is also monotone, then \eqref{eqn:riskmin} is a convex problem when $X$ is convex and $f(\xi,\cdot\,)$ is convex for all $\xi\in \Xi$; see, e.g., \cite[Example 8.13]{primer}. Such monotonicity holds if $\cV$ is monotone as can be seen directly. Then, \eqref{eqn:regmin} is also a convex problem when $X$ is convex and $f(\xi,\cdot\,)$ is convex for all $\xi\in \Xi$.

The seemingly complicated situation for mixed risk $\cR$ from Example \ref{eMeasRegRisk}(d) simplifies even further. In that case, \eqref{eqn:riskmin} is equivalent to
\[
\nnmin_{x\in X, \gamma_0, \gamma_1, \dots, \gamma_q} \gamma_0 + \nsum_{i=1}^q \lambda_i \cV_i\big(f(\bfxi,x) - \gamma_0 - \gamma_i\big) ~\mbox{ subject to }~ \nsum_{i=1}^q \lambda_i \gamma_i = 0,
\]
where $\cV_1, \dots, \cV_q$ are the regular measures of regret that pair with the regular measures of risk $\cR_1, \dots, \cR_q$, which in turn produce the mixed risk $\cR$. Reformulations of this kind extend to situations with multiple quantities of interest, possibly some emerging as constraints.

Optimality conditions are important in algorithmic development and to interpret solutions of optimization problems. For a discussion of this subject, we refer to \cite{RuszczynskiShapiro.06,RockafellarUryasevZabarankin.06b,KouriSurowiec.18}; see also the subgradient formulas in Subsection \ref{subsec:riskidentifiers}. Algorithms for computing stationary points of nonconvex, nonsmooth risk-minimization problems emerge from \cite{LiuCuiPang.22}.\\

\noindent {\bf Risk Measures from Expectation Regret.} A main approach to constructing a regular measure of risk is to start with a lsc convex function $v:\reals\to (-\infty,\infty]$ that also satisfies the properties
\[
v(0) = 0, ~\mbox{ but } ~ v(\xi)> \xi ~~~~\forall \xi\in\reals\setminus\{0\}.
\]
The value $v(\xi)$ could reflect a decision maker's displeasure with an outcome $\xi$. The expectation theorem in \cite{RockafellarUryasev.13} shows that the functional on $\cL^2$ defined by $\cV(\bfxi) = \Ex[v(\bfxi)]$ is a regular measure of regret and thus defines a regular measure of risk via Theorem \ref{tRiskregret}: $\cR(\bfxi) = \min_{\gamma\in \reals} ~\gamma + \Ex[v(\bfxi - \gamma)]$. Superquantiles fall within this class with $v(\xi) = \max\{0,\xi\}/(1-\alpha)$ for $\alpha \in (0,1)$. If $v$ is not bounded from above by a quadratic function, then these regular measures of risk and regret may not be real-valued (on $\cL^2)$. Generally, one can view $v$ as a normalized disutility function; see \cite{BentalTeboulle.07,RockafellarUryasev.13,Rockafellar.20}. Minimization of regular measures of risk constructed in this manner results in expectation minimization problems after passing from \eqref{eqn:riskmin} to \eqref{eqn:regmin}. Thus, one can leverage the vast array of algorithms for such problems including the SGD method (cf. Subsection \ref{eqn:algosuper}), the broad theory of $M$-estimators \cite{HessSeri.19}, and also tackle extensions in the context of PDE-constrained optimization using local approximations of function values and gradients \cite{ZouKouriAquino.22} or low-rank tensor approximations of random fields \cite{AntilDolgovOnwunta.21}.\\

\noindent {\bf Approximations using Epi-Regularization.} As exemplified by superquantile risk and penalty regret (see  Example \ref{eMeasRegRisk}(a)), a regular measure of regret $\cV$ when applied to a quantity of interest $f(\xi,x)$ might result in a nonsmooth function $x\mapsto \cV(f(\bfxi,x))$, which is less ideal when solving \eqref{eqn:regmin}. Likewise, $x\mapsto \cR(f(\bfxi,x))$ is often nonsmooth. However, it is possible to construct approximations using epi-regularization \cite{KouriSurowiec.20}; see also \cite{BurkeHoheisel.17}. Specifically, one can replace a regular measure of risk $\cR$ by the approximation
\[
\cR_\epsilon^\Phi(\bfxi) = \inf \big\{\cR(\bfxi - \bfeta) + \epsilon \Phi(\bfeta/\epsilon)~\big|~\bfeta \in \cL^2\big\},
\]
where $\epsilon \in (0,\infty)$ and $\Phi:\cL^2\to (-\infty,\infty]$ is lsc and convex, with $\Phi(\bfxi)<\infty$ for at least one $\bfxi\in \cL^2$. Since there is much flexibility, we can choose $\Phi$ such that $\cR_\epsilon^\Phi$ has desirable properties. For example, $\Phi(\bfxi) = \half\Ex[ \bfxi^2]$ produces in the case of $\alpha$-superquantiles (i.e., $\cR = \srsk_\alpha$) with $\alpha \in (0,1)$ that
\[
\cR_\epsilon^\Phi(\bfxi) = \min_{\gamma\in \reals} ~\gamma + \frac{1}{1-\alpha} \Ex\big[v_\epsilon(\bfxi - \gamma)\big], ~~\mbox{ where } ~v_\epsilon(\xi) = \begin{cases}
  0 & \mbox{ if } \xi \leq 0\\
  \tfrac{1}{2\epsilon}\xi^2 & \mbox{ if } \xi\in \big(0, \tfrac{\epsilon}{1-\alpha}\big)\\
  \tfrac{1}{1-\alpha}\big(\xi - \tfrac{\epsilon}{2(1-\alpha)}\big) & \mbox{ otherwise}.
\end{cases}
\]
Since $v_\epsilon$ is continuously differentiable, we may approximate $\cV$ in \eqref{eqn:regmin} by the functional $\cV_\epsilon(\bfxi) = \Ex[v_\epsilon(\bfxi)]$ and reap computational benefits. The approximation error will be small when $\epsilon$ is small because then $v_\epsilon$ accurately approximates $\xi\mapsto \max\{0, \xi\}/(1-\alpha)$; see \cite{KouriSurowiec.20} for details. However, the functional $\cV_\epsilon$ is not a regular measure of regret; $\cV_\epsilon(\bfxi)$ might not be strictly larger than $\Ex[\bfxi]$.\\

\noindent {\bf Approximations using Sampling.} Almost always in practice, the ``true'' probability distribution of a random variable $\bfxi$ is unknown but one might have the ability to generate a sample $\xi_1, \dots, \xi_\nu$, which produces an empirical distribution. Let $\bfxi^\nu$ be a random variable with this empirical distribution, i.e., each of the values $\xi_i$, $i = 1, \dots, \nu$, occurs with probability $1/\nu$. This means that while the actual risk  $\cR(\bfxi)$ might not be computable because the distribution of $\bfxi$ is unknown, that of $\bfxi^\nu$ tends to be available but of course is also a bit different. From \cite{Shapiro.13} we know that if $\cR$ is a law-invariant, real-valued, monotone, and regular risk measure, then $\cR(\bfxi^\nu)\to \cR(\bfxi)$ as the sample size $\nu\to \infty$ almost surely\footnote{Since the sample that generates $\bfxi^\nu$ is random, the ``almost surely'' refers to the fact that the convergence holds for all samples $\{\xi_1, \xi_2, \dots\}$ in a set of probability one.}.  In fact, this also holds beyond random variables defined on $\cL^2$; see \cite{Shapiro.13} for details and \cite{Delbaen.21} for recent refinements.  (We note that the results from \cite{Shapiro.13} and those in the following theorem rely on the minor technical assumption that the probability space $(\Omega,\cA,\bbP)$ is atomless and complete.)

\begin{theorem}\label{tEpiconv}{\rm (epi-convergence under sampling)}.
For a law-invariant, monotone, and regular measure of risk $\cR:\cL^2\to \reals$, consider $f:\reals^m\times \reals^n\to \reals$, an $m$-dimensional random vector $\bfxi$ with support $\Xi$, and its sample-based approximation $\bfxi^\nu$ (see the previous paragraph). Suppose that the following hold:
\begin{enumerate}[(a)]

\item $f(\bfxi,x)\in \cL^2$ for all $x\in \reals^n$.

\item $f$ is a random lsc function\footnote{A sufficient (but not necessary) condition for $f$ to be random lsc is that it is measurable in its first argument and continuous in its second argument; see, e.g., \cite[Section 8.G]{primer}.}.

\item For every $x\in \reals^n$, there is neighborhood $N$ and $g:\reals^m\to [0,\infty)$ such that $g(\bfxi) \in \cL^2$ and $f(\xi,x) \geq g(\xi)$ for all $x\in N$ and $\xi \in \Xi$.

\end{enumerate}
Consider the functions $\phi,\phi^\nu:\reals^n\to \reals$ given by $\phi(x) = \cR(f(\bfxi,x))$  and $\phi^\nu(x) = \cR(f(\bfxi^\nu,x))$. Then, $\phi$ is lsc and $\phi^\nu$ epi-converges to $\phi$ as $\nu\to \infty$ almost surely.

Moreover, if $f(\xi, \cdot\,)$ is convex for all $\xi\in \Xi$, then assumptions (b,c) can be dropped and the conclusion strengthened to $\phi^\nu$ converging to $\phi$ uniformly on compact sets as $\nu\to \infty$ almost surely.
\end{theorem}
\state Proof. The result appears in \cite{Shapiro.13}, which also handles integrable random variables.\eop

The main consequence of $\phi^\nu$ epi-converging to $\phi$ is that all cluster points of a sequence of minimizers of $\phi^\nu$ are minimizers of $\phi$. For an introduction to epi-convergence, we refer to \cite[Section 4.C]{primer}. Thus, under the assumptions of the theorem, we are on solid footing if we approximate a distribution by an empirical distribution in a risk-minimization problem.
For the effect of more general changes to a probability distribution, we refer to \cite{Pichler.13b,KieselRuhlickeStahlZheng.16,ClausKratschmerSchultz.17,EmbrechtsSchiedWang.22}. In
\cite{PflugWozabal.10}, we find expressions for the asymptotic distributions of sample-based approximations of law-invariant risk measures; far reaching extensions emerge in \cite{DentchevaPenevRuszczynski.17} for nested expectations and, recently, in \cite{DentchevaLinPenev.22} for kernels, wavelets, and other estimators. Asymptotics for minimum values and minimizers appear in \cite{GuiguesKratschmerShapiro.18}. We also have uniform error bounds for broad classes of risk measures \cite{LeeParkShin.20,KhimLeqiPrasadRavikumar.20,LeqiHuangLiptonAzizzadenesheli.22}. In the specific setting of a derivative-free stochastic mirror descent algorithm and spectral risk measures, \cite{HollandHaress.21} develops generalization errors of the expectation and high-probability kinds. The paper \cite{KalogeriasPowell.18} establishes convergence rates of stochastic subgradient methods for mean-plus-semideviation risk measures, while \cite{GilesHajiali.19} obtains rates for multi-level nested simulations of risk measures.

\subsection{Mixed Superquantiles and Law Invariance}\label{subsec:mixed}

We can extend the mixing of risk measures in Example \ref{eMeasRegRisk}(d) to integrals over an uncountable collection of risk measures, especially superquantiles, to produce new regular measures of risk as pioneered by Kusuoka \cite{Kusuoka.01}. This subsection confirms that a large class of reasonable risk measures stem from superquantiles in this manner, which further highlights the centrality of superquantiles in decision making under uncertainty.

\begin{definition}{\rm (mixed superquantile measure of risk).} A weighting measure\footnote{A weighting measure is a probability measure on the Borel subsets of $[0,1)$.} $\lambda$ produces a {\em mixed superquantile measure of risk} $\cR:\cL^2\to (-\infty, \infty]$ given by
\begin{equation}\label{eqn:R}
  \cR(\bfxi) = \int_0^1 \bar Q(\beta) ~d\lambda(\beta),
\end{equation}
where $\bar Q(\beta)$ is the $\beta$-superquantile of $\bfxi$.
\end{definition}

If the weighting measure is supported on $\{\alpha_i \in [0, 1), ~i = 1, \dots, q\}$, then the resulting mixed superquantile measure of risk reduces to Example \ref{eMeasRegRisk}(d) with $\cR_i = \srsk_{\alpha_i}$, $i=1, \dots, q$. A common choice is $\alpha_1 = 0$ and $\alpha_2 = \alpha$, which produces $\cR(\bfxi) = \lambda_1 \Ex[\bfxi] + \lambda_2 \srsk_\alpha(\bfxi)$, where $\lambda_1,\lambda_2$ are nonnegative weights summing to 1; see for example \cite{FrohlichWilliamson.22} for its use in machine learning and Section \ref{sec:addexamples}.

The definition of mixed superquantile measures of risk allows for ``averaging'' more than a finite number
of superquantiles. Given $\alpha\in [0,1)$, one possibility places zero weight on superquantiles $\bar Q(\beta)$ for $\beta \in [0,\alpha)$ and weight $1/(1-\alpha)$ on $\bar Q(\beta)$ for $\beta \in [\alpha, 1)$. This produces the $\alpha$-{\it second-order superquantile} of a random variable $\bfxi$  \cite{RockafellarRoyset.18}:
\begin{equation}\label{eqn:barbar}
   \bar{\bar Q}(\alpha)= \frac{1}{1-\alpha}\int_\alpha^1{\bar Q}(\beta)d\beta.
\end{equation}
Second-order superquantiles define measures of risk, which are more conservative than those based on quantiles and superquantiles. We see that $Q(\alpha) \leq \bar Q(\alpha) \leq \bar{\bar Q}(\alpha)$ regardless of $\alpha\in (0,1)$. In addition to producing new measures of risk, second-order superquantiles also underpin the argmin-formula for superquantiles \eqref{eqn:argminsuper}; see \cite{RockafellarRoyset.18,Kouri.19b}.

Properties of mixed superquantile risk measures are traced back to \cite{Acerbi.02,RockafellarUryasevZabarankin.02,RockafellarUryasevZabarankin.06a}; here we summarize key facts as given in \cite{RockafellarRoyset.18}.

\begin{proposition}\label{prop:mixed} {\rm (mixed superquantile properties).} A mixed superquantile risk measure $\cR$ as defined in \eqref{eqn:R} is well-defined, monotone, and positively homogeneous. It is regular if $\lambda(\{0\}) < 1$, but lacking averseness if $\lambda(\{0\}) = 1$. It is real-valued on $\cL^2$ whenever $\lambda$ satisfies $\int_0^1 (1-\beta)^{-1/2}\,d\lambda(\beta) <\infty$ and, regardless of the weighting measure, has $\cR(\bfxi)<\infty$ whenever $\sup \bfxi<\infty$.

In terms of the quantile function $Q$ of $\bfxi$, one has the alternative expression
\begin{equation}\label{eqn:riskprofile}
  \cR(\bfxi) = \int_0^1 Q(\beta)\phi(\beta)d\beta, \mbox{ where } \phi(\beta)
      = \int_{0\leq \alpha < \beta} \frac{1}{1-\alpha} ~d\lambda(\alpha),
     ~\beta\in [0,1].
\end{equation}
The {\em risk profile function} $\phi$ is right-continuous and nondecreasing on $[0, 1]$ with $\phi(0) = 0$ and satisfies $\int_0^1 (1-\alpha)d\phi(\alpha) = 1$.  Conversely, any $\phi$ with these properties arises from a unique weighting measure $\lambda$ given by $d\lambda(\alpha) = (1-\alpha)d\phi(\alpha)$.
\end{proposition}

The alternative expression in \eqref{eqn:riskprofile} makes a deep connection with the integral formula \eqref{eqn:superquantAcerbi} for a (single) superquantile. The latter averages all quantiles above $\alpha$. In contrast, \eqref{eqn:riskprofile} potentially weighs the quantiles differently using a risk profile function $\phi$. The resulting risk measures are also known as spectral risk measures \cite{Acerbi.02}. They connect fundamentally with distortion functionals \cite{Pflug.06} common in insurance applications. There are further connections with dual utility theory \cite{Yaari.87}; see also \cite{DentchevaRuszczynski.13}.

We refer to \cite{Kouri.19b} for ways to compute mixed superquantile risk using numerical integration. The paper \cite{VijayanPrashanth.21} exemplifies recent efforts to use mixed superquantile risk measures in reinforcement learning. Since the risk profile function models a decision maker's preferences, it is often unknown and \cite{WangXu.20} examines worst-case models over a class of such functions; see also \cite{GuoXu.22}. The paper \cite{Li.18} examines worst-case values of law-invariant measures of risk under only partial information about random variables.

While the class of mixed superquantile risk measures is clearly large, the pioneering contribution in  \cite{Kusuoka.01} was to show that {\em every} real-valued, law-invariant, positively homogeneous, monotone, regular measure of risk $\cR$ can be written in the form
\[
\cR(\bfxi)  = \sup_{\lambda \in \Lambda} \int_0^1 \bar Q(\beta)\, d\lambda(\beta)
\]
for some set $\Lambda$ of weighting measures. This is the {\em Kusuoka representation} of $\cR$. (The paper \cite{Kusuoka.01} showed the existence of such $\Lambda$ for bounded random variables, with extensions to $p$-integrable random variables appearing in \cite{PflugRomisch.07}, refinements about uniqueness being established by \cite{Shapiro.13b}, and insight in the case of atomic probability spaces emerging from \cite{NoyanRudolf.15}.) Thus, superquantiles are the fundamental building blocks of a large class of meaningful measures of risk.

An immediate consequence of the Kusuoka representation is the following. For two random variables $\bfxi_0$ and $\bfxi_1$ and {\em any} real-valued, law-invariant, positively homogeneous, monotone, regular measure of risk $\cR$, one has that
\[
\srsk_\alpha(\bfxi_0) \leq \srsk_\alpha(\bfxi_1) ~~\forall \alpha \in [0,1) ~~~~\Longrightarrow~~~~ \cR(\bfxi_0) \leq \cR(\bfxi_1).
\]
Thus, if $\bfxi_0$ dominates $\bfxi_1$ in the sense of the left-hand side of the implication, then it becomes less critical to determine the ``right'' measure of risk; most meaningful risk measures will prefer $\bfxi_0$ over $\bfxi_1$. (This notion of dominance is equivalent to second-order stochastic dominance; see \cite{DentchevaMartinez.12,RockafellarRoyset.14}.) Even without such complete dominance, we can plot $\srsk_\alpha(\bfxi_0)$ and $\srsk_\alpha(\bfxi_1)$ as functions of $\alpha$ to highlight the pros and cons with each decision (random variable). In the context of machine learning, \cite{FrohlichWilliamson.22} generates such plots for the purpose of comparing different statistical models and their resulting errors.

\subsection{Fairness}

Fairness is an emerging area for application of risk measures. As pioneered in \cite{WilliamsonMenon.19} (see also \cite{LeqiPrasadRavikumar.19,FrohlichWilliamson.22}), one can even out the performance of a classifier across various subgroups by replacing the usual expectation with another risk measure during the training. Specifically, consider a quantity of interest $f((\xi,\eta),c)$ representing loss or estimation error, where $c\in\reals^n$ is a vector of coefficients specifying a statistical model (e.g., a neural network) and $(\xi,\eta)\in\reals^m$ is a vector consisting of input (features) $\xi$ and output (labels) $\eta$. An algorithm for supervised learning may seek to find $c$ that minimizes $\Ex[f((\bfxi,\bfeta),c)]$, where $(\bfxi,\bfeta)$ is a random vector with some assumed probability distribution typically taken as the empirical distribution of a training data set. The choice of expectation as risk measure may produce highly uneven performance of the resulting statistical model across subgroups. For certain types of outcomes of $(\bfxi,\bfeta)$, for instance those corresponding to male customers, the statistical model may over-estimate the output $\eta$, while for other outcomes, corresponding to female customers, the model may under-estimate the output. Biases of these kinds are problematic and sometimes illegal.

Suppose that we augment the input-output vector by also including an input $\zeta$, which represents sensitive characteristics such as gender. Thus, the random vector $(\bfxi,\bfeta,\bfzeta)$ has now three parts. As proposed in \cite{WilliamsonMenon.19}, we may change the training problem from that of minimizing $\Ex[f((\bfxi,\bfeta,\bfzeta),c)]$ to solving
\begin{equation}\label{eqn:fairnessWM}
\nnmin_{c\in \reals^n}\, \cR\big(g(\bfzeta,c)\big), ~~\mbox{ where } ~g(\zeta,c) = \Ex\Big[ f\big((\bfxi,\bfeta,\bfzeta),c\big)~\big|~\bfzeta = \zeta\big]
\end{equation}
and $\cR$ is a measure of risk that acts on a new quantity of interest $g(\zeta,c)$, which is uncertain because the value of $\zeta$ is governed by the distribution of $\bfzeta$. For a fixed $\zeta$ (say ``male'') and coefficient vector $c$ describing a statistical model, $g(\zeta,c)$ is the conditional expectation of the original quantity of interest, given $\bfzeta = \zeta$. Thus, it quantifies the average performance of the statistical model when applied to subjects with sensitive feature $\zeta$. If $g(\zeta,c)$ is much higher for some values of $\zeta$ than others, then the statistical model specified by $c$ might be perceived as unfair. The role of the risk measure $\cR$ is to assess the random variable $g(\bfzeta,c)$ and to safeguard against high values. Thus, a solution $c^\star$ of \eqref{eqn:fairnessWM} tends to specify a statistical model that performs reasonably well regardless of sensitive characteristics. We avoid an upper tail in the distribution of $g(\bfzeta,c^\star)$ that extends far to the right. Since quantities of interest in machine learning are often nonnegative, the left tail may also be ``light.'' This means that $g(\bfzeta,c^\star)$ has little variability, which might be our goal. (The paper \cite{WilliamsonMenon.19} defines {\em perfect fairness} of a statistical model given by $c$ as having $g(\zeta,c)= g(\zeta',c)$ for all possible values $\zeta,\zeta'$.)

In certain settings, one is not allowed (by law or company policy) to use sensitive data during the training and testing of an artificial intelligence system \cite{HashimotoSrivastavaNamkoongLiang.18}. This precludes a direct application of the above approach. Still, by switching from expected loss to a risk measure, we can control the upper tail of the resulting loss distribution across all individuals. This safeguards against statistical models with poor performance for some individuals and excellent performance for others. Consequently, even if the sensitive characteristic for each individual is unknown, we can achieve more consistent performance by using risk measures; see also \cite{HashimotoSrivastavaNamkoongLiang.18} for a related approach based on distributionally robust optimization.

\section{Error and Deviation}\label{sec:errordev}

The broad approaches to decision making based on measures of risk and regret as developed in the previous section extend to statistical concepts such as mean-squared error, standard deviation, and generalized regression. Already in the pioneering work \cite{Markowitz.59} on portfolio optimization, we find a discussion of semivariation, semideviation, and expected absolute deviation as means of quantifying ``variability'' in a random variable and the ``skewness'' of distributions; see, e.g., \cite{Speranza.93,KonnoShirakawa.94,OgryczakRuszczynski.99,OgryczakRuszczynski.02} for related quantities and computational methods. The first systematic studies of general classes of {\em measures of deviation} and {\em measures of error} appear to be \cite{RockafellarUryasevZabarankin.02,RockafellarUryasevZabarankin.06a,RockafellarUryasevZabarankin.08}, which approach the concepts axiomatically. Further connections with measures of risk and regret appear in \cite{RockafellarUryasev.13}; see \cite{RockafellarRoyset.15} for technical refinements. This produces quadrangles of risk, regret, error, and deviation measures with new insights and computational possibilities.

\subsection{Measures of Error}

We start by examining how to quantify the ``nonzeroness'' of a random variable and follow \cite{RockafellarUryasevZabarankin.08,RockafellarUryasev.13,RockafellarRoyset.15}.

\begin{definition}{\rm (regular measure of error and its statistic).}\label{dErrorStat} A {\em regular measure of error} $\cE$ is a functional from $\cL^2$ to $[0,\infty]$ that is lsc and convex, with
\[
\cE(\bfnull) = 0, ~\mbox{ but }~ \cE(\bfxi) > 0 ~~\forall \bfxi \neq \bfnull.
\]
The quantity $\cE(\bfxi)$ is the {\em error} of $\bfxi$. The corresponding {\em statistic} is $\cS(\bfxi)  = \nargmin_{\gamma\in\reals} \cE(\bfxi-\gamma)$.
\end{definition}

Measures of error are important in regression analysis where the goal is to minimize the ``nonzeroness'' of certain random variables representing the ``residual'' difference between observations and model predications. The statistic of $\bfxi$ is the set of scalars, possibly a single number, that best approximate a random variable $\bfxi$ in the sense of the corresponding regular measure of error. (Inherited from \cite{RockafellarUryasevZabarankin.08}, we use ``statistic'' in a somewhat different meaning than in the Statistics literature.)
We note that the definition in \cite{RockafellarUryasev.13} includes a limiting condition that is shown to be superfluous in \cite{RockafellarRoyset.15}.

\begin{example}{\rm (regular measures of error).}\label{eMeasRegError} We have the following examples of regular measures of error $\cE$ and corresponding statistics $\cS$; see, e.g., \cite[Example 8.15]{primer} for justifications of regularity.

\begin{enumerate}[{\rm (a)}]
\item {\em Koenker-Bassett error} and {\em quantile statistic}\footnote{For the definition of $Q^\lplus(\alpha)$, see \eqref{eqn:argminquant}.} with $\alpha \in (0,1)$:
\begin{flalign*}
& \cE(\bfxi) = \frac{1}{1-\alpha}\Ex\big[\max\{0, \bfxi\}\big]-\Ex[\bfxi] & \cS(\bfxi)= \big[Q(\alpha), Q^\lplus(\alpha)\big].
\end{flalign*}

\item {\em Worst-case error} and {\em statistic}:
\begin{flalign*}
& \cE(\bfxi) = \begin{cases}
-\Ex[\bfxi] & \mbox{ if } \nsup \bfxi \leq 0\\
\infty & \mbox{ otherwise}
\end{cases}
&  \cS(\bfxi) = \begin{cases}
\{\sup \bfxi\} & \mbox{ if } \sup \bfxi<\infty\\
\emptyset & \mbox{ otherwise}.
\end{cases}
\end{flalign*}

\item $\cL^2$-{\em error} and {\em statistic} with $\lambda>0$:
\begin{flalign*}
& \cE(\bfxi) = \lambda \sqrt{\Ex[\bfxi^2]} = \lambda\|\bfxi\|_{\cL^2} &  \cS(\bfxi)= \big\{\Ex[\bfxi]\big\}.
\end{flalign*}

\item {\em Mixed error} and {\em statistic} with $\lambda_1, \dots, \lambda_q>0$ and $\nsum_{i=1}^q \lambda_i = 1$:
\begin{flalign*}
& \cE(\bfxi) = \inf_{\gamma_1, \dots, \gamma_q} \Big\{\nsum_{i=1}^q \lambda_i \cE_i(\bfxi-\gamma_i)~\Big|~\nsum_{i=1}^q \lambda_i \gamma_i = 0 \Big\} &
\cS(\bfxi)  = \nsum_{i=1}^q \lambda_i \cS_i(\bfxi),
\end{flalign*}
where $\cS_i(\bfxi)  = \nargmin_{\gamma\in\reals} \cE_i(\bfxi-\gamma)$ and $\cE_i$ is a regular measure of error.
\end{enumerate}
\end{example}

In regression analysis, we seek to predict (or forecast) the value of a random variable $\bfeta$ from the values of some other random variables $\bfxi_1, \dots, \bfxi_n$. For instance, these other random variables could be input to a system, which we observe or perhaps even control, and $\bfeta$ is the output of the system, which we hope to predict. Let $\bfxi = (\bfxi_1, \dots, \bfxi_n)$ and suppose that $\bfeta,\bfxi_1, \dots, \bfxi_n\in \cL^2$. We may seek a ``best'' statistical model of the form $\gamma + \langle c, \bfxi\rangle$ using a regular measure of error $\cE$. This leads to the {\em generalized regression problem}
\begin{equation}\label{eqn:regProblem}
\nnmin_{(\gamma,c)\in \reals^{1+n}} ~\cE\big(\bfeta - \gamma - \langle c, \bfxi\rangle\big).
\end{equation}
A minimizer $(\gamma^\star, c^\star)$ of the problem defines the random variable
$\gamma^\star + \langle c^\star, \bfxi\rangle$, which then makes the ``nonzeroness'' of $\bfeta - (\gamma + \langle c, \bfxi\rangle)$ as low as possible in the sense of the selected measure of error. Thus, $\gamma^\star + \langle c^\star, \bfxi\rangle$ is a best possible approximation of $\bfeta$ in this sense.
Since $\cE$ is convex, it follows that the function $(\gamma,c)\mapsto \cE(\bfeta - \gamma - \langle c, \bfxi\rangle)$ is convex. Thus, \eqref{eqn:regProblem} is computationally appealing.

The choice in Example \ref{eMeasRegError}(c) leads to {\em least-squares regression} regardless of $\lambda>0$. Since the measure of error can be squared without affecting the set of minimizers, in this case one has
\begin{equation*}
\nargmin_{(\gamma, c)\in\reals^{1+n}} \cE\big(\bfeta - \gamma - \langle c, \bfxi\rangle\big) = \nargmin_{(\gamma, c)\in\reals^{1+n}} \Ex\Big [\big(\bfeta - \gamma - \langle c, \bfxi\rangle\big)^2\Big].
\end{equation*}
In fact, $\bfeta\mapsto \Ex[\bfeta^2]$ is also a regular measure of error and we could just as well have adopted it from the start. 

There is a long tradition for considering measures of error beyond least-squares, especially in robust statistics; see \cite{Huber.81}. The Koenker-Bassett error in Example \ref{eMeasRegError}(a) leads to {\em quantile regression}. Figure \ref{fig:lift} illustrates statistical models for predicting lift force using least-squares regression (solid line) and quantile regression with $\alpha = 0.5, 0.75$ (dashed lines) and $\alpha = 0.95, 0.995$ (dotted lines). While the lines are quite similar, quantile regression with higher values of $\alpha$ results in more conservative estimates in the sense that the values of $\gamma + c \, \bfxi$ tend to exceed those of $\bfeta$.

As eluded to in Subsection \ref{subsec:stat}, one can go much beyond linear statistical models in \eqref{eqn:regProblem}. We refer to \cite{RockafellarUryasevZabarankin.08} for several refinements and to  \cite{GrechukZabarankin.19} for a discussion of the connection between regular measures of error and inter-regenerative relationships via log-likelihood and entropy maximization. Innovative thinking about alternative regression approaches based on constrained residuals, leveraging superquantiles, appears in \cite{TrindadeUryasevShapiroZrazhevsky.07}.

\subsection{Measures of Deviation}

We next quantify ``nonconstancy'' of a random variable. This can be accomplished using the standard deviation, but there are many other possibilities as well; see the axiomatic development in \cite{RockafellarUryasevZabarankin.02,RockafellarUryasevZabarankin.06a}. Here, we adopt the definition in \cite{RockafellarUryasev.13}.

\begin{definition}{\rm (regular measure of deviation).} A {\em regular measure of deviation} $\cD$ is a functional from $\cL^2$ to $[0,\infty]$ that is lsc and convex, with
\[
\cD(\bfxi) = 0 ~\mbox{ if }~ \bfxi ~\mbox{ is constant; }~~~~ \cD(\bfxi) > 0 \,\mbox{ otherwise}.
\]
The quantity $\cD(\bfxi)$ is the {\em deviation} of $\bfxi$.
\end{definition}

\begin{example}{\rm (regular measures of deviation).}\label{eMeasRegDev} We have the following examples of regular measures of deviation $\cD$; see, e.g.,  \cite[Example 8.18]{primer} for supporting arguments.
\begin{enumerate}[{\rm (a)}]

\item {\em Superquantile deviation} with $\alpha \in (0,1)$:
\begin{flalign*}
& \cD(\bfxi) = \srsk_\alpha(\bfxi) - \Ex[\bfxi].&
\end{flalign*}

\item {\em Worst-case deviation}:
\begin{flalign*}
& \cD(\bfxi) = \nsup \bfxi - \Ex[\bfxi].&
\end{flalign*}

\item {\em Scaled standard deviation} with $\lambda>0$:
\begin{flalign*}
& \cD(\bfxi) = \lambda \std(\bfxi).&
\end{flalign*}

\item {\em Mixed deviation} with $\lambda_1, \dots, \lambda_q>0$ and $\nsum_{i=1}^q \lambda_i = 1$:
\begin{flalign*}
& \cD(\bfxi) = \nsum_{i=1}^q \lambda_i \cD_i(\bfxi),~ \mbox{ for regular measures of deviation } \cD_1, \dots, \cD_q.&
\end{flalign*}
\end{enumerate}

\end{example}

The regular measures of deviation in Example \ref{eMeasRegDev}(a,b,c) share the common property that they are obtained from the regular measures of risk in Example \ref{eMeasRegRisk}(a,b,c), respectively, by subtracting $\Ex[\bfxi]$. Similarly, the regular measures of error in Example \ref{eMeasRegError}(a,b,c) can be constructed from the regular measures of regret in Example \ref{eMeasRegRisk}(a,b,c), respectively, by subtracting $\Ex[\bfxi]$. These connections hold in general, as the next theorem (adapted from \cite{RockafellarRoyset.15}) asserts.

\begin{theorem}\label{thm:quadrangle} {\rm (expectation translations).} Every regular measure of deviation $\cD$ defines a regular measure of risk $\cR$ and vice versa through the
relations:
\[
  \cR(\bfxi) = \cD(\bfxi) + \Ex[\bfxi]~~~ \mbox{ and } ~~~  \cD(\bfxi) = \cR(\bfxi) - \Ex[\bfxi].
\]

Similarly, every regular measure of error $\cE$ defines a regular measure of regret $\cV$ and vice versa
through the relations:
\[
  \cV(\bfxi) = \cE(\bfxi) + \Ex[\bfxi] ~~~\mbox{ and }~~~ \cE(\bfxi) = \cV(\bfxi) - \Ex[\bfxi].
\]
\end{theorem}

We complete the picture by connecting regular measures of error and deviation in a manner that resembles the relation between regular measures of regret and risk; see Theorem \ref{tRiskregret}.

\begin{theorem}\label{thm:errdev} {\rm (error-deviation).} For a regular measure of error $\cE$, a regular measure of deviation $\cD$ is obtained by
\[
\cD(\bfxi) = \min_{\gamma\in \reals} \,\cE(\bfxi-\gamma).
\]
For every regular measure of deviation $\cD$ there is a regular measure of error $\cE$, not necessarily unique, that constructs $\cD$ through this minimization formula.

Moreover, if $\cE$ is paired with a regular measure of regret $\cV$ via Theorem \ref{thm:quadrangle}, then
\begin{equation}\label{eqn:argminboth}
  \cS(\bfxi)  = \nargmin_{\gamma\in\reals} \cE(\bfxi-\gamma) = \nargmin_{\gamma\in\reals} \big\{\gamma + \cV(\bfxi-\gamma)\big\}
\end{equation}
and this set is a nonempty compact interval as long as $\cV(\bfxi-\gamma)$, or
equivalently $\cE(\bfxi-\gamma)$, is finite for some $\gamma\in\reals$.
\end{theorem}
\state Proof. This fact stems from \cite{RockafellarUryasevZabarankin.08}, but here stated as it appears in \cite[Theorem 8.21]{primer} which incorporates the refinements of \cite{RockafellarRoyset.15}.\eop

\drawing{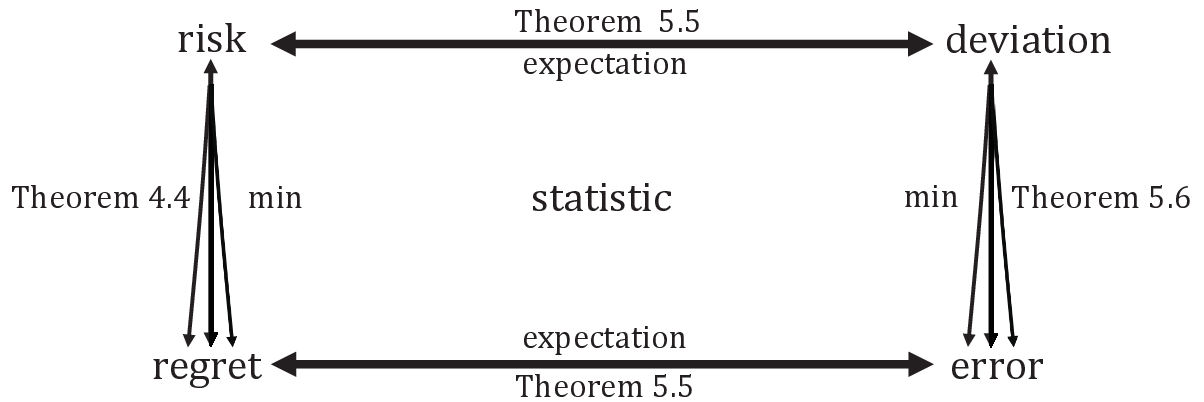}{4.6in} {Relations between regular measures of risk, regret, deviation, and error that form a risk quadrangle.} {fig:quad}

The key insight of \cite{RockafellarUryasev.13} (with technical refinements in \cite{RockafellarRoyset.15}) as summarized in Theorems \ref{tRiskregret}, \ref{thm:quadrangle}, and \ref{thm:errdev} is that regular measures of risk, regret, deviation, and error can be connected in a {\em risk quadrangle}. Subsection \ref{subsec:riskmin} shows how passing to a regular measure of regret may simplify risk minimization. Similarly, a regular measure of deviation can be evaluated by first minimizing regret as in \eqref{eqn:argminboth}. One can construct new regular measures of risk and deviation by starting from either a regular measure of regret or a regular measure of error using the relations summarized in Figure \ref{fig:quad}. The ``horizontal'' connections in the figure are one-to-one as seen in Theorem \ref{thm:quadrangle}, while the ``vertical'' connections are not unique in the sense that multiple regular measures of regret produce the same regular measure of risk, with a similar situation taking place when passing from error to deviation.

The regular measures of risk, regret, deviation, and error labeled (a) in Examples \ref{eMeasRegRisk}, \ref{eMeasRegError}, and \ref{eMeasRegDev} are connected in a risk quadrangle in the sense of Figure \ref{fig:quad}. Likewise, the regular measures of risk, regret, deviation, and error labeled (b,c,d) in Examples \ref{eMeasRegRisk}, \ref{eMeasRegError}, and \ref{eMeasRegDev} are connected, respectively, and then also form risk quadrangles. For many risk quadrangles, we refer to Section \ref{sec:addexamples} and \cite{RockafellarUryasev.13,Kouri.19b}. Further connections with maximum entropy is developed in
\cite{GrechukMolybohaZabarankin.09}.

The relations between regular measures of error and deviation lead to a decomposition of the generalized  regression problem \eqref{eqn:regProblem} first identified in \cite{RockafellarUryasevZabarankin.08} for positively homogeneous functionals. Here, we state the fact as it appears in \cite{RockafellarRoyset.15}.

\begin{proposition}\label{prop:errdev} {\rm (decomposition).} For a regular measure of error $\cE$ and a regular measure of deviation $\cD$ paired by Theorem \ref{thm:errdev}, one has for any $\bfeta, \bfxi_1, \dots, \bfxi_n\in \cL^2$,
\[
(\gamma^\star, c^\star) \in \nargmin_{\gamma,c} \cE\big(\bfeta - \gamma - \langle c, \bfxi\rangle\big) ~~\Longleftrightarrow~~ c^\star \in \nargmin_{c} \cD\big(\bfeta - \langle c, \bfxi\rangle\big) \,\mbox{ and }\, \gamma^\star \in \cS\big(\bfeta - \langle c^\star, \bfxi\rangle\big),
\]
where $\bfxi = (\bfxi_1, \dots, \bfxi_n)$ and $\cS$ is the statistic corresponding to $\cE$.
\end{proposition}

The proposition shows that we can determine the coefficients $(\gamma^\star, c^\star)$ in the generalized  regression problem \eqref{eqn:regProblem} using two steps: first, determine $c^\star$ by minimizing deviation and, second, fix $\gamma^\star$ by computing a statistic. From Theorem \ref{thm:quadrangle} we see that minimizing deviation is equivalent to minimizing risk modified by an expectation. Thus, $c^\star$ can often be computed by algorithms for minimizing risk, with minor adjustments.

Separate roles for the ``slope'' $c$ and the ``intercept'' $\gamma$ emerge from the proposition. The former is selected to minimize the ``nonconstancy'' of $\bfeta-\langle c, \bfxi\rangle$, while the latter translates the functional $\bfxi \mapsto \langle c^\star, \bfxi\rangle$ up or down to match the correct statistic.

The relative quality of a statistical model obtained by solving the generalized  regression problem \eqref{eqn:regProblem} can be assessed using the following concept proposed in \cite{RockafellarRoyset.15}; see also precursors in \cite{RockafellarRoysetMiranda.14}.
For a regular measure of error $\cE$ paired with a regular measure of deviation $\cD$ via Theorem \ref{thm:errdev}, the {\em coefficient of determination}, or {\em R-squared}, of the random vector $\bfxi = (\bfxi_1, \dots, \bfxi_n)$ relative to $\bfeta$ is given as\footnote{Here, $\infty/\infty$ is interpreted as 1 and $0/0$ as 0.}
\begin{equation*}\label{eqn:hatR}
R^2 = 1 - \frac{\inf_{(\gamma, c) \in \reals^{1+n}}\cE\big(\bfeta - \gamma - \langle c, \bfxi\rangle\big)}{\cD(\bfeta)}.
\end{equation*}
The definition extends a classical concept from least-squares regression, where $\cD$ is the variance and $\cE$ is the mean-squared error. It is immediate from the definition that if $(\gamma^\star,c^\star)$ is a minimizer of the generalized  regression problem and $\cE(\bfeta - \gamma^\star - \langle c^\star, \bfxi\rangle) = 0$, then $R^2=1$. However, such vanishing error only takes place when the residual $\bfeta - \gamma^\star - \langle c^\star, \bfxi\rangle= \bfnull$; the model $\gamma^\star + \langle c^\star, \bfxi\rangle$ predicts $\bfeta$ perfectly. Under less ideal circumstances, one would assess the quality of $(\gamma^\star,c^\star)$ by seeing how close $R^2$ is to 1. In general, $R^2\geq 0$ and $R^2=0$ when $(\gamma,0)$ is a minimizer in \eqref{eqn:regProblem}. In that case, $\bfxi$ provides ``no information'' about $\bfeta$.

In Example \ref{eSurrogate}, we obtain five statistical models of the form $\gamma + c\,\bfxi$; see the lines in Figure \ref{fig:lift}. The quality of the models can be assessed using $R^2$. Quantile regression with $\alpha = 0.75$ (highest dashed line) produces a deviation of 7.97 and an error of 1.83 and thus $R^2 =  0.77$. For quantile regression with $\alpha = 0.95$ (lowest dotted line), we obtain a deviation of 12.32 and an  error of 3.29 and thus $R^2 =  0.73$.

The paper \cite{RockafellarUryasevZabarankin.08} brought to the forefront the possibilities of constructing conservative statistical models through the use of ``nonstandard'' measures of error and other adjustments.  Later developments include theoretical refinements in \cite{RockafellarRoyset.15}, applications to reliability engineering in \cite{JakemanKouriHuerta.22}, and applications to naval architecture in \cite{RoysetBonfiglioVernengoBrizzolara.17,BonfiglioRoyset.19}. Support vector regression is the topic of \cite{MalandiiUryasev.22}.

\subsection{Superquantile Regression}\label{subsec:superregression}

The formula \eqref{eqn:argminsuper} for superquantiles as minimizers of some functionals gives rise to {\em superquantile regression} \cite{RockafellarRoysetMiranda.14} that supplements least-squares and quantile regression. It stems from the insight that the functionals $\cR_\alpha, \cV_\alpha, \cD_\alpha, \cE_\alpha$ on $\cL^2$ given by
\begin{align*}
\cR_\alpha(\bfxi) & = \frac{1}{1-\alpha} \int_{\alpha}^1 \srsk_\beta(\bfxi) \, d\beta & \cD_\alpha(\bfxi) & = \cR_\alpha(\bfxi) - \Ex[\bfxi]\\
\cV_\alpha(\bfxi) & = \frac{1}{1-\alpha} \int_{0}^1 \max\big\{0, \srsk_\beta(\bfxi)\big\} \, d\beta & \cE_\alpha(\bfxi) & = \cV_\alpha(\bfxi) - \Ex[\bfxi]
\end{align*}
for $\alpha \in [0,1)$ are regular measures of risk, regret, deviation, and error, respectively, and form a risk quadrangle \cite{RockafellarRoyset.14}, i.e., they are connected in the sense of Theorems \ref{tRiskregret}, \ref{thm:quadrangle}, and \ref{thm:errdev}. The corresponding statistic is $\cS_\alpha(\bfxi) = \{\srsk_\alpha(\bfxi)\}$ so $\gamma\mapsto \cE_\alpha(\bfxi-\gamma)$ has the $\alpha$-superquantile of $\bfxi$ as its unique minimizer. We observed that {\em no} error measure of the expectation kind achieves an $\alpha$-superquantile as its minimizer when $\alpha \in (0,1)$; see \cite[Theorem 11]{Gneiting.11}. Thus, one needs to consider more complicated measures of error with $\cE_\alpha$ being one possibility, others are derived in \cite{Kouri.19b}. For an introduction to the area of {\em elicitation}, we refer to \cite{LambertPennockShoham.08}, with a machine learning perspective, and \cite{Steinwart2014ElicitationAI,BelliniBignozzi.15,Ziegel.16}; see also the recent paper \cite{FrongilloKash.21}.

Superquantile regression predicts a random variable $\bfeta$ from the random vector $\bfxi = (\bfxi_1, \dots, \bfxi_n)$ using the statistical model $\gamma^\star + \langle c^\star, \bfxi\rangle$, where $(\gamma^\star,c^\star)$ is a minimizer of \eqref{eqn:regProblem} with $\cE = \cE_\alpha$ for some $\alpha \in [0,1)$. In view of Proposition \ref{prop:errdev}, it might be computationally more efficient to obtain $(\gamma^\star, c^\star)$ in two steps: First, solve
\begin{equation}\label{eqn:superregproblem}
c^\star \in \nargmin_{c\in\reals^n} \bigg\{\frac{1}{1-\alpha} \int_{\alpha}^1 \srsk_\beta \hspace{-0.07cm}  \big(\bfeta - \langle c, \bfxi\rangle\big) \, d\beta - \Ex\big[ \bfeta - \langle c, \bfxi\rangle \big]\bigg\}.
\end{equation}
Second, set $\gamma^\star = \srsk_\alpha(\bfeta - \langle c^\star, \bfxi\rangle)$. 

If the probability distribution of $(\bfeta,\bfxi)$ is finite with equal probability at each atom, then \eqref{eqn:superregproblem} is equivalent to a linear optimization problem. Let $\{\eta_j\in\reals, \xi^j\in\reals^n, j = 1, \dots, s\}$ be the possible outcomes of $(\bfeta,\bfxi)$. As shown in \cite{RockafellarRoysetMiranda.14}, a minimizer of \eqref{eqn:superregproblem} is equivalently obtained by solving 
\begin{align*}
\nnmin_{c,u,v,\tau} \frac{1}{1 - \alpha}\sum_{i=s_\alpha}^{s-1} (\beta_i-\beta_{i-1})u_i +&
\frac{1}{s(1 - \alpha)} \sum_{i=s_\alpha}^{s-1} \sum_{j=1}^s a_i v_{ij} + \frac{\tau}{s(1 - \alpha)} - \frac{1}{s} \sum_{j=1}^s \big(\eta_j - \langle c,\xi^j\rangle\big)\\
\mbox{subject to } ~~\eta_j - \langle c,\xi^j\rangle - u_i & \leq v_{ij}, ~~i=s_\alpha, \ldots, s-1, ~j=1, \ldots, s\\
                                            0 & \leq v_{ij}, ~~i=s_\alpha, \ldots, s-1, ~j=1, \ldots, s\\
                   \eta_j - \langle c,\xi^j\rangle & \leq \tau, ~~\,~j=1, \ldots, s\\
 c \in \reals^n, ~~~ u = (u_{s_\alpha}, \ldots, u_{s-1}) & \in \reals^{s - s_\alpha}, ~~~ v = (v_{s_\alpha,1}, \ldots, v_{s-1,s}) \in \reals^{(s - s_\alpha)s}, ~~~ \tau \in \reals,
\end{align*}
where $s_\alpha = \lceil s\alpha \rceil$ is the smallest integer no smaller than $s\alpha$, $\beta_{s_\alpha -1} = \alpha$, $\beta_i = i/s$, for $i=s_\alpha, s_\alpha+1, \dots, s$, and $a_i = \ln(1-\beta_{i-1}) - \ln(1-\beta_{i})$ for $i=s_\alpha, s_\alpha+1, \dots, s-1$.
For a related formulation that makes explicit connection with mixed risk, we refer to \cite{GolodnikovKuzmenkoUryasev.19}. Alternatively, we can solve \eqref{eqn:superregproblem} by approximating the integral using standard numerical integration techniques \cite{RockafellarRoysetMiranda.14,Kouri.19b}. As laid out in \cite[Example 8.41]{primer}, one can also solve \eqref{eqn:superregproblem} using the subgradient method.

Superquantile regression as define above is {\em different} from minimizing $\srsk_\alpha(\bfeta - \gamma - \langle c, \bfxi\rangle)$. However, this distinction is not always made in the literature; see for example \cite[Example 2.2]{LaguelPillutlaMalickHarchaoui.21}. It is also different than attempts to estimate conditional superquantiles, i.e., superquantiles of random variables that depend on a vector of predictors; see \cite{Scaillet.05,CaiWang.08,Kato.12} for kernel-based approaches and \cite{PeracchiTanase.08,LeoratoPeracchiTanase.12,ChunShapiroUryasev.12} for methods passing through \eqref{eqn:superquantAcerbi}. Still, we can view \cite{ChunShapiroUryasev.12} as approaching superquantile regression via an approximation based on the mixed error measure in Example \ref{eMeasRegError}(d) with $\cE_i$, $i = 1, \dots, q$, being Koenker-Bassett error measures at different probability levels $\alpha_i$. That reference also estimates conditional superquantiles using an appropriately shifted least-squares regression curve based on superquantiles of the resulting residuals.

\section{Duality Theory}\label{sec:duality}

A fundamental fact from convex analysis is that a convex function, with minor exceptions, can be expressed as the supremum over a collection of affine functions. Since regular measures of risk, regret, error, and deviation are convex by definition and also satisfy the required technical conditions, every such functional has an alternative representation. This offers computational possibilities and help with interpretation. Most significantly, it makes connections with distributionally robust optimization as previewed in the discussion of Mr.~Averse and Ms.~Ambiguous (cf. Subsection \ref{subsec:equivalence} and \eqref{eqn:maxsuper}). We now examine such connections broadly, and go much beyond superquantiles.

\subsection{Conjugacy}

On the space $\cL^2$ of random variables, we adopt the inner product $(\bfxi,\bfpi)\mapsto \Ex[\bfxi\bfpi]$ and this leads to a definition of conjugates.

\begin{definition}{\rm (conjugate).}\label{dConjugateFcnsRV}
For $\cF:\cL^2\to [-\infty,\infty]$, the functional $\cF^*:\cL^2\to [-\infty,\infty]$ defined by
\[
\cF^*(\bfpi) = \nsup\big\{\Ex[\bfxi\bfpi] - \cF(\bfxi)  ~\big|~ \bfxi\in \cL^2\big\}
\]
is the {\em conjugate} of $\cF$.
\end{definition}

A main motivation for restricting the attention to square-integrable random variables is that a conjugate is then defined on $\cL^2$ as well. In contrast, if $\bfxi$ were only integrable, then $\bfpi$ would have to be bounded for $\Ex[\bfxi\bfpi]$ to be finite; see \cite{RuszczynskiShapiro.06} for such extensions. A central result from convex analysis is the Fenchel-Moreau theorem; see, e.g., \cite[Theorem 5]{Rockafellar.74}.

\begin{theorem}{\rm (Fenchel-Moreau).}\label{pConjugateFcnal}
For a lsc convex functional $\cF:\cL^2\to (-\infty,\infty]$ with\footnote{The notation $\cF\not\equiv \infty$ simply indicates that the functional with $\cF(\bfxi) = \infty$ for all $\bfxi\in \cL^2$ is ruled out.} $\cF\not\equiv \infty$, the conjugate $\cF^*$ is also lsc and convex with $\cF^*\not\equiv \infty$, $\cF^*(\bfxi)>-\infty$ for all $\bfxi\in \cL^2$, and $(\cF^*)^* = \cF$. Thus,
\[
\cF(\bfxi) = \nsup\big\{\Ex[\bfxi\bfpi] - \cF^*(\bfpi)  ~\big|~ \bfpi\in \cL^2\big\}.
\]
\end{theorem}

Since regular measures of risk, regret, deviation, and error are lsc and convex by definition and also finite at some point, the Fenchel-Moreau theorem \ref{pConjugateFcnal} applies and furnishes alternative expressions. This observation is most useful if we can characterize the corresponding conjugates. While this is possible in many situations, we concentrate on the attractive case of positively homogeneous functionals; see \cite{Rockafellar.74,BauschkeCombettes.11} for comprehensive treatments. A summary of properties appears in \cite[Section 2]{RockafellarRoyset.15}.

We recall that the {\em domain} of $\cF:\cL^2\to [-\infty,\infty]$ is the set $\dom \cF = \{\bfxi \in \cL^2~|~\cF(\bfxi) < \infty\}$. With this notation, we have the following direct consequence of the Fenchel-Moreau theorem \ref{pConjugateFcnal}.

\begin{proposition}{\rm (conjugacy under positive homogeneity).}\label{pPosHom}
For a positively homogeneous, lsc, and convex functional $\cF:\cL^2\to (-\infty,\infty]$ with $\cF\not\equiv \infty$, one has
\begin{equation*}\label{eqn:conjugacyPosHom}
\cF(\bfxi) = \nsup\big\{\Ex[\bfxi\bfpi] ~\big|~ \bfpi\in \dom \cF^*\big\}.
\end{equation*}
A nonempty, closed, and convex set $\cC\subset \cL^2$ is the domain of $\cF^*$ for some positively homogeneous, lsc, and convex $\cF:\cL^2\to (-\infty,\infty]$ with $\cF\not\equiv \infty$ and then
\begin{equation*}\label{eqn:conjugacyPosHom1}
\cC = \big\{\bfpi \in \cL^2~\big|~ \Ex[\bfxi\bfpi] \leq \cF(\bfxi) ~~\forall \bfxi\in \cL^2\big\}.
\end{equation*}
\end{proposition}

The proposition asserts that positively homogeneous, lsc, and convex functionals are fully characterized by a subset of $\cL^2$. 

The restriction to positively homogeneous functionals still captures many important situations. For example, if $\cV $ is a positively homogeneous regular measure of regret, then the regular measure of risk $\cR$ constructed by Theorem \ref{tRiskregret} is also positively homogeneous. If $\cE $ is a positively homogeneous regular measure of error, then the regular measure of deviation $\cD$ constructed by Theorem \ref{thm:errdev} is also positively homogeneous. These facts follow immediately from the theorems. Since $\bfxi\mapsto \Ex[\max\{0, \bfxi\}]$ is positively homogeneous, $\cV$ in Example \ref{eMeasRegRisk}(a) is positively homogeneous and then also $\srsk_\alpha$. The measures of regret and risk in Example \ref{eMeasRegRisk}(b,c) are also positively homogeneous. For Example \ref{eMeasRegRisk}(d), $\cV$ is positively homogeneous when $\cV_1, \dots, \cV_q$ are positively homogeneous and then $\cR$ also has this property.

\subsection{Risk Envelopes and Dual Algorithms}

With the development of risk measures and related concepts, it was quickly realized that dual expressions were available via the Fenchel-Moreau theorem \ref{pConjugateFcnal}; see \cite{Delbaen.00,Delbaen.02,FollmerSchied.02c,RockafellarUryasevZabarankin.02}. In particular, this insight led to the dual formula \eqref{eqn:maxsuper} for superquantiles in \cite{Delbaen.00}. While these papers concentrate on bounded random variables or random variables in $\cL^2$, subsequent efforts \cite{Pflug.06b,RuszczynskiShapiro.06,RuszczynskiShapiro.06c} address more general spaces of random variables and also furnish many examples. We follow the development in \cite{RockafellarUryasev.13} and concentrate on random variables in $\cL^2$ for simplicity.

As a direct application of Proposition \ref{pPosHom}, we obtain that for a positively homogeneous regular measure of risk $\cR$, one has
\begin{equation}\label{eqn:poshomRisk}
\cR(\bfxi) = \nsup\big\{\Ex[\bfxi\bfpi] ~\big|~ \bfpi\in \Pi\big\},
\end{equation}
where, following the terminology of \cite{RockafellarUryasevZabarankin.02}, $\Pi = \dom \cR^*$ is the {\em risk envelope}  of $\cR$. Thus, a measure of risk of this kind is fully characterized by the domain of its conjugate. Convex analysis furnishes details about such domains. The following fact is taken from \cite[Proposition 8.30]{primer}.

\begin{proposition}{\rm (properties of risk envelopes).}\label{pPropRiskEnvelope} The risk envelope $\Pi$ of a positively homogeneous regular measure of risk $\cR$ is closed and convex and, for every $\bfpi \in \Pi$, one has $\Ex[\bfpi] = 1$. Moreover, it has the expression
\[
\Pi = \big\{\bfpi\in \cL^2~\big|~\Ex[\bfxi\bfpi] \leq \cR(\bfxi) ~~\forall \bfxi\in \cL^2\big\}.
\]

A  nonempty, closed, and convex set $\cC\subset \cL^2$ is the risk envelope of some positively homogeneous regular measure of risk on $\cL^2$ provided that $\cC$ also satisfies:
\begin{align*}
&\Ex[\bfpi] = 1 ~~~~\forall \bfpi \in \cC;\\
&\mbox{for each nonconstant } \bfxi\in \cL^2, \mbox{ there exists } \bfpi \in \cC \mbox{ such that } \,\Ex[\bfxi\bfpi] > \Ex[\bfxi].
\end{align*}
\end{proposition}

\begin{example}{\rm (risk envelopes).}\label{eRiskEnv} The risk envelopes corresponding to the regular measures of risk in Example \ref{eMeasRegRisk} are as follows; see, e.g., \cite[Example 8.31]{primer} for a justification:
\begin{enumerate}[{\rm (a)}]
\item {\em Superquantile risk envelope}:
\begin{flalign}\label{eqn:superriskenvelope}
&\Pi = \big\{\bfpi ~\big|~ 0 \leq \bfpi(\omega) \leq 1/(1-\alpha)~~\forall \omega, ~~~\Ex[\bfpi] = 1\big\}.&
\end{flalign}

\item {\em Worst-case risk envelope}:
\begin{flalign*}
&\Pi = \big\{\bfpi ~\big|~ \bfpi(\omega)\geq 0 ~~\forall \omega, ~~~ \Ex[\bfpi] = 1\big\}.&
\end{flalign*}

\item {\em Mean-plus-standard-deviation risk envelope}\footnote{The constant random variable with value 1 is denoted by $\bfone$.}:
\begin{flalign*}
&\Pi = \big\{\bfone + \lambda \bfpi ~\big|~ \Ex[\bfpi^2] \leq 1, ~~\Ex[\bfpi] = 0\big\}.&
\end{flalign*}

\item {\em Mixed risk envelope}:
\begin{flalign*}
&\Pi = \Big\{\nsum_{i=1}^q \lambda_i \bfpi_i~\Big|~\bfpi_i\in \Pi_i \Big\},&
\end{flalign*}
where $\Pi_i$ is the risk envelope of the positively homogeneous regular measure of risk $\cR_i$.
\end{enumerate}
\end{example}

The dual representation of risk measures in  \eqref{eqn:poshomRisk} offers several computational possibilities. In particular, it brings forward an expectation, which can subsequently be approximated using Monte Carlo sampling or other techniques.\\

\noindent {\bf Dual Algorithms for Risk Minimization.} The discussion in Subsection \ref{eqn:algosuper} about dual algorithms for superquantile minimization hints to possibilities that derive from \eqref{eqn:poshomRisk}. Concretely, let us consider a quantity of interest represented by $f:\reals^m\times\reals^n\to \reals$, a constraint set $X\subset \reals^n$, a positively homogeneous regular measure of risk $\cR$, and a finitely distributed random vector $\bfxi$ with support $\Xi\subset\reals^m$ of cardinality $s$ and probabilities $\{p_{\xi}>0, \xi\in \Xi\}$. This leads to the problem
\begin{equation}\label{eqn:riskmindual}
\nnmin_{x\in X} ~\cR\big(f(\bfxi,x)\big) = \nsup\big\{\Ex\big[ f(\bfxi,x) \bfpi \big] ~\big|~ \bfpi\in \Pi\big\} = \nsup \Big\{  \nsum_{\xi\in \Xi} p_{\xi} f(\xi, x) q_ {\xi} ~\Big|~q  \in \cQ\Big\},
\end{equation}
where $\Pi$ is the risk envelope of $\cR$. The second equality recognizes that the underlying probability space $(\Omega, \cA, \bbP)$ defining $\cL^2$ can in this case have $\Omega = \Xi$ and $\bbP(\{\xi\}) = p_{\xi}$. Thus, all random variables $\bfpi\in \cL^2$ are represented by $s$-dimensional vectors of the form $q = (q_{\xi}, \xi\in \Xi)\in \reals^s$ and the maximization over $\Pi$ in \eqref{eqn:riskmindual} is equivalent to a maximization over a subset $\cQ$ of $\reals^s$. Specifically, $\cQ = \{q\in \reals^s~|~ \exists \bfpi \in \Pi \mbox{ such that } q_\xi = \bfpi(\xi), \xi\in \Xi\}$. In the case of superquantiles, \eqref{eqn:superriskenvelope} gives $\Pi$ and this yields
\begin{equation}\label{eqn:cQ}
\cQ = \Big\{q\in \reals^s~\Big|~ 0 \leq q_\xi \leq \frac{1}{1-\alpha}, \xi\in \Xi, ~\nsum_{\xi\in\Xi} p_\xi q_\xi = 1\Big\}.
\end{equation}
By setting $\bar p_\xi = p_\xi q_\xi$, we reconcile the present development with the formula \eqref{eqn:maxsuper}. 

Regardless of the circumstances, the risk-minimization problem \eqref{eqn:riskmindual} is equivalently stated as a minsup problem for which there are many algorithmic approaches. These include subgradient-type methods which apply, especially, if $f(\xi,\cdot\,)$ is convex for all $\xi\in \Xi$ (see, e.g., \cite[Section 2.I]{primer}), the outer approximation algorithm (see, e.g., \cite[Section 6.C]{primer} and the discussion in Subsection \ref{eqn:algosuper}), and ``primal-dual'' algorithms (see, e.g., \cite{KouriSurowiec.22}). While these approaches provide important stepping stones, there are bound to be implementation challenges depending on the specific risk measure under consideration. In the absence of a finite distribution for $\bfxi$, the first equality in \eqref{eqn:riskmindual} still holds but it becomes harder to compute a maximizer, which may not even exist, causing additional challenges; see our discussion in Subsection \ref{subsec:riskidentifiers}.\\

For recent derivations of risk envelope formulas, we refer to \cite{AngSunYao.18,SunYangYaoZhang.20}, which focus on regular measures of deviation and the use of set operations, and \cite{RockafellarRoyset.18}, which concentrates on mixed risk measures and second-order superquantiles; see also \cite{Shapiro.13,FollmerSchied.16}. Dual expressions for measures of risk and related quantities remain a crucial stepping stone in several contexts. The papers \cite{GrechukZabarankin.17,ErnstPichlerSprungk.22} utilize them in sensitivity analysis and \cite{KouriSurowiec.18} in derivation of optimality condition; see also \cite{RockafellarUryasevZabarankin.06b}. Dual expressions are key in developing algorithms, including in the computationally challenging areas of PDE-constrained optimization \cite{KouriSurowiec.16,KouriSurowiec.20,KouriSurowiec.22} and statistical learning
\cite{CuriLevyJegelkaKrause.20,LevyCarmonDuchiSidford.20,LaguelMalickHarchaoui.20,LaguelPillutlaMalickHarchaoui.21}.

\subsection{Connections with Distributionally Robust Optimization}\label{subsec:distributionallyrobust}

It has long been recognized that a probability distribution adopted in an application is just a {\em model} of the uncertainty associated with unsettled parameters; see \cite{Dupacova.66} for the first recorded work in the area of stochastic programming but the perspective can be traced back to the early studies of games as summarized in \cite{vonNeumannMorgenstern.44}. This leads to {\em distributionally robust optimization problems} of the form
\begin{equation}\label{eqn:dro}
\nnmin_{x\in X} \,\sup_{P\in \cP} \int f(\xi,x)\, dP(\xi),
\end{equation}
where $f:\reals^m\times \reals^n\to \reals$ represents a quantity of interest, $X$ is a constraint set, and $\cP$ is an {\em ambiguity set} of candidate distributions $P$ on $\reals^m$. (Here, $P$ is a probability distribution on $\reals^m$, which clashes slightly with the earlier use of $P$ for a cumulative distribution function.) Distributionally robust optimization problems may be perceived to ``optimally'' address distributional ambiguity \cite{VanparysEsfahaniKuhn.21}. Key developments include \cite{WiesemannKuhnSim.14,Shapiro.17,ZhenKuhnWiesemann.21}; see also the recent reviews  \cite{PostekDenhertogMelenberg.16,RahimianMehrotra.22,LinFangGao.22}. A main focus in this area is the development of tractable reformulations of the minsup problem \eqref{eqn:dro} as a minimization problem. We see next that such tractable reformulations are immediately available for certain $\cP$.

The following fact is well known from convex analysis; see, e.g., \cite{RockafellarUryasev.13}.

\begin{proposition}{\rm (conjugacy under monotonicity).}\label{pConjMono}
For a lsc convex functional $\cF:\cL^2\to (-\infty,\infty]$ with $\cF\not\equiv \infty$, one has that $\cF$ is monotone if and only if $\bfpi\geq 0$ for all $\bfpi \in \dom \cF^*$.
\end{proposition}

The proposition together with Example \ref{eRiskEnv} confirm that superquantile risk and worst-case risk measures are monotone, while mean-plus-standard-deviation risk measures are not.

Combining Propositions \ref{pPropRiskEnvelope} and \ref{pConjMono}, we conclude that a risk envelope $\Pi$ of a positively homogeneous, monotone, regular measure of risk contains exclusively nonnegative random variables with expectation one. This has the consequence that we can interpret the expression $\Ex[\bfxi\bfpi] = \int \bfxi(\omega)\bfpi(\omega) \, d\bbP(\omega)$ as the expectation of $\bfxi$ under a different probability distribution $\bar \bbP$. From this perspective, $\bfpi$ is the density (Radon-Nikodym derivative) $d\bar \bbP/d\bbP$. Consequently, a monotone, positively homogeneous, regular measure of risk can be expressed by \eqref{eqn:poshomRisk} using a risk envelope $\Pi$ or, equivalently, as the supremum of $\int \bfxi(\omega) d\bar\bbP(\omega)$ over an ambiguity set with probability distributions $\bar\bbP$, which brings us back to \eqref{eqn:dro}. We refer to \cite{RuszczynskiShapiro.06} for a formal treatment and to \cite[Section 6]{RockafellarUryasev.13} for a brief introduction.

This insight implies that {\em every} monotone, positively homogeneous, regular measure of risk inherently addresses ambiguity leading to some distributionally robust problem akin to  \eqref{eqn:dro}. Conversely, a  distributionally robust problem \eqref{eqn:dro} corresponds to a risk-minimization problem provided that the densities defined by $\cP$ form a set $\Pi$ satisfying the properties in the second half of Proposition \ref{pPropRiskEnvelope}, which implicitly requires that the distributions in $\cP$ are absolutely continuous with respect to the underlying ``base'' distribution. The risk-minimization problem is, in turn, equivalently to a regret minimization problem; cf. Subsection \ref{subsec:riskmin}. The risk-minimization and the regret-minimization problems might be tractable alternatives to the original distributionally robust problem \eqref{eqn:dro}.

In summary, the introductory discussion in Subsection \ref{subsec:equivalence} about Mr.~Averse, a risk-averse decision maker relying on superquantiles, and Ms.~Ambiguity, which is risk-neutral but uncertain about the underlying probability distribution, extends to all monotone, positively homogeneous, regular measures of risk. These risk measures can therefore be used to model {\em both} risk-averseness and distributional ambiguity. Regardless of the initial motivation, the duality between risk-averseness and distributional ambiguity allows us to switch between the two perspectives, adopting the one that is computationally attractive or affords other advantages.

Further connections appear in \cite{GotohKimLim.18}, which shows that an ambiguity set centered at an empirical distribution makes the distributionally robust problem \eqref{eqn:dro} in some sense close to that of minimizing a {\em mean-plus-standard-deviation} risk measure; see also \cite{GotohKimLim.21} and references therein. Recent reviews of distributionally robust models in engineering and portfolio optimization include \cite{KapteynWillcoxPhilpott.19,PflugPohl.18}. Adversarial learning \cite{MadryMakelovSchmidtTsiprasVladu.18} is also of the form \eqref{eqn:dro} with a specific ambiguity set that only shifts the support.

\subsection{Risk Identifiers and Subgradients}\label{subsec:riskidentifiers}

While solving distributionally robust optimization problems such as \eqref{eqn:dro} or utilizing dual formulas in risk minimization, it might be important to identify which probability distribution $P\in \cP$ or which random variable $\bfpi\in \Pi$, if any, attains the maximum. These issues are also intimately tied to subgradients of functionals. Convex analysis provides key insights; see, e.g., \cite{Rockafellar.74}.

\begin{definition}{\rm (subgradients of functional).}\label{subgradFucnt}
For a convex functional $\cF:\cL^2\to [-\infty,\infty]$ and a point $\bfxi_0\in \cL^2$ at which $\cF$ is finite, $\bfpi\in \cL^2$ is a {\em subgradient} of $\cF$ at $\bfxi_0$ when
\[
\cF(\bfxi)\geq \cF(\bfxi_0) + \Ex\big[\bfpi(\bfxi - \bfxi_0) \big] ~~~~\forall \bfxi,\bfxi_0\in \cL^2.
\]
The set of all subgradients of $\cF$ at $\bfxi_0$ is denoted by $\partial \cF(\bfxi_0)$.
\end{definition}

The {\em interior} of $\cC\subset \cL^2$, denoted by $\nt \cC$, consists of every $\bfxi\in \cC$ for which there exists $\rho>0$  such that $\{\bfxi_0~|~ \|\bfxi_0-\bfxi\|_{\cL^2} \leq \rho\} \subset \cC$. In particular, $\nt(\dom \cF) = \cL^2$ when $\cF$ is real-valued. From \cite[Proposition 3.1]{RuszczynskiShapiro.06}, we know that a monotone convex functional $\cF:\cL^2\to (-\infty,\infty]$ has at least one subgradient at every point in $\nt (\dom \cF)$. The following classical fact from convex analysis provides a means to calculate subgradients; see, for example, \cite[Proposition 8.36]{primer} for a short proof based on the Fenchel-Moreau theorem \ref{pConjugateFcnal}. The claim about nonemptiness follows from \cite[Corollary 8B and Theorem 11]{Rockafellar.74}.

\begin{proposition}{\rm (subgradients from conjugates).}\label{lConjugateInverseFunc} For a lsc convex functional $\cF:\cL^2\to (-\infty,\infty]$ and a point $\bfxi$ at which $\cF$ is finite, one has
\[
\partial \cF(\bfxi) = \nargmax\big\{ \Ex[\bfxi\bfpi] - \cF^*(\bfpi)~\big|~ \bfpi\in \cL^2 \big\}.
\]
This subset of $\cL^2$ is nonempty provided that $\bfxi \in  \nt (\dom \cF)$.
\end{proposition}

For a positively homogeneous regular measure of risk, we say that $\hat\bfpi \in \nargmax \{ \Ex[\bfxi\bfpi] ~|~ \bfpi \in \Pi\}$ is a {\em risk identifier} of $\cR$ at $\bfxi$. Thus, the maximum in \eqref{eqn:poshomRisk} is indeed attained as long as $\bfxi \in \nt (\dom \cR)$ and $\cR$ is positively homogeneous and regular. The term risk identifier was coined in \cite{RockafellarUryasevZabarankin.06b} but the quantity was known much earlier simply as a subgradient of $\cR$. A more general treatment of existence of risk identifiers appears in \cite{Kouri.17}, with a particular focus on distributionally robust optimization. Expressions for risk identifiers in the case of mixed measures of risk appear in \cite{RockafellarRoyset.18}.

We now have tools to compute subgradients of functions involving risk measures and related functionals as already recognized by \cite{RuszczynskiShapiro.06,RockafellarUryasevZabarankin.06b}. For concreteness, we limit the focus to compositions with linear functions; see, e.g., \cite[Section 4]{RockafellarRoyset.18} for more general cases.

Given a positively homogeneous regular measure of risk $\cR$ and $\bfeta,\bfxi_1, \dots, \bfxi_n\in \cL^2$, we consider the problem
\begin{equation}\label{eqn:Riskminproblem}
\nnmin_{x\in X\subset\reals^n} ~\phi(x) = \cR\big( \bfeta-\langle \bfxi, x\rangle\big),
\end{equation}
where $\bfxi = (\bfxi_1, \dots, \bfxi_n)$. For example,  $\bfeta-\langle \bfxi, x\rangle$ might represent the (random) shortfall of ``production'' $\langle \bfxi, x\rangle$ relative to ``demand'' $\bfeta$. We also consider the generalized  regression problem \eqref{eqn:regProblem}, which in view of Proposition \ref{prop:errdev} effectively reduces to 
\begin{equation}\label{eqn:devminproblem}
\nnmin_{c\in C\subset\reals^n} ~\psi(c) = \cD\big(\bfeta - \langle c, \bfxi\rangle\big),
\end{equation}
where $\cD$ is the regular measure of deviation that pairs (in the sense of Theorem \ref{thm:errdev}) with the regular measure of error of interest in \eqref{eqn:regProblem}. There are accessible subgradient formulas for these objective functions, which allow us to bring in subgradient methods, cutting-plane methods, and related algorithms. The following proposition is taken from \cite[Proposition 8.38]{primer}; see \cite{RuszczynskiShapiro.06,RockafellarUryasevZabarankin.06b} for similar formulas.

\begin{proposition}{\rm (subgradients in risk and deviation minimization).}\label{pCompRuleRisk} For a real-valued, positively homogeneous, and regular measure of risk $\cR$, with risk envelope $\Pi$, the measure of deviation $\cD$ paired with $\cR$ in Theorem \ref{thm:quadrangle}, and random variables $\bfeta, \bfxi_1, \dots, \bfxi_n \in \cL^2$, with $\bfxi = (\bfxi_1, \dots, \bfxi_n)$, consider $\phi:\reals^n\to \reals$ in \eqref{eqn:Riskminproblem} and $\psi:\reals^n\to \reals$ in \eqref{eqn:devminproblem}. Then,
\begin{align*}
\partial \phi(x) & = \Big\{ - \Ex[\bfxi \hat\bfpi]  ~\Big|~ \hat\bfpi \in \nargmax\big\{\Ex\big[\big(\bfeta -\langle \bfxi,x\rangle\big) \bfpi\big] ~\big|~ \bfpi \in \Pi \big\} \Big\}~~~~\forall x\in \reals^n\\
\partial \psi(c) & = \Big\{ \Ex[\bfxi] - \Ex[\bfxi \hat\bfpi]  ~\Big|~ \hat\bfpi \in \nargmax\big\{\Ex\big[\big(\bfeta -\langle c, \bfxi\rangle\big) \bfpi\big] ~\big|~ \bfpi \in \Pi \big\} \Big\}~~~~\forall c\in \reals^n.
\end{align*}
\end{proposition}

In the case of superquantiles, we obtain the following example. 

\begin{example}\label{eSuperquantileSubgrad}{\rm (subgradients for superquantile functions).} For $\alpha\in (0,1)$ and finitely distributed random variables $\bfeta, \bfxi_1, \dots, \bfxi_n$, with values $\{\eta^i, \xi_1^i$, $\dots$, $\xi_n^i$, $i=1, \dots, s\}$ and probabilities $\{p_i>0, i=1, \dots, s\}$, consider the function $\phi:\reals^n\to \reals$ given by
\[
\phi(x) = \srsk_\alpha\big(\bfeta - \langle \bfxi,x\rangle \big),
\]
where $\bfxi = (\bfxi_1, \dots, \bfxi_n)$. Let $\xi^i = (\xi^i_1, \dots, \xi^i_n)$. A subgradient of $\phi$ at $x$ is then given by
\[
-\nsum_{i=1}^s p_i \xi^i \hat q_i, ~~~\mbox{ where }~ \hat q_i = \begin{cases}
\frac{1}{1-\alpha} &\mbox{ if }   i\in \bbI_+(x)\\
\frac{\sigma-\alpha}{(1-\alpha)\tau} & \mbox{ if } i\in \bbI_0(x)\\
0 &\mbox{ otherwise},
\end{cases}
\]
with $\sigma = 1 - \nsum_{i \in \bbI_+(x)} p_i$, $\tau = \nsum_{i \in \bbI_0(x)} p_i$,  and, denoting the $\alpha$-quantile of $\bfeta - \langle \bfxi, x\rangle$ by $Q(\alpha)$, also 
\[
\bbI_+(x) = \big\{i ~\big|~\eta^i - \langle \xi^i, x\rangle > Q(\alpha)\big\} ~~\mbox{ and } ~~\bbI_0(x) = \big\{i ~\big|~\eta^i - \langle \xi^i, x\rangle = Q(\alpha)\big\}.
\]
\end{example}
\state Detail. Using $\cQ$ in \eqref{eqn:cQ} and the arguments around that equation, we seek
\[
\hat q \in\nargmax_{q \in \cQ} \nsum_{i=1}^s p_i \big(\eta^i - \langle \xi^i, x\rangle\big) q_i,
\]
which is essentially available explicitly. The asserted subgradient follows by Proposition \ref{pCompRuleRisk}. 

The computational work required to obtain a subgradient is of order $O(s\ln s)$ as it essentially requires sorting the numbers $\{\eta^i- \langle \xi^i, x\rangle$, $i=1, \dots, s\}$. We also obtain the alternative formula: $\srsk_\alpha(\bfeta - \langle \bfxi, x\rangle) = \nsum_{i=1}^s p_i (\eta^i - \langle \xi^i, x\rangle) \hat q_i$.\eop

\section{Additional Examples}\label{sec:addexamples}

There is a vast landscape of risk measures beyond the instances listed in Example \ref{eMeasRegRisk}; \cite{RockafellarUryasev.13} and \cite[Chapter 6]{ShapiroDentchevaRuszczynski.21} review many others. 
This section discusses a few additional risk measures and provides recipes for constructing even more. We also touch on other measures of regret, deviation, and error. 

We refer to \cite{Pichler.17} for a discussion of how to compare risk measures quantitatively, and this might help in the process of choosing a suitable one. 

\begin{example}{\rm (mean-plus-upper-semideviation)}. For $\lambda\in (0,\infty)$, the measures of risk, regret, deviation, and error given by  
\begin{align*} 
\cR(\bfxi) &= \Ex[\bfxi] + \lambda\sqrt{\Ex\big[\max\{0,\bfxi - \Ex[\bfxi]\}^2\big]}    ~~~~~  &&\cD(\bfxi)  = \lambda\sqrt{\Ex\big[\max\{0,\bfxi - \Ex[\bfxi]\}^2\big]}\\
\cV(\bfxi) &= \lambda\sqrt{\Ex\big[\max\{0,\bfxi - \Ex[\bfxi]\}^2\big]} + \max\big\{0, 2\Ex[\bfxi]\big\} ~~~                              &&\cE(\bfxi) = \lambda\sqrt{\Ex\big[\max\{0,\bfxi - \Ex[\bfxi]\}^2\big]} + \big|\Ex[\bfxi]\big|
\end{align*}
are regular and form a risk quadrangle in the sense of Figure \ref{fig:quad}. Moreover, the functionals are all positively homogenous. The corresponding statistic is $\cS(\bfxi) = \Ex[\bfxi]$. The risk envelope of $\cR$, furnishing the alternative expression \eqref{eqn:poshomRisk}, is given by
\[
\Pi = \big\{\bfone + \lambda\bfpi - \lambda\Ex[\bfpi] ~\big|~\bfpi(\omega) \geq 0 ~\forall\omega, ~\Ex[\bfpi^2] \leq 1\big\}. 
\]
We refer to $\cR$ as a mean-plus-upper-semideviation risk measure.
\end{example}
\state Detail. Most of these properties can be found in \cite[Section 6.3]{ShapiroDentchevaRuszczynski.21} or follow immediately from the various definitions. The construction of measures of regret and error is motivated by the approach taken in the proofs of Theorems 8.9 and 8.21 in \cite{primer}. Extensions beyond square-integrable random variables, to $p$-integrable random variables, are discussed in \cite[Section 6.3]{ShapiroDentchevaRuszczynski.21}. That reference also covers when the resulting risk measures are monotone. 

Markowitz already in \cite{Markowitz.59} studied semideviations. Later efforts include \cite{OgryczakRuszczynski.99,OgryczakRuszczynski.01,OgryczakRuszczynski.02}, with a focus on connections with stochastic dominance, and the algorithmic developments in \cite{KalogeriasPowell.18,KalogeriasPowell.22}.\eop

\begin{example}{\rm (entropic risk)}.\label{eEntropicrisk} The measures of risk, regret, deviation, and error given by 
\begin{align*} 
\cR(\bfxi) &= \ln \Ex\big[\exp(\bfxi)\big]   ~~~~~  &&\cD(\bfxi) = \ln \Ex\big[\exp(\bfxi - \Ex[\bfxi])\big]\\
\cV(\bfxi) &= \Ex\big[\exp(\bfxi)-1\big] ~~~~&&                              \cE(\bfxi)  = \Ex\big[ \exp(\bfxi)-\bfxi-1\big]
\end{align*}
are regular  and form a risk quadrangle in the sense of Figure \ref{fig:quad}. The corresponding statistic is $\cS(\bfxi) = \ln \Ex[\exp(\bfxi)]$. While $\cR$ is monotone, it is not positively homogeneous and therefore lacks an expression of the form \eqref{eqn:poshomRisk}.  Still, the Fenchel-Moreau theorem \ref{pConjugateFcnal} applies and yields the alternative expression 
\[
\cR(\bfxi) = \nsup\big\{ \Ex[\bfxi\bfpi] - \Ex[\bfpi \ln \bfpi]  ~\big|~ \bfpi(\omega) \geq 0 ~\forall \omega, ~\Ex[\bfpi] = 1 \big\},
\] 
where $0 \ln 0$ is defined as 0. We refer to $\cR(\bfxi)$ as entropic risk. 
\end{example}
\state Detail. These facts can be deduced from \cite[Section 6.3]{ShapiroDentchevaRuszczynski.21} and \cite{RockafellarUryasev.13}. We note that $\cR$ may not be real-valued unless, for example, the underlying probability space is finite; see \cite[Section 6.3]{ShapiroDentchevaRuszczynski.21} for a detailed discussion. The measure of regret $\cV$ stems from a ``normalized'' exponential disutility function and thus has a long history within expected utility theory. We refer to \cite{LaevenStadje.13} for extensions.\eop

Examples \ref{eMeasRegRisk}(d), \ref{eMeasRegError}(d), \ref{eMeasRegDev}(d), and the supporting discussion show how we can construct new regular measures of risk, regret, error, and deviation by combining existing ones. Subsection \ref{subsec:mixed} relies on the same principle, but with the narrower focus on using superquantiles as the ``base'' measures of risk from which others emerge. Simple scaling may also be used. If $\cR_0, \cV_0, \cD_0, \cE_0$ are regular measures of risk, regret, deviation, and error forming a risk quadrangle in the sense of Figure \ref{fig:quad}, with statistic $\cS_0$, and $\lambda \in (0,\infty)$, then  the functionals $\cR, \cV, \cD ,\cE$ given by
\begin{align*} 
\cR(\bfxi) &= \lambda \cR_0(\lambda^{-1}\bfxi)    ~~&& \cD(\bfxi)  = \lambda\cD_0(\lambda^{-1}\bfxi)\\
\cV(\bfxi) &= \lambda \cV_0(\lambda^{-1}\bfxi)    ~~&&     \cE(\bfxi)  = \lambda \cE_0(\lambda^{-1}\bfxi)
\end{align*}
are regular and form a risk quadrangle, with statistic $\cS(\bfxi)  = \lambda\cS_0(\lambda^{-1}\bfxi)$. Here, monotonicity is preserved by this scaling process. Naturally, the scaling is only interesting when a functional is {\em not} positively homogeneous such as in Example \ref{eEntropicrisk}. 

Adding a weighted functional to the expected value is also meaningful; see, for example, \cite{FrohlichWilliamson.22} for usage in machine learning. Again starting with $\cR_0, \cV_0, \cD_0, \cE_0$, all regular and forming a risk quadrangle, and $\lambda \in (0,\infty)$, we obtain that the functionals given by 
\begin{align*} 
\cR(\bfxi) &= (1-\lambda) \Ex[\bfxi] + \lambda \cR_0(\bfxi)    ~~&&  \cD(\bfxi) = \lambda\cD_0(\bfxi)\\
\cV(\bfxi) &= (1 - \lambda)\Ex[\bfxi] + \lambda \cV_0(\bfxi)   ~~&&  \cE(\bfxi) = \lambda \cE_0(\bfxi)
\end{align*}
are regular and form a risk quadrangle, with statistic $\cS(\bfxi) = \cS_0(\bfxi)$. Positive homogeneity is preserved by this process regardless of $\lambda \in (0,\infty)$ and the same holds for monotonicity in the case of $\cD$ and $\cE$. Monotonicity of $\cR$ and $\cV$ is preserved when $\lambda \leq 1$. The properties of scaling and expectation-mixing follow straightforwardly from the various definitions; see \cite{RockafellarUryasev.13} for a recording of these facts. 

For additional means to construct new functionals from existing ones, we refer to \cite{FrohlichWilliamson.22,FrohlichWilliamson.22b} and \cite{LiuMaoWangWei.22}, which combines risk measures using epi-sums (inf-convolutions) as motivated by risk sharing. The paper \cite{WangWeiWillmot.20}  studies signed Choquet integrals and their convexity and \cite{LiuSchiedWang.21} utilizes distributional transforms.

\begin{example}{\rm (entropic value-at-risk)}. For $\alpha \in [0,1)$, the measure of risk given by 
\[
\cR_\alpha(\bfxi) = \inf_{\gamma>0} \gamma^{-1} \ln\bigg(\frac{\Ex[\exp(\gamma\bfxi)]}{1-\alpha} \bigg)
\]
is convex, positively homogeneous, monotone, and also satisfies the constancy property. 
The risk envelope of $\cR_\alpha$ furnishing the alternative expression \eqref{eqn:poshomRisk} is given by
\[
\Pi = \big\{ \bfpi  ~\big|~ \Ex[\bfpi \ln \bfpi] \leq -\ln(1-\alpha), ~\Ex[\bfpi] = 1, ~\bfpi(\omega) \geq 0 ~\forall \omega\big\},
\]
where again $0 \ln 0$ is defined as 0. Thus, this risk measure yields the worst-case expected value with respect to probability distributions that has Kullback-Leibler divergence from $\bbP$ of at most $-\ln(1-\alpha)$. In particular, $\cR_0(\bfxi) = \Ex[\bfxi]$ and $\lim_{\alpha\upto 1} \cR_\alpha(\bfxi) = \sup \bfxi$. The $\alpha$-superquantile of $\bfxi$ never exceeds $\cR_\alpha(\bfxi)$. We refer to $\cR_\alpha(\bfxi)$ as the entropic value-at-risk (at level $\alpha$). 
\end{example}
\state Detail. These properties are given in \cite{Ahmadi.11,Ahmadi.12,Ahmadi.12b}, which pioneered entropic value-at-risk. With the involvement of the moment-generating function $\gamma\mapsto \Ex[\exp(\gamma\bfxi)]$ in the definition of $\cR_\alpha$, it is natural that these claims are restricted to random variables for which this function is finite; see \cite{Ahmadi.12,Ahmadi.12b} for details. An advantage of entropic value-at-risk is the ease by which it can be computed under various independence assumptions. 

For far-reaching extensions leveraging more general divergences in the construction of risk envelopes, we refer to \cite{PichlerSchlotter.20b}, where, for example, Kusuoka representations appear as well. Efforts in this direction already emerged in \cite{Ahmadi.12}.\eop 

For additional measures of risk, we refer to \cite{CastagnoliCattelanMaccheroniTebaldiWang.22}, which considers a broad class of star-shaped risk measures that extend beyond the convex ones. Connections between risk measures and scoring rules emerge from \cite{SmithBickel.22}. The paper \cite{LeqiPrasadRavikumar.19} discusses risk measures motivated by cumulative prospect theory and \cite{LiuWang.21} examines tail-focused generalizations of quantiles and superquantiles. Further constructions involving divergences can be found in \cite{FollmerWeber.15,DommelPichler.20}.

\section{Computational Tools}\label{sec:comptool}

The number and capabilities of computational tools for risk minimization are rapidly expanding. Standard packages for optimization address many problems, especially after reformulations to common formats such as linear and nonlinear programs; see for example Subsections \ref{subsec:superinopt} and \ref{subsec:superregression}. There are also specialized packages that address risk minimization directly.

{\em Portfolio Safeguard} \cite{Aorda.22} allows users to specify measures of risk, regret, deviation, and error using a high-level language. It also has built-in custom algorithms that address large-scale problems and those involving nonsmoothness, for example caused by the formula \eqref{eqn:superquantRU}. It also has a tailored format for superquantile regression (referred to as CVaR regression in the documentation). However, being a commercial product, the details about some of the algorithms are unavailable. Portfolio Safeguard has Matlab and R interfaces.

The Python package {\em SPQR} \cite{spqr.22} implements first-order algorithms for superquantile minimization from \cite{LaguelMalickHarchaoui.20}; see also the discussion in Subsection \ref{eqn:algosuper}. Another Python package {\em SQwash} \cite{sqwash.22} provides connections with PyTorch to solve superquantile problems involving neural networks. It also implements dual smoothing algorithms from \cite{LaguelPillutlaMalickHarchaoui.21,PillutlaLaguelMalickHarchaoui.22} as surveyed in Subsection \ref{eqn:algosuper}. These two packages are especially tailored to statistical estimation problems.

There are several packages motivated by financial applications, an area where superquantiles are known by the alternative names conditional value-at-risk, tail value-at-risk, and expected shortfall. To solve risk-minimization problems, the Python package {\em Riskfolio-Lib} \cite{riskfolio.22} leverages the optimization package {\em CvxPy} and closely integrates with pandas data structures. {\em AzaPy} \cite{azapy.22} is another Python package that handles mixed superquantile risk measures.
{\em Cvar-Portfolio} \cite{cvarport.22} is a portfolio optimization package in Python that addresses high-dimensional problems.

Beyond Python, we find packages in R \cite{port.22}, C++ \cite{tailrisk.22}, and Julia \cite{emp.22,cvar.22,riskmeasures.22}, with the latter reference furnishing implementation of an exceptionally wide range of risk measures. The Rapid Optimization Library \cite{rol.22}, developed in C++, focuses on large-scale engineering applications and has capabilities to address problems with uncertainty including risk minimization.

\section{Extensions and Challenges}\label{sec:extensions}

Since decision making under uncertainty arise in nearly all human activity, this survey is unable to review all aspects in detail. We have ignored the fact that some optimization problems involve infinite-dimensional $x$ such as in PDE-constrained optimization; see, e.g., \cite{KouriShapiro.18}. While this introduces numerous challenges, it does {\em not} change the definition of risk measures or their desirable properties. We still apply a risk measure to $f(\bfxi,x)$, which is a random variable for any (infinite-dimensional) $x$ as long as $f(\,\cdot\,,x)$ is measurable, just as before. In fact, many articles cited in the survey deal with such general cases. This section mentions two other domains: measures of reliability and multi-stage problems. We end with a discussion of research opportunities.

\subsection{Measures of Reliability}\label{subsec:reliability}

The risk of a random variable has the same unit as the random variable itself; it is a number that represents (conservatively) the unknown future value. Definition \ref{dRiskmeasure} reflects this through its treatment of constant random variables. In reliability engineering and operations research, it is common to consider other ways of quantifying uncertainty and ensuring safety. Here, we propose terminology for these alternative ways as a counterpart to Definition \ref{dRiskmeasure}.

\begin{definition}{\rm (reliability measure)}.
A {\it measure of reliability} $\cP$ assigns to a random variable $\bfxi$ a number $\cP(\bfxi) \in [0,1]$ as a quantification of its reliability, with this number being either 0 or 1 if $\bfxi$ is a constant random variable.
\end{definition}

A measure of reliability is a functional on a space of random variables, with (hopefully) desirable properties just as measures of risk, regret, deviation, and error. Since $\cP(\bfxi) \in [0,1]$, it can be interpreted as a probability with the following prominent examples.\\

\noindent{\bf Probability of exceedance.} For $\tau\in\reals$, the choice $\cP_\tau(\bfxi) = \prob\{\bfxi > \tau\}$ utilizes the {\em probability of exceedance} of the threshold $\tau$. This measure of reliability is in one-to-one correspondence with the quantile risk measure $\cR_\alpha(\bfxi) = Q(\alpha)$, where $Q(\alpha)$ is the $\alpha$-quantile of $\bfxi$, $\alpha \in (0,1)$. Specifically, $\cR_\alpha(\bfxi) \leq \tau$ if and only if $\cP_\tau(\bfxi) \leq 1-\alpha$. If the threshold $\tau = 0$, then $\cP_\tau(\bfxi)$ is called the {\em failure probability} of $\bfxi$; see, e.g., \cite{RockafellarRoyset.15b}. Probabilities of exceedance give rise to chance constraints as seen in \cite[Examples 2.57 and 3.20]{primer}.\\

\noindent{\bf Buffered probability of exceedance.} For an integrable random variable $\bfxi$ with superquantiles $\bar Q(\alpha)$, $\alpha \in [0,1]$, the {\em buffered probability} of exceeding $\tau\in\reals$ is
\[
\bprob\{\bfxi> \tau\} = \begin{cases}
0 & \mbox{ if } \tau \geq \bar Q(1)\\
1-\bar Q^{-1}(\tau) & \mbox{ if } \bar Q(0) < \tau < \bar Q(1)\\
1 & \mbox{ otherwise.}
\end{cases}
\]
Here, $\bar Q^{-1}(\tau)$ is the solution to the equation $\tau = \bar Q(\alpha)$, which is unique when $\bar Q(0) < \tau < \bar Q(1)$ because $\bar Q$ is a continuous function on $(0,1)$ and also increasing on $(0,\,\prob\{\bfxi<\bar Q(1)\})$; see \cite[Theorem 2]{RockafellarRoyset.14} and \cite[Proposition A.1]{MafusalovUryasev.18}. For each $\tau$, this defines the measure of reliability $\cP_\tau(\bfxi) = \bprob\{\bfxi> \tau\}$ first studied in \cite{RockafellarRoyset.10} for the case $\tau =0$ under the name {\em buffered failure probability} and later extended to address natural hazards \cite{DavisUryasev.16}, classification problems \cite{NortonMafusalovUryasev.17,NortonUryasev.19}, and cash-flow matching \cite{ShangKuzmenkoUryasev.18}; see \cite{MafusalovUryasev.18,MafusalovShaprioUryasev.18} for theoretical foundations and \cite{RockafellarUryasev.20} for a decomposition method to solve resulting optimization problems. In addition to  \cite{RockafellarRoyset.10}, other application papers in reliability engineering include \cite{MinguezCastilloLara.13,ZrazhevskyGolodnikovUryasevZrazhevsky.20,ByunDeoliveiraRoyset.22,ByunRoyset.22}. Sensitivity analysis of buffered probabilities appears in \cite{ZhangUryasevGuan.19,RoysetByun.21}. The paper \cite{Kouri.19} defines higher-order buffered probabilities. The measure of reliability $\cP_\tau$ is in one-to-one correspondence with the $\alpha$-superquantile measure of risk: For integrable $\bfxi$, $\alpha \in (0,1]$, and $\tau\in \reals$, one has $\srsk_\alpha(\bfxi) \leq \tau$ if and only if $\cP_\tau(\bfxi) \leq 1-\alpha$.

\subsection{Dynamic and Multi-Stage Optimization}\label{subsec:multi}

Many real-world problems involve intricate and gradual revelation of the ``true'' values of uncertain parameters, usually intertwined with decisions, which even might be continuously adjusted. Modeling of such situations lead to optimal control problems, Markov decision processes, and multi-stage stochastic programming. Decisions now need to account for the uncertainty associated with stochastic processes as well as the opportunities for later recourse actions. This complicates notions of risk-averseness and how to model it in a meaningful and computationally tractable manner.

Again, the initial efforts were motivated by financial applications \cite{ArtznerDelbaenEberHeathKu.02} and dealt with cadlag processes \cite{CheriditoDelbaenKupper.04,CheriditoDelbaenKupper.05}; see also \cite{CheriditoDelbaenKupper.06}. The papers \cite{Riedel.04,RuszczynskiShapiro.06b} sever connections with the ``static'' case and develop dynamic and conditional risk measures from an axiomatic point of view. Around the same time, \cite{EichhornRomisch.05} defines polyhedral risk measures for multi-stage problems; see also \cite{GuiguesRomisch.12}. Informally, a conditional risk measure at time $t$ assigns to a sequence of random variables (representing future ``costs'') a random variable modeling the amount we would be willing to pay at time $t$ to avoid facing the future costs. A collection of such conditional risk measures, one for each time period, amounts to a dynamic measure of risk. Time consistency of a dynamic measure of risk is then the property that it will deem, at the present time, a future sequence of random costs $\bfxi^1$ no worse than another sequence $\bfxi^2$ whenever it deems $\bfxi^1$ at least as good as $\bfxi^2$ at a future time $t$ and $\bfxi^1$ is identical to $\bfxi^2$ between now and $t$. We refer to the tutorials \cite{Ruszczynski.13,Shapiro.21}, the monograph \cite{PflugPichler.14}, and the recent papers \cite{DommelPichler.21,FanRuszczynski.22} for many more details. Connections with distributionally robust optimization appear in \cite{PichlerShapiro.21}. The paper \cite{IancuPetrikSubramanian.15} discusses the difference between quantifying the risk using an end state as compared to quantifying it using a ``composite'' risk incurred at each stage in a decision process. Risk measures in the context of Markov decision processes are pioneered in \cite{Ruszczynski.10,Ruszczynski.14}, which furnish new thinking about time consistency, dynamic programming equations, and a value iteration method. Later efforts include \cite{ChowTamarMannorPavone.15,ShaoGuptaHaskell.22,XiaGlynn.22,BieleckiCialencoRuszczynski.22}.

The papers \cite{GarcFern.15,ChowGhavamzadehJansonPavone.18} address reinforcement learning, while  \cite{MillerYang.17,PichlerSchlotter.22,WangChapman.22} study optimal control and \cite{ShapiroTekayaDacostaSoares.13} examines dual dynamic programming from a risk-averse perspective.

\subsection{Challenges and Open Problems}

Despite 25 years of rapid development, the area of risk measures and related concepts remains relevant with many opportunities and open problems.

A wealth of opportunities emerge in connection with {\em autonomous systems} such as driverless cars and delivery vehicles; see the recent review \cite{WangChapman.22}. The resulting problems are often ``dynamic'' with the need for a combination of learning and decision making, but with conservativeness and reliability as a main focus.

As seen above, there have been many efforts to utilize superquantiles in {\em statistical learning} but fewer to leverage other risk measures; see \cite{FrohlichWilliamson.22,FrohlichWilliamson.22b} for notable exceptions. In view of Section \ref{sec:errordev}, it appears that measures of error and deviation should be central to this area, in fact more than measures of risk. Besides \cite{Yang.21}, there seem to be few attempts at utilizing measures of risk for {\em adversarial training} and similar approaches to``robust'' learning. Since adversarial training corresponds to the choice of a worst-case risk measure, it only represents one of many possibilities for promoting conservativeness and robustness.

There are surprisingly few studies of {\em Bayesian formulations} of risk-averse optimization problems under uncertainty. A recent effort is \cite{WuZhuZhou.18}, which applies a risk measure to a posterior distribution, but this can be extended in several directions and presents alternatives to distributionally robust formulations. Situations with {\em decision-dependent probability measures} and {\em partially observable Markov decision processes} are likewise fertile ground; see \cite{DentchevaRuszczynski.20} for emerging ideas.

As discussed in Subsection \ref{subsec:distributionallyrobust}, risk measures may arise in response to ambiguity about the ``true'' probability distribution. A source of such ambiguity could be a quantity of interest $f(\bfxi,x)$ for which we know only the marginal distributions of the random vector $\bfxi = (\bfxi_1, \dots, \bfxi_m)$. The paper \cite{GhossoubHallSaunders.22} examines the situation in the context of distributionally robust optimization. It would be useful to develop {\em risk envelopes} in concrete situations with such ``missing information.'' 

A solution $(\gamma^\star, c^\star)$ of the generalized  regression problem \eqref{eqn:regProblem} produces a ``best'' approximation $\gamma^\star + \langle c^\star, \bfxi\rangle$ of the random variable $\bfeta$ in the sense of a regular measure of error $\cE$. Still, $\gamma^\star + \langle c^\star, \xi\rangle$ may not approximate well the conditional statistic $\cS(\bfeta_\xi)$ for all $\xi$. Here, $\cS$ is the statistic corresponding to $\cE$ in the sense of Definition \ref{dErrorStat} and $\bfeta_\xi$ is a random variable with the distribution of $\bfeta$ conditional on $\bfxi$ having the value $\xi$. Theorem 5.1 of \cite{RockafellarRoyset.15} addresses such {\em statistic tracking} for cases with additive noise and certain risk measures. However, we would like to understand this issue in broader settings.

The construction of {\em surrogates} (cf. Example \ref{eSurrogate}) and supporting {\em designs of experiments} may need to be carried out conservatively. As in \cite{BonfiglioRoyset.19}, we often seek surrogates that predict response quantities from more easily available parameters and models. With our orientation toward response quantities such as ``cost'' and ``damage,'' where high values are undesirable, we tend to be more concerned about underestimating the response than overestimating. Design of experiments may need to be carry out with these factors in mind, while also considering risk-averseness relative to the outcomes of the experiments. For example, with few experiments left, the common expected improvement criterion for selecting the next design point might be replaced with a risk-averse alternative. We refer to \cite{KouriJakemanHuerta.22} for initial efforts in this direction, but numerous opportunities remain with the possibility of improving how complex physical systems are assessed and designed. Ideas from distributionally robust Bayesian optimization \cite{KirschnerBogunovicJegelkaKrause.20} may provide guidance for efforts in this direction.

Risk-based formulations also arise for {\em games and equilibrium problems} as exemplified by \cite{RalphSmeers.15,LunaSagastizabalSolodov.16,PangSenShanbhag.17,LianeasNikolovaStiermoses.19,FerrisPhilpott.22}. The consideration of risk in these settings is complicated by the two sources of uncertainty for an agent: lack of knowledge about the action of the other agents and inherent randomness in the environment. Moreover, each agent might have their own tolerance for risk, which raises the question whether there is a central planner, utilizing {\em some} measure of risk, that achieves a distributed solution. Since games and equilibrium problems tend to be large scale, nonsmooth, and nonconvex, there are also numerous computational challenges in this area.

While extending back to \cite{RockafellarRoyset.10}, {\em buffered probabilities} have not been studied in the context of surrogate models. The recent sensitivity results in \cite{RoysetByun.21} may provide a useful stepping stone toward surrogates for buffered probabilities of exceedance. While a chance constraint naturally extends from a single random variable not exceeding a threshold to a random vector taking values in a multi-dimensional set, a buffered probability constraint is presently limited to a single random variable. With the introduction of {\em measures of reliability}, we hope to spur the development of many alternatives to failure probabilities and buffered probabilities as a parallel to the numerous measures of risk.\\

\noindent {\bf Acknowledgement.} The author is thankful for valuable input and encouragement from Amir Ahmadi-Javid, Harbir Antil, Christian Fr\"{o}hlich, Drew Kouri, Boris Kramer, and Ruodu Wang and for the insights provided by two reviewers. This work is supported in part by ONR (Mathematical and Resource Optimization), ONR (Science of Autonomy), and AFOSR (Mathematical Optimization). 

\addcontentsline{toc}{section}{References}

\bibliographystyle{plain}
\bibliography{refs}

\end{document}